\pdfoutput=1

\documentclass[12pt]{amsart}
\usepackage[margin=1in]{geometry} 

\usepackage{amsmath,amssymb,amsthm,bm,colonequals,graphicx,mathrsfs,mathtools,microtype,stmaryrd}
\usepackage[shortlabels]{enumitem}
\usepackage[colorinlistoftodos]{todonotes}

\numberwithin{equation}{section}

\theoremstyle{plain}

\newtheorem{theorem}{Theorem}[section] 

\newtheorem{lemma}[theorem]{Lemma} 

\newtheorem{proposition}[theorem]{Proposition} 

\newtheorem{proposition-definition}[theorem]{Proposition-Definition} 

\newtheorem{corollary}[theorem]{Corollary} 

\theoremstyle{definition}

\newtheorem{definition}[theorem]{Definition} 

\newtheorem{example}[theorem]{Example} 

\newtheorem{question}[theorem]{Question} 

\theoremstyle{remark}

\newtheorem{remark}[theorem]{Remark} 

\newtheorem{observation}[theorem]{Observation} 

\newcommand{\Aff}{\mathbb{A}} 

\newcommand{\CC}{\mathbb{C}}
\newcommand{\EE}{\mathbb{E}}
\newcommand{\FF}{\mathbb{F}}

\newcommand{\PP}{\mathbb{P}}
\newcommand{\QQ}{\mathbb{Q}}
\newcommand{\RR}{\mathbb{R}}

\newcommand{\ZZ}{\mathbb{Z}}

\newcommand{\im}{\operatorname{im}}

\newcommand{\res}{\overline}

\newcommand{\abs}[1]{\lvert #1 \rvert}
\newcommand{\card}[1]{\lvert #1 \rvert}
\newcommand{\norm}[1]{\lVert #1 \rVert}
\newcommand{\floor}[1]{\lfloor #1 \rfloor}

\newcommand{\eps}{\epsilon}

\newcommand{\Supp}{\operatorname{Supp}}

\newcommand{\disc}{\operatorname{disc}}
\newcommand{\Aut}{\operatorname{Aut}}

\newcommand{\Pic}{\operatorname{Pic}}

\newcommand{\map}{\operatorname}
\newcommand{\mscr}{\mathscr}
\newcommand{\mcal}{\mathcal}
\newcommand{\mf}{\mathfrak}

\newcommand{\ol}{\overline}

\newcommand{\inmat}[1]{\left[\begin{smallmatrix} #1 \end{smallmatrix}\right]}

\newcommand{\texpdf}{\texorpdfstring}

\newcommand{\defeq}{\colonequals}

\newcommand{\maps}{\colon}

\newcommand{\belongs}{\subseteq}

\newcommand{\set}[1]{\{#1\}}

\newcommand{\grad}{\nabla}

\usepackage[pagebackref=true]{hyperref}
\usepackage[alphabetic,initials]{amsrefs}

\begin{document}

\title{Special cubic zeros and the dual variety}


\date{}
\author{Victor Y. Wang}
\address{Fine Hall, 304 Washington Road, Princeton, NJ 08540, USA}
\address{Courant Institute, 251 Mercer Street, New York, NY 10012, USA}
\address{IST Austria, Am Campus 1, 3400 Klosterneuburg, Austria}
\email{vywang@alum.mit.edu}

\subjclass{Primary 11D45; Secondary 11D25, 11G25, 11P55, 14B05}
\keywords{Cubic form, circle method, special subvarieties, biases, dual variety}

\begin{abstract}
Let $F$ be a diagonal cubic form over $\mathbb{Z}$ in $6$ variables.
From the dual variety in the delta method of Duke--Friedlander--Iwaniec and Heath-Brown,
we unconditionally extract a weighted count of certain special integral zeros of $F$ in regions of diameter $X \to \infty$.
Heath-Brown did the same in $4$ variables, but our analysis differs and captures some novel features.
We also put forth an axiomatic framework for more general $F$.
\end{abstract}

\maketitle

\setcounter{tocdepth}{3}

\section{Introduction}
\label{SEC:motivating-background-and-intro-to-structured-part-detection}

Let $F\in \ZZ[\bm{x}]=\ZZ[x_1,\dots,x_m]$ be a cubic form in $m\geq 4$ variables.
Let $V$ be the hypersurface $F=0$ in $\PP^{m-1}_\QQ$.
Let $w\maps \RR^m\to \RR$ be a smooth, compactly supported function.
Let $\Supp{w}$ be the closure of the set $\{\bm{x}\in \RR^m: w(\bm{x})\neq 0\}$ in $\RR^m$.
Assume that $V$ is smooth and (for convenience) that $\bm{0}\notin \map{Supp} w$.
For reals $X\geq 1$, let
\begin{equation}
\label{EQN:define-N(X)}
N_{F,w}(X) \defeq \sum_{\bm{x}\in \ZZ^m:\, F(\bm{x})=0} w(\bm{x}/X).
\end{equation}
Let $\Upsilon$ be the set of vector spaces $L\belongs \QQ^m$ over $\QQ$ of dimension $\floor{m/2}$ with $F\vert_L = 0$.
The ``Hardy--Littlewood model'' for $N_{F,w}(X)$ can fail via $\Upsilon$ if $m\leq 6$.\footnote{For singular $V$, a similar failure can occur for larger $m$; see \cite{brudern2019instance} for a nice example with $m=8$.}
But based on \cite{hooley1986some}*{Conjecture~2}, \cite{franke1989rational}, \cite{vaughan1995certain}*{Appendix}, \cite{peyre1995hauteurs}, et al., one may conjecture
\begin{equation}
\label{EQN:Hooley--Manin--Peyre-smooth-cubic-case}
N_{F,w}(X)
- \sum_{\bm{x}\in \ZZ^m:\, \bm{x}\in \bigcup_{L\in \Upsilon} L} w(\bm{x}/X)
= (c_{F,w}+o_{X\to\infty}(1))\cdot X^{m-3}(\log{X})^{r-1+\bm{1}_{m=4}}
\end{equation}
for a certain predicted constant $c_{F,w}\in \RR$ and integer $r\geq 1$, where $r=1$ if $m\geq 5$.
Here $r$ is the rank of the Picard group of $V$.

For further discussion and references on \eqref{EQN:Hooley--Manin--Peyre-smooth-cubic-case}, see \cite{bombieri2009problems}*{\S2} or \cite{wang2022thesis}*{\S6.5}.
The natural analog of \eqref{EQN:Hooley--Manin--Peyre-smooth-cubic-case} for rational points is closely related to \eqref{EQN:Hooley--Manin--Peyre-smooth-cubic-case}, and is in fact equivalent when $m\ge 5$, but the precise sieve-theoretic relationship is delicate when $m=4$;
cf.~\cite{heath1996new}*{paragraph preceding Corollary~2} and \cite{heath1996new}*{Theorem~8 versus Corollary~2} on quadrics.

In the present paper, we seek to organically detect, within the delta method of \cites{duke1993bounds,heath1996new}, the locus $\bigcup_{L\in \Upsilon} L$ featured in \eqref{EQN:Hooley--Manin--Peyre-smooth-cubic-case}.
In the delta method, as in \cite{kloosterman1926representation}'s circle method,
one averages over numerators to a given denominator (modulus), and often then invokes Poisson summation.
Set $Y\defeq X^{3/2}$ as on \cite{heath1998circle}*{p.~676}; then
\begin{equation}
\label{EQN:delta-method}
(1+O_A(Y^{-A}))\cdot N_{F,w}(X)
= Y^{-2}
\sum_{n\geq1}
\sum_{\bm{c}\in\ZZ^m}
n^{-m}S_{\bm{c}}(n)I_{\bm{c}}(n)
\end{equation}
(for all $A>0$) by \cite{heath1996new}*{Theorem~2, (1.2), up to easy manipulations from \S3},
where
\begin{equation}
\label{EQN:define-I_c(n)}
I_{\bm{c}}(n)
\defeq \int_{\bm{x}\in \RR^m}d\bm{x}\,w(\bm{x}/X)
h(n/Y,F(\bm{x})/Y^2)
e(-\bm{c}\cdot\bm{x}/n)
\end{equation}
for a certain fixed smooth function $h\maps (0,\infty) \times \RR \to \RR$
usually left in the background
(see \cite{heath1998circle}*{(2.3)} for the precise definition of $h$),
and where
\begin{equation}
\label{EQN:define-S_c(n)}
S_{\bm{c}}(n) \defeq
\sum_{1\leq a\leq n:\, \gcd(a,n)=1}\,
\sum_{1\leq \bm{x}\leq n} e_n(aF(\bm{x}) + \bm{c}\cdot\bm{x}).
\end{equation}
Here $\bm{c}\cdot \bm{x}\defeq c_1x_1+\dots+c_mx_m$, and in $S_{\bm{c}}(n)$ the variable $\bm{x}$ runs over $1\leq x_1,\dots,x_m\leq n$.
See Proposition~\ref{PROP:standard-delta-method-n,c-cutoffs} for the basic analytic properties
of $I_{\bm{c}}(n)$.

Given $\bm{c}\in\ZZ^m$,
the sum $S_{\bm{c}}(n)$ is multiplicative in $n$,
and roughly governed by
the solutions to $F(\bm{x}) = \bm{c}\cdot \bm{x} = 0$ over $\ZZ/n\ZZ$.
Let $V_{\bm{c}}$ be the intersection $F(\bm{x}) = \bm{c}\cdot \bm{x} = 0$ in $\PP^{m-1}_\QQ$.
There is a classical \emph{discriminant polynomial} $F^\vee\in \ZZ[\bm{c}]$ associated to the family $\bm{c}\mapsto V_{\bm{c}}$.
The equation $F^\vee(\bm{c})=0$ cuts out the projective \emph{dual variety} $V^\vee$ of $V$.
See \S\ref{SEC:alg-geom-background} for details on $F^\vee$, $V^\vee$.


It is natural to analyze the right-hand side of \eqref{EQN:delta-method} separately over $F^\vee(\bm{c})\neq 0$ and $F^\vee(\bm{c})=0$.
In another paper (\cite{wang2023_HLH_vs_RMT}), we conditionally address $F^\vee(\bm{c})\neq 0$ for certain $F$, $w$.
In the present (self-contained) paper, we focus on $F^\vee(\bm{c})=0$.
From this locus in \eqref{EQN:delta-method},
our main theorem unconditionally isolates $\bigcup_{L\in \Upsilon} L$, for certain $F$:

\begin{theorem}
\label{THM:contribution-from-generically-singular-c's}
Suppose $F$ is diagonal and $m\in \set{4,6}$.
For $Y=X^{3/2}$ as above, we have
\begin{equation}
\label{EQN:main-theorem-singular-c's}
Y^{-2}\sum_{\bm{c}\in \ZZ^m:\, F^\vee(\bm{c}) = 0,\; \bm{c}\ne \bm{0}}\,
\sum_{n\geq 1} n^{-m}S_{\bm{c}}(n)I_{\bm{c}}(n)
= O_{F,w,\eps}(X^{m/2-1/4+\eps})
+ \sum_{L\in \Upsilon} \sum_{\bm{x}\in L\cap \ZZ^m} w(\bm{x}/X).
\end{equation}
\end{theorem}

For $m=4$,
Theorem~\ref{THM:contribution-from-generically-singular-c's}
follows from \cite{heath1998circle}*{Lemmas 7.2 and 8.1},
with a better error of $O_{F,w,\eps}(X^{3/2+\eps})$.
As we will explain in \S\ref{SEC:full-proof-outline},
the case $m=6$ seems to present new difficulties, perhaps most easily resolved using our new, more geometric, strategy.
But \cite{heath1998circle} does bound the left-hand side of \eqref{EQN:main-theorem-singular-c's} by $O_{F,w,\eps}(X^{3+\eps})$ for $m=6$,
giving us a good foundation to build on.
In view of \cite{glas2022question}, one might also hope to adapt our results to function fields.\footnote{Very recently, our work has largely been extended to function fields \cite{BGW2024forthcoming}.}

Our main results assume $F$ is diagonal, but we do also discuss non-singular $F$ to a nontrivial extent.
In particular, we will isolate explicit technical ingredients
(listed in Remark~\ref{RMK:non-diagonal-F})
that---if true more generally---would
allow one to generalize Theorem~\ref{THM:contribution-from-generically-singular-c's}.

If $2\mid m$ (as in Theorem~\ref{THM:contribution-from-generically-singular-c's}), the set $\Upsilon$ is known to be finite
for general reasons (recalled in \S\ref{SEC:maximal-linear-under-duality} below).
Furthermore, we can relate \eqref{EQN:main-theorem-singular-c's} to \eqref{EQN:Hooley--Manin--Peyre-smooth-cubic-case} by inclusion-exclusion:
\begin{equation*}
\sum_{L\in \Upsilon} \sum_{\bm{x}\in L\cap \ZZ^m} w(\bm{x}/X)
= O_{F,w}(X^{m/2-1})
+ \sum_{\bm{x}\in \ZZ^m:\, \bm{x}\in \bigcup_{L\in \Upsilon} L} w(\bm{x}/X).
\end{equation*}
For each $L\in \Upsilon$, if we define $\sigma_{\infty,L^\perp,w}$ as in \eqref{EQN:real-density-linear-spaces}, we also have (for all $A>0$)
\begin{equation}
\label{EQN:intro-Lambda-sum-re-interpretation-formula}
\sum_{\bm{x}\in L\cap \ZZ^m} w(\bm{x}/X)
= \sigma_{\infty,L^\perp,w} X^{m/2} + O_{L,w,A}(X^{-A}).
\end{equation}

If $F$ is diagonal,
$\Upsilon$ has an explicit classical description:
see Proposition~\ref{PROP:characterize-diagonal-L/Q}.
In general, linear subvarieties of cubics
are closely tied to the \emph{$h$-invariant} introduced by \cites{davenport1962exponential,davenport1964non}.
(See \cite{dietmann2017h}*{Lemma~1.1} for a more general relationship.)
In fact,
the present \S\ref{SEC:proving-bias-in-linear-S_c's} relies on
a convenient choice of
``$h$-decompositions of $F$''
corresponding to the elements of $\Upsilon$.

We believe that with enough work, the power saving $1/4$ in Theorem~\ref{THM:contribution-from-generically-singular-c's} could be improved.
It would be very interesting to go past $1/2$ for $m=4$.

Theorem~\ref{THM:contribution-from-generically-singular-c's} has the following corollary, which we need for subsequent work (\cite{wang2023_HLH_vs_RMT}).

\begin{corollary}
\label{COR:Disc-vanishing-sum-with-explicit-linear-density}
Let $m=6$.
Define $\sigma_{\infty,F,w}$, $\mf{S}_F$ as in \S\ref{SEC:c=0-singular-series}.
Then in the setting of Theorem~\ref{THM:contribution-from-generically-singular-c's},
\begin{equation*}
Y^{-2}\sum_{\bm{c}\in \ZZ^m:\, F^\vee(\bm{c}) = 0}\,
\sum_{n\geq 1} n^{-m}S_{\bm{c}}(n)I_{\bm{c}}(n)
= O_{F,w,\eps}(X^{2.75+\eps})
+ \sigma_{\infty,F,w} \mathfrak{S}_F X^3
+ \sum_{L\in \Upsilon} \sigma_{\infty,L^\perp,w} X^3.
\end{equation*}
\end{corollary}

\begin{proof}
Combine \eqref{EQN:main-theorem-singular-c's} (on $\bm{c}\ne \bm{0}$) with
\eqref{INEQ:isolate-singular-series-from-center} (on $\bm{c}=\bm{0}$).
Note that $m/2=m-3=3$.
\end{proof}

Theorem~\ref{THM:contribution-from-generically-singular-c's} and Corollary~\ref{COR:Disc-vanishing-sum-with-explicit-linear-density} are completely unconditional.
When $m=6$ and $F$ is diagonal,
these results let us reformulate
the conjecture \eqref{EQN:Hooley--Manin--Peyre-smooth-cubic-case}
as a statement purely about \emph{cancellation} in
the sum $\sum_{\bm{c}\in \ZZ^m:\, F^\vee(\bm{c})\neq 0}
\sum_{n\geq 1} n^{-m} S_{\bm{c}}(n) I_{\bm{c}}(n)$.
A similar reformulation might be possible much more generally.
But at least when $m=4$,
subtleties in the constant $c_{F,w}$ in \eqref{EQN:Hooley--Manin--Peyre-smooth-cubic-case}
would likely demand
a recipe
beyond ``restriction to $F^\vee=0$'' in \eqref{EQN:delta-method}.

Before proceeding, we make some convenient definitions.
Let
\begin{equation}
\label{DEFN:normalize-S_c(n)}
    S^\natural_{\bm{c}}(n)\defeq n^{-(1+m)/2}S_{\bm{c}}(n);
\end{equation}
then in particular,
$n^{-m} S_{\bm{c}}(n) I_{\bm{c}}(n)
= n^{(1-m)/2} S^\natural_{\bm{c}}(n) I_{\bm{c}}(n)$.
The following is also convenient:

\begin{definition}
\label{DEFN:primitive-sublattices}
Given a vector space $L\belongs \QQ^m$ over $\QQ$,
let $L^\perp$ be the orthogonal complement of $L$
with respect to $\bm{c}\cdot \bm{x}$.
Then let $\Lambda\defeq L\cap \ZZ^m$ and $\Lambda^\perp\defeq L^\perp \cap \ZZ^m$.
\end{definition}



We now sketch the proof of Theorem~\ref{THM:contribution-from-generically-singular-c's}.
The proof starts generally,
observing that $F^\vee\vert_{L^\perp}=0$ for all $L\in \Upsilon$
(see Proposition~\ref{PROP:dual-linear-subvariety}).
Conversely, at least for diagonal $F$,
most $\bm{c}$'s on the left-hand side of \eqref{EQN:main-theorem-singular-c's}
are in fact \emph{linear}, in the sense of the following definition:

\begin{definition}
\label{DEFN:orthogonally-linear-c's}
Call a solution $\bm{c}\in\ZZ^m$ to $F^\vee(\bm{c})=0$ \emph{linear}
if $\bm{c}\in \bigcup_{L\in \Upsilon} \Lambda^\perp$.
\end{definition}

We expect ``typical'' linear $\bm{c}$'s to be the simplest.
Proposition~\ref{PROP:dual-linear-subvariety}(2), which establishes a vanishing baseline for the jets $j^\bullet{F^\vee}$ over $\bigcup_{L\in \Upsilon} L^\perp$,
thus inspires the following definition:

\begin{definition}
\label{DEFN:unsurprising-Disc-zero-locus}
Call $F^\vee$ \emph{unsurprising} if uniformly over reals $C\geq 1$,
the box $[-C,C]^m$ contains at most $O_\eps(C^{m/2-1+\eps})$ points in the union of the following two sets:
\begin{equation}
\label{EQN:define-exceptional-sets-E_1,E_2}
{\textstyle
\mcal{E}_1\defeq \set{\bm{c}\in \ZZ^m: F^\vee(\bm{c})=0}\setminus \bigcup_{L\in \Upsilon} \Lambda^\perp,
\quad \mcal{E}_2\defeq \set{\bm{c}\in \bigcup_{L\in \Upsilon} \Lambda^\perp: j^{2^{m/2-1}}{F^\vee}(\bm{c})=\bm{0}}.
}
\end{equation}
\end{definition}

We prove (in \S\ref{SEC:full-proof-outline}) that when $F$ is diagonal,
$F^\vee$ is unsurprising,
so that if $\Upsilon\neq \emptyset$, then almost all solutions to $F^\vee=0$ are linear with nonzero $2^{m/2-1}$-jet.\footnote{An amusing corollary is that the Gauss map $\gamma\maps V\to V^\vee$ introduced in \S\ref{SEC:alg-geom-background} can be far from surjective on $\QQ$-points, since $\gamma^{-1}([\bm{c}])\cap V(\QQ)$ may be empty for a typical nonzero $\bm{c}\in \bigcup_{L\in \Upsilon} L^\perp$.}
A weaker version of Definition~\ref{DEFN:unsurprising-Disc-zero-locus}, with $C^{m/2-\delta}$ in place of $C^{m/2-1+\eps}$, would suffice for qualitative purposes.


For the ``least degenerate'' linear $\bm{c}$'s,
Lemma~\ref{LEM:stating-bias-for-generic-linear-c's} isolates an explicit positive bias
\begin{equation}
\label{EQN:illustrate-bias-philosophy}
S^\natural_{\bm{c}}(p^l)
= [A_{p^l}(\bm{c}) + O(p^{-l/2})]
\cdot (1-p^{-1})\cdot p^{l/2}
\end{equation}
for most primes $p$,
with $A_p(\bm{c})=1$ and $A_{p^l}(\bm{c})\ll 1$.
The resulting reduction in arithmetic complexity of $S_{\bm{c}}(n)$ lets us
dramatically simplify $\sum_{\bm{c}\in\Lambda^\perp}S_{\bm{c}}(n)I_{\bm{c}}(n)$
by Poisson summation over various individual residue classes $\bm{c}\equiv \bm{b}\bmod{n_0\Lambda^\perp}$ with $n_0\ll n/X$ dividing $n$.
When $m=6$, we must also carefully separate $\bm{c}=\bm{0}$ from $\bm{c}\neq\bm{0}$;
we use Lemma~\ref{LEM:sum-S_0(n)-trivially} (decay of the singular series over large moduli).
Eventually, $\sigma_{\infty,L^\perp,w} X^{m/2}$ appears, yielding \eqref{EQN:main-theorem-singular-c's}.


As the positivity of $S_{\bm{c}}(n)$ for linear $\bm{c}$'s might suggest,
we do not need cancellation over $n$ in \eqref{EQN:main-theorem-singular-c's}.
The deepest result we use on $L$-functions (when $m=6$) is
the (purely local) Weil bound for hyperelliptic curves of genus $\leq 2$.
However, it might be possible to improve the error term in Theorem~\ref{THM:contribution-from-generically-singular-c's} using deeper results on $L$-functions.

The following remark proposes axioms that would let us go beyond diagonal $F$.

\begin{remark}
\label{RMK:non-diagonal-F}
In general, \eqref{EQN:main-theorem-singular-c's} holds, provided $m\in \set{4,6,8}$ and (1)--(4) hold:
\begin{enumerate}
    \item $F^\vee$ is unsurprising (in the sense of Definition~\ref{DEFN:unsurprising-Disc-zero-locus}).
    
    \item $\Supp{w}\belongs \set{\bm{x}\in \RR^m: \det(\map{Hess}{F}(\bm{x}))\neq 0}$.
    
    \item There exists a homogeneous polynomial $W\in \ZZ[\bm{c}]$, satisfying
    \begin{equation}
    \label{COND:linear-locus-degeneracy-contains-W=0}
    {\textstyle \set{\bm{c}\in \bigcup_{L\in \Upsilon} \Lambda^\perp: W(\bm{c})=0}}
    \belongs \set{\bm{c}\in \ZZ^m: j^{2^{m/2-1}}{F^\vee}(\bm{c})=\bm{0}},
    \end{equation}
    such that uniformly over integers $C,n\geq 1$, we have
    \begin{align}
    \sum_{\bm{c}\in [-C,C]^m:\, F^\vee(\bm{c})=0,\; W(\bm{c})\neq 0} n^{-1}\abs{S^\natural_{\bm{c}}(n)}^2
    &\ll_\eps (Cn)^\eps (C^{m/2} + n^{m/6}),
    \label{COND:second-moment-axiom-over-nonzero-W} \\
    \sum_{\bm{c}\in [-C,C]^m:\, F^\vee(\bm{c})=0,\; W(\bm{c})=0}
    n^{-1/2}\abs{S^\natural_{\bm{c}}(n)}
    &\ll_\eps (Cn)^\eps (C^{(m-1)/2} n^{1/6} + n^{m/6}).
    \label{COND:first-moment-axiom-over-W=0}
    \end{align}
    
    \item In Lemma~\ref{LEM:stating-bias-for-generic-linear-c's},
    the formula for $S_{\bm{c}}(p)$ and upper bound for $S_{\bm{c}}(p^{\geq 2})$ remain true,
    provided $p$ exceeds some constant depending on $m$ and $F$.
\end{enumerate}
(There could be alternatives to (3), but (3) would be convenient.
See \S\ref{SEC:worst-case-Disc-0-sparse-bound} and \S\ref{SEC:full-proof-outline} for details.)
\end{remark}

In Remark~\ref{RMK:non-diagonal-F},
we expect (4) to be the most tractable in general out of (1), (3), and (4).
In fact, \cites{wang2023dichotomous,BGW2024forthcoming} already make progress on (4).
Also, \cite{salberger2023counting} could be helpful for (1), at least if $m=4$.
Furthermore, the case where $m=8$ and $F$ is diagonal might be fully accessible, and might provide insight on secondary terms in the circle method.
This would be similar in spirit to \cite{getz2018secondary},
where secondary terms were obtained for quadrics in an even number of variables.
The most mysterious axiom might be (3), but it does at least hold for diagonal $F$, and could plausibly follow from some undiscovered geometric stratification.

\begin{proposition}
Suppose $F$ is diagonal, and $m\geq 4$ is arbitrary.
Then \ref{RMK:non-diagonal-F}(3) holds.
\end{proposition}

\begin{proof}
[Proof sketch]
Let $W\defeq c_1\cdots c_m$.
One can prove \eqref{COND:second-moment-axiom-over-nonzero-W} and \eqref{COND:first-moment-axiom-over-W=0} using \eqref{INEQ:crude-pointwise-bound-on-diagonal-S_c(n)}.
(See \cite{wang2023_large_sieve_diagonal_cubic_forms}*{\S7};
\eqref{COND:second-moment-axiom-over-nonzero-W} holds with $C^{m/2}$ in place of $C^{m/2} + n^{m/6}$.)
Also, Observation~\ref{OBS:higher-order-diagonal-Disc-vanishing-critera} implies \eqref{COND:linear-locus-degeneracy-contains-W=0}.
\end{proof}

\begin{remark}
It would be interesting to extend our analysis to $m=5$,
even just for diagonal $F$.
See \cite{bombieri2009problems}*{\S3} for a discussion of
the potentially infinite family of lines on $V\belongs \PP^4$ if $m=5$.
Can one see these lines via $V^\vee$
(cf.~Proposition~\ref{PROP:dual-linear-subvariety})?
It is conceivable that for $m=5$, one might have to look past the locus $F^\vee=0$ in \eqref{EQN:delta-method}, and perhaps include all $\bm{c}$ for which $V_{\bm{c}}$ has Picard rank $\geq 2$.
The ``bias'' philosophy behind \eqref{EQN:illustrate-bias-philosophy} may shed some light here.
\end{remark}

We now outline the rest of the paper.

\S\ref{SEC:alg-geom-background} collects some classical algebraic geometry used
as a black box in parts of \S\ref{SEC:maximal-linear-under-duality}.

\S\ref{SEC:maximal-linear-under-duality} suggests a geometric backbone in \eqref{EQN:delta-method} for
the eventual harmonic detection of $\Upsilon$.

\S\ref{SEC:reverse-sums-by-duality} proves some identities
to be used in \S\ref{SEC:full-proof-outline} to obtain the expected ``main terms'' in \eqref{EQN:main-theorem-singular-c's}.

\S\ref{SEC:worst-case-Disc-0-sparse-bound} provides general upper bounds to be used at several points in \S\ref{SEC:full-proof-outline}.

\S\ref{SEC:c=0-singular-series} analyzes the $\bm{c}=\bm{0}$ contribution in \eqref{EQN:delta-method},
while also proving Lemma~\ref{LEM:sum-S_0(n)-trivially}.

\S\ref{SEC:proving-bias-in-linear-S_c's} proves some new asymptotic formulas for $S_{\bm{c}}(p^l)$ (see \eqref{EQN:illustrate-bias-philosophy} and Lemma~\ref{LEM:stating-bias-for-generic-linear-c's}).

\S\ref{SEC:full-proof-outline} uses Lemmas~\ref{LEM:sum-S_0(n)-trivially} and~\ref{LEM:stating-bias-for-generic-linear-c's},
plus results from \S\S\ref{SEC:maximal-linear-under-duality}--\ref{SEC:worst-case-Disc-0-sparse-bound} and \cite{heath1998circle},
to prove Theorem~\ref{THM:contribution-from-generically-singular-c's}.


\subsection{Conventions}
\label{SUBSEC:conventions}


We let $\disc(F)$ denote the discriminant of $F$, so that $\disc(F)\ne 0$ if and only if $V$ is smooth.
We always assume $V$ is smooth, unless specified otherwise.

We define $\bm{1}_E$ to be $1$ if $E$ holds, and $0$ if $E$ does not hold.
We let $e(t)\defeq e^{2\pi it}$,
and $e_r(t)\defeq e(t/r)$.
Also, we use the following notation for integrals:
\begin{equation*}
\int_X dx\,f(x)\defeq \int_X f(x)\,dx,
\quad \int_{X\times Y} dx\,dy\,f(x,y)\defeq \int_X dx \left(\int_Y dy\,f(x,y)\right).
\end{equation*}

We let $\EE_{b\in S}[f(b)]\defeq \card{S}^{-1} \sum_{b\in S} f(b)$ (if $S$ is finite and nonempty).

If $F=F_1x_1^3+\dots+F_mx_m^3$, then we assume (after scaling $F^\vee$ if necessary) that
\begin{equation}
\label{INEQ:diagonal-Fdual-primitive-scaling-convention}
F^\vee \in \ZZ[\bm{c}]
\qquad \textnormal{and} \qquad
\gcd(\textnormal{coefficients of }F^\vee) = (6^m)!F_1\cdots F_m.
\end{equation}

For a vector space $U$, we let $\PP{U}$ denote the projectivization of $U$.

We let $\map{Hess}{F}=\map{Hess}{F}(\bm{x})$ denote the usual $m\times m$ Hessian matrix of $F$.

We let $\ZZ_{\geq 0}\defeq \set{n\in \ZZ: n\geq 0}$.
We let $\partial_x\defeq \partial/\partial x$;
we let $\partial_{\bm{c}}^{\bm{r}}\defeq \partial_{c_1}^{r_1}\cdots\partial_{c_m}^{r_m}$ (for $\bm{r}\in \ZZ_{\geq 0}^m$) and $\abs{\bm{b}}\defeq \sum_{i\in S} b_i$ (for $\bm{b}\in \ZZ^S$).
We let
\begin{equation}
\label{EQN:define-r-jet}
j^r{F^\vee}\defeq (\partial_{\bm{c}}^{\bm{r}}{F^\vee})_{\bm{r}\in \ZZ_{\geq 0}^m:\, \abs{\bm{r}}\leq r}
\end{equation}
denote the \emph{$r$-jet} of $F^\vee$,
recording all partial derivatives of $F^\vee$ of order $\leq r$.

We write $f\ll_S g$, or $g\gg_S f$, to mean ``$\abs{f}\leq Bg$ for some real $B = B(S) > 0$''.
We let $O_S(g)$ denote a quantity that is $\ll_S g$.
We write $f\asymp_S g$ if $f\ll_S g\ll_S f$.

When making estimates, we think of $m$, $F$, $w$ as fixed constants.

\section{Algebraic geometry background}
\label{SEC:alg-geom-background}

We do algebraic geometry in the language of schemes, with the symbols $=$ and $\belongs$ interpreted scheme-theoretically unless specified otherwise.
We call a reduced, locally closed subscheme of projective space over a field a \emph{variety}.
Ultimately, geometry lets us cut out subsets of $\ZZ^m$, and count points over finite rings, in rigorous ways (see e.g.~the key Lemma~\ref{LEM:stating-bias-for-generic-linear-c's}).
Most of our work can be interpreted classically, in terms of points over algebraically closed fields.

The rest of this section reviews the geometry of the gradient map
\begin{equation*}
    \grad{F} = (\partial F/\partial x_1,\dots,\partial F/\partial x_m).
\end{equation*}
For diagonal $F$,
a more explicit analysis is possible
(see \S\ref{SUBSEC:diagonal-dual-example}).

Since $V$ is smooth (and $3F = \bm{x}\cdot \grad{F}$), the gradient $\grad{F}$ defines a morphism $[\grad{F}]\maps \PP^{m-1}\to (\PP^{m-1})^\vee$.
(We will often identify the dual projective space $(\PP^{m-1})^\vee$ with $\PP^{m-1}$, using $c_1,\dots,c_m$ as projective coordinates.)
The map $[\grad{F}]$, known as the \emph{polar map} of $V$, is finite surjective of degree $2^{m-1}$, by dimension and intersection theory;
cf.~\cite{dolgachev2012classical}*{p.~29}.
Since $\PP^{m-1}$ is smooth, $[\grad{F}]$ must then be flat (by ``miracle flatness'').

\begin{definition}
\label{DEFN:define-Vdual-as-scheme-theoretic-image}
Let $V^\vee\belongs (\PP^{m-1}_\QQ)^\vee$ be the scheme-theoretic image of $V$ under $[\grad{F}]$.
\end{definition}

Since $V$ is a smooth projective hypersurface, $V^\vee$ is the \emph{dual variety} of $V$.
Upon restricting $[\grad{F}]$ to $V$,
we get the finite surjective \emph{Gauss map} $\gamma\maps V\to V^\vee$.

Here $V$ is geometrically irreducible, so $V^\vee = \im\gamma$ must be too.
Since $\gamma$ is finite, $V^\vee$ must therefore be a geometrically integral hypersurface, i.e.~the zero scheme of some \emph{absolutely irreducible form} $F^\vee\in \QQ[\bm{c}]$.
The definition of $V^\vee$ then implies that if $\bm{c}\in \CC^m\setminus \set{\bm{0}}$, then $F^\vee(\bm{c})=0$ if and only if the projective scheme $F(\bm{x})=\bm{c}\cdot\bm{x}=0$ over $\CC$ is singular.

(At least for diagonal $F$,
one can explicitly compute $F^\vee$;
see \eqref{EQN:diagonal-discriminant-factorization} in \S\ref{SUBSEC:diagonal-dual-example}.)

\begin{remark}
\label{RMK:equivalent-ways-to-define-dual-variety}
There are many equivalent ways to define $V^\vee$.
For example, $V^\vee$ is the locus of zeros $\bm{c}$ of a certain polynomial $\disc(F,\bm{c})$ in $c_1,\dots,c_m$ and the coefficients of $F$; see e.g.~\cite{wang2023dichotomous}*{Proposition~4.4}.
It should also be possible to interpret $V^\vee$ as a norm of $V$ under $[\grad{F}]$, in the sense of \cite{stacks-project}*{\href{https://stacks.math.columbia.edu/tag/0BD2}{Tag 0BD2}}.
\end{remark}

The polar map $[\grad{F^\vee}]\maps \PP^{m-1}\dashrightarrow \PP^{m-1}$ of $V^\vee$ is a rational map defined away from $\map{Sing}(V^\vee)$, the \emph{set of singular points} of $V^\vee$.
(Here $\map{Sing}(V^\vee)$ is a proper closed subset of $V^\vee$.
It is known that $V^\vee$ is singular, since $\deg{F}\geq 3$; see e.g.~\cite{wang2023dichotomous}*{Proposition~4.4}.)

The reflexivity theorem says that $(V^\vee)^\vee = V$.
The biduality theorem says that if $[\bm{x}]\in V$ and $[\bm{c}]\in V^\vee$ are smooth points,
then $[\grad{F}(\bm{x})] = [\bm{c}]$ if and only if $[\grad{F^\vee}(\bm{c})] = [\bm{x}]$.
(Both facts are on \cite{dolgachev2012classical}*{p.~30}.)
For us, $V$ is smooth,
so biduality implies that the polar maps $[\grad{F}]$, $[\grad{F^\vee}]$ restrict to
\emph{inverse morphisms} between $V\setminus[\grad{F}]^{-1}(\map{Sing}(V^\vee))$ and $V^\vee\setminus\map{Sing}(V^\vee)$.

It is known that $\deg{F^\vee} = 3\cdot 2^{m-2}$ \cite{dolgachev2012classical}*{p.~33, (1.47)}.

Since $[\grad{F}]\maps \PP^{m-1}\to (\PP^{m-1})^\vee$ is a finite surjective morphism of smooth varieties (and in particular, is generically \'{e}tale),
its ramification theory is well-behaved.
Following \cite{stacks-project}*{\href{https://stacks.math.columbia.edu/tag/0BWJ}{Tag 0BWJ}},
let $R_{[\grad{F}]}$ be the closed subscheme of $\PP^{m-1}$ cut out by the different ideal $\mf{D}_{[\grad{F}]}\belongs \mcal{O}_{\PP^{m-1}}$ of $[\grad{F}]$.
Following \cite{stacks-project}*{\href{https://stacks.math.columbia.edu/tag/0BW8}{Tag 0BW8} and \href{https://stacks.math.columbia.edu/tag/0BWA}{Tag 0BWA}}, let $B_{[\grad{F}]}$ be the norm of $R_{[\grad{F}]}$ (or equivalently, the discriminant of $[\grad{F}]$).
In our setting, $R_{[\grad{F}]}$ is an effective Cartier divisor in $\PP^{m-1}$, and $B_{[\grad{F}]}$ is thus an effective Cartier divisor in $(\PP^{m-1})^\vee$.





The points of $R_{[\grad{F}]}$ are precisely those $x\in \PP^{m-1}$ at which $[\grad{F}]$ is unramified.
Furthermore, we have a set-theoretic equality $B_{[\grad{F}]} = [\grad{F}](R_{[\grad{F}]})$.

\begin{definition}
\label{DEFN:polar-ramification-branch-divisors}
Call $R_{[\grad{F}]}$ and $B_{[\grad{F}]}$ the \emph{ramification divisor} and \emph{branch divisor} of $[\grad{F}]$, respectively.
Let $H_R\in \ZZ[\bm{x}]$ and $H_B\in \ZZ[\bm{c}]$ be homogeneous polynomials defining $R_{[\grad{F}]}$ in $\PP^{m-1}$ and $B_{[\grad{F}]}$ in $(\PP^{m-1})^\vee$, respectively.
\end{definition}

Since $\deg{F}\geq 3$, one can show that $R_{[\grad{F}]}$ and $B_{[\grad{F}]}$ are nonempty, and thus hypersurfaces.
In fact,
by \cite{dolgachev2012classical}*{p.~29, Proposition~1.2.1},
$R_{[\grad{F}]} = \map{hess}(V)$,
where $\map{hess}(V)$ denotes the subscheme $\det(\map{Hess}{F}(\bm{x}))=0$ of $\PP^{m-1}$.

\begin{proposition}
\label{PROP:universality-of-dual-variety-and-branch-divisor}
Say we let $F$ vary over the locus $\disc(F)\neq 0$.
Then one can choose $F^\vee$ and $H_B$ to be polynomials in $c_1,\dots,c_m$ and the coefficients of $F$.
\end{proposition}

\begin{proof}
This is possible by standard ``universal'' constructions compatible with our definitions of $V^\vee$ and $B_{[\grad{F}]}$.
For $B_{[\grad{F}]}$, one can appeal to \cite{stacks-project}*{\href{https://stacks.math.columbia.edu/tag/0BD2}{Tag 0BD2}}; norms respect base change (and thus vary nicely in families).
For $V^\vee$,
see e.g.~Remark~\ref{RMK:equivalent-ways-to-define-dual-variety}.
\end{proof}

It is known that $V\not\belongs \map{hess}(V)$ \cite{hooley1988nonary}*{Lemma~1}.
How about after applying $[\grad{F}]$?

\begin{question}
\label{QUES:dual-variety-vs-branch-locus}
Is it necessarily true that $V^\vee \not\belongs B_{[\grad{F}]}$?
\end{question}

Question~\ref{QUES:dual-variety-vs-branch-locus} comes up in Proposition~\ref{PROP:dual-linear-subvariety},
but we happen to be able to sidestep it there.

\section{Maximal linear subvarieties under duality}
\label{SEC:maximal-linear-under-duality}

Fix a smooth cubic $V/\QQ$ as in \S\ref{SEC:motivating-background-and-intro-to-structured-part-detection}.
The reader only interested in diagonal $F$ can skim forwards to \S\ref{SUBSEC:diagonal-dual-example},
which explicitly analyzes $\Upsilon$
through the lens of $F^\vee$.

\subsection{A preliminary general analysis}
\label{SUBSEC:prelim-general-analysis}


If $L\in \Upsilon$,
then differentiating $F$ along $L$ implies $L\perp\grad{F}(\bm{x})$
for all $\bm{x}\in L$.
So the restriction $\gamma\vert_{\PP{L}} = [\grad{F}]\vert_{\PP{L}}$ maps into $\PP{L^\perp}$.

Since $\deg{F}\geq 3$,
it is also known by \cite{debarre2003lines}*{Lemma~3} (or Starr \cite{starr2005fact_in_browning2006density}*{Appendix})
that if $m$ is \emph{even},
then $\Upsilon$ is finite.
We would like to understand $\Upsilon$ in terms of \eqref{EQN:delta-method}.
Proposition~\ref{PROP:dual-linear-subvariety} suggests one plausible starting route:
duality (i.e.~detecting $L^\perp$ through $F^\vee$).

\begin{proposition}
\label{PROP:dual-linear-subvariety}
Suppose $2\mid m\geq 4$,
and fix $L$ in $\Upsilon$.
\begin{enumerate}
    \item $\gamma\vert_{\PP{L}}$ is a finite flat surjective morphism $\PP{L}\to \PP{L^\perp}$ of degree $2^{m/2-1}\geq 2$.
    
    \item The jet $j^{2^{m/2-1}-1}{F^\vee}$, defined as in \eqref{EQN:define-r-jet}, vanishes over $\PP{L^\perp}$.
\end{enumerate}
\end{proposition}

For diagonal $F$,
we provide an explicit proof of Proposition~\ref{PROP:dual-linear-subvariety} in \S\ref{SUBSEC:diagonal-dual-example}.
For more general $F$,
we instead rely on \S\ref{SEC:alg-geom-background}, plus some progress on Question~\ref{QUES:fiber-degree-to-jet-vanishing}.
We would not be surprised if better technique allowed one to answer Question~\ref{QUES:fiber-degree-to-jet-vanishing} in general.


\begin{question}
\label{QUES:fiber-degree-to-jet-vanishing}
Let $y=[\bm{c}]\in (\PP^{m-1})^\vee$.
If the scheme-theoretic fiber $\PP^{m-1} \times_{[\grad{F}]} y$ of $[\grad{F}]$ over $y$ has degree $d$, is it necessarily true that $j^{d-1}{F^\vee}(\bm{c})=\bm{0}$?
\end{question}

We now begin the proof of Proposition~\ref{PROP:dual-linear-subvariety}.
Since $[\grad{F}]$ is finite,
the restriction $\gamma\vert_{\PP{L}}\maps \PP{L}\to \PP{L^\perp}$ is itself finite,
and thus surjective by dimension theory.
So $\PP{L^\perp}\belongs\im\gamma\vert_{\PP{L}}\belongs\im\gamma\belongs V^\vee$.

\begin{proof}
[Proof of (1)]

Since $\gamma\vert_{\PP{L}}$ is finite surjective and $\PP{L}$, $\PP{L^\perp}$ are smooth,
``miracle flatness'' implies flatness of $\gamma\vert_{\PP{L}}$.
Also,
$\gamma\vert_{\PP{L}}$ has degree $2^{m/2-1}$ (cf.~\cite{dolgachev2012classical}*{top of p.~29}),
since it is a morphism given by quadratic polynomials, between projective spaces of dimension $m/2-1$.
\end{proof}

Proposition~\ref{PROP:dual-linear-subvariety}(2) is inspired by
the factorization of $F^\vee$ over $\ol{\QQ}[\bm{c}^{1/2}]$ when $F$ is diagonal (see \eqref{EQN:diagonal-discriminant-factorization} in \S\ref{SUBSEC:diagonal-dual-example} below).
However,
giving a rigorous ``factorization'' of $F^\vee$ seems to require a bit of setup,
since the map $[\grad{F}]$ presumably need not be Galois in general.
Furthermore, our factorization will only be useful for (2) if
\begin{equation}
\label{COND:L^perp-not-contained-in-branch-locus}
    \PP{L^\perp} \not\belongs B_{[\grad{F}]}
    \quad(\textnormal{or equivalently, $H_B\vert_{L^\perp} \neq 0$}).
\end{equation}

\begin{question}
\label{QUES:can-L^perp-be-contained-in-branch-locus?}
Is $\PP{L^\perp} \belongs B_{[\grad{F}]}$ possible
(if $V$ is smooth and $L\in \Upsilon$)?
\end{question}

We do not know the answer to Question~\ref{QUES:can-L^perp-be-contained-in-branch-locus?}.
Fortunately, if we fix $m$ and $L$, then the set $\mscr{U}$ of all $m$-variable cubic forms $P/\QQ$ with $\disc(P)\neq 0$ and $P\vert_L=0$ is a dense open set in a copy of $\Aff^{N(m)}_\QQ$, where $N(m) = \binom{m+2}{3} - \binom{\frac{m}{2}+2}{3}$.
Furthermore,
Proposition~\ref{PROP:universality-of-dual-variety-and-branch-divisor} implies that
\ref{PROP:dual-linear-subvariety}(2) is a closed condition on $F\in \mscr{U}$,
and that \eqref{COND:L^perp-not-contained-in-branch-locus} is an open condition on $F\in \mscr{U}$.
Since \eqref{COND:L^perp-not-contained-in-branch-locus} holds when $L$ is $x_1+x_2=x_3+x_4=x_5+x_6=0$ (as we may assume, by a linear change of variables) and $F$ is $x_1^3+\dots+x_6^3$, it thus suffices to prove \ref{PROP:dual-linear-subvariety}(2) assuming \eqref{COND:L^perp-not-contained-in-branch-locus}.

Hence, for the rest of \S\ref{SUBSEC:prelim-general-analysis}, we assume \eqref{COND:L^perp-not-contained-in-branch-locus},
though \eqref{COND:L^perp-not-contained-in-branch-locus} will only come into play after some initial work.
Consider the hypersurface complements $S\defeq \PP^{m-1}\setminus B_{[\grad{F}]}$ and $X\defeq [\grad{F}]^{-1}S\belongs \PP^{m-1}$.
Then $[\grad{F}]\vert_X\maps X\to S$ is finite \'{e}tale of degree $2^{m-1}$.
Write $\phi\defeq [\grad{F}]\vert_X$.
By Grothendieck's Galois theory,
there exists a finite \'{e}tale Galois cover $\pi\maps X'\to X$ with $X'$ connected and $\phi\circ\pi\maps X'\to S$ (finite \'{e}tale) Galois.
Let $G\defeq\Aut_S(X')$ and $H\defeq\Aut_X(X')$.

The group $G$ acts transitively on $X'$.
So for any geometric points $[\bm{c}]\in S(\ol{\QQ})$ and $p\in X'_{[\bm{c}]}$, we can characterize the fiber $X_{[\bm{c}]}$ as the set $\set{\pi(gp):g\in H\backslash G}\belongs X(\ol{\QQ})$.

\subsubsection{Constructing a product ``divisible'' by \texpdf{$F^\vee$}{Fvee}}

View $F=F(\bm{x})$ as a section of $\mcal{O}_X(3)$.
Consider the $G$-equivariant line bundle $\mcal{L}\defeq \bigotimes_{g\in H\backslash G} g^\ast(\pi^\ast\mcal{O}_X(3))$ on $X'$.
The product
\begin{equation}
\label{EQN:define-G-invariant-product-alpha}
    \alpha\defeq \prod_{g\in H\backslash G} (g^\ast\pi^\ast{F})
\end{equation}
defines a $G$-invariant section of $\mcal{L}$ on $X'$.
(A \emph{$G$-invariant section} $\alpha\in\Gamma(X',\mcal{L})$ is equivalent in data to a $G$-equivariant morphism $\alpha\maps\mcal{O}_{X'}\to\mcal{L}$.)

For every geometric point $p\in X'(\ol{\QQ})$ with $\phi\pi(p)\in S\cap V^\vee$, there exists $g\in H\backslash G$ with $\pi(gp)\in X\cap V$, so that $(g^\ast\pi^\ast{F})(p) = F(\pi gp) = 0$.
So
\begin{equation*}
    \alpha\vert_{(\phi\pi)^{-1}(S\cap V^\vee)} = 0.
\end{equation*}
Therefore, by Galois descent (in the form of an equivalence of categories),
there exist a line bundle $\mcal{D}$ on $S$ with $\mcal{L}\cong (\phi\pi)^\ast\mcal{D}$,
and a section $\delta\in \Gamma(S, \mcal{D})$ vanishing along $S\cap V^\vee$,
with $\alpha=(\phi\pi)^\ast\delta$.
Let $\mcal{F}\defeq \phi^\ast\mcal{D}$ and $\beta\defeq \phi^\ast\delta\in \Gamma(X, \mcal{F})$;
then $\mcal{L}\cong \pi^\ast\mcal{F}$ and $\alpha=\pi^\ast\beta$.

But $S$, $X$ are hypersurface complements in $\PP^{m-1}$,
so $\Pic(\PP^{m-1})\to \Pic(S)$ and $\Pic(\PP^{m-1})\to \Pic(X)$ are surjective
and we may identify $\mcal{F}$, $\mcal{D}$ with suitable powers of $\mcal{O}_X(1)$, $\mcal{O}_S(1)$, respectively.
Then up to a choice of nonzero constant factors,
we may view $\beta$, $\delta$ as homogeneous rational functions
(i.e.~ratios of homogeneous $m$-variable polynomials)
satisfying $F^\vee(\bm{c})\mid\delta$ on $S$ and $F^\vee(\grad{F}(\bm{x}))=\phi^\ast{F^\vee}\mid\beta$ on $X$.
Here we interpret divisibility of two sections on a scheme to mean their ratio is a global section of the obvious ``tensor-quotient'' line bundle.

\subsubsection{``Factoring'' \texpdf{$F^\vee$}{Fvee}}

By the definition of $V^\vee$ and $F^\vee$ in \S\ref{SEC:alg-geom-background},
we have $F(\bm{x})\mid F^\vee(\grad{F}(\bm{x}))$ in $\QQ[\bm{x}]$.
So $F\mid \phi^\ast{F^\vee}$ on $X$.
By \eqref{EQN:define-G-invariant-product-alpha}, it follows that $\alpha\mid (\pi^\ast\phi^\ast{F^\vee})^{\card{H\backslash G}}$ on $X'$, since $g^\ast\pi^\ast\phi^\ast{F^\vee} = \pi^\ast\phi^\ast{F^\vee}$ for all $g\in G$.
Since $\alpha = \pi^\ast\phi^\ast\delta$, it follows (by Galois descent) that $\delta\mid (F^\vee)^{\card{H\backslash G}}$ on $S$.
Since $F^\vee\mid \delta$ on $S$, and $F^\vee$ is prime in $\QQ[\bm{c}]$, we conclude that there exists an integer $e\geq 1$ satisfying $(F^\vee)^e\mid \delta$ and $\delta\mid (F^\vee)^e$ on $S$.

We need to restrict to $V^\vee$.
Luckily,
\eqref{COND:L^perp-not-contained-in-branch-locus} and $\PP{L^\perp}\belongs V^\vee$ imply that $S\cap V^\vee \ne \emptyset$ (but see Question~\ref{QUES:dual-variety-vs-branch-locus}).
Furthermore, $V^\vee$ is geometrically reduced, so $V^\vee\setminus \map{Sing}(V^\vee)$ is a dense open subvariety of $V^\vee$ (by ``generic smoothness'').
Thus $(S\cap V^\vee)\setminus \map{Sing}(V^\vee)\neq \emptyset$.

Choose a geometric point $[\bm{c}]$ of $(S\cap V^\vee)\setminus\map{Sing}(V^\vee)$.
Biduality furnishes a \emph{unique} point $[\bm{x}]\in X_{[\bm{c}]}$ with $F(\bm{x})=0$.
So if $p\in X'_{[\bm{x}]}$ and $g\in G$,
then the section $g^\ast\pi^\ast{F}$ on $X'$ evaluates to $0$ at $p$ \emph{if and only if} $g\in H$.
Thus $\pi^\ast{F}\nmid\prod_{[1]\neq g\in H\backslash G}(g^\ast\pi^\ast{F})$ on $X'$.
It follows that
$(\pi^\ast{F})^2\nmid \alpha$,
whence $F^2\nmid \beta$;
whence $(F^\vee)^2\nmid \delta$.

Thus $e=1$.
In particular, $\delta\mid F^\vee$ on $S$, so $\alpha\mid (\phi\pi)^\ast{F^\vee}$ on $X'$.

\subsubsection{Differentiating the product}

Using \eqref{COND:L^perp-not-contained-in-branch-locus} one last time
(more seriously than before),
we will now complete the proof of the second part of Proposition~\ref{PROP:dual-linear-subvariety}.

\begin{proof}
[Proof of (2)]
By \eqref{COND:L^perp-not-contained-in-branch-locus},
$S\cap\PP{L^\perp}\neq\emptyset$.
Yet $\phi\vert_{X\cap\PP{L}} = \gamma\vert_{X\cap\PP{L}}\maps X\cap\PP{L}\to S\cap\PP{L^\perp}$ is finite \'{e}tale of degree $2^{m/2-1}$, by part~(1) and the definition of $S$.
Let $[\bm{c}]\in (S\cap \PP{L^\perp})(\ol{\QQ})$,
and fix $p\in X'_{[\bm{c}]}$.
Then there exist at least $2^{m/2-1}$ cosets $g\in H\backslash G$ with $\pi gp\in (X\cap\PP{L})(\ol{\QQ})\belongs V(\ol{\QQ})$.
Applying the Leibniz rule to \eqref{EQN:define-G-invariant-product-alpha}, after restricting to a small affine neighborhood of $p$, we thus get $j_p^r{\alpha}(p)=\bm{0}$ for $r\defeq 2^{m/2-1}-1$,
where $j^r\maps\mcal{L}\to J^r\mcal{L}$ denotes the $r$th-order jet map
``along'' $\mcal{L}$ (from $\mcal{L}$ to its $r$th jet bundle $J^r\mcal{L}$).

Since $\alpha\mid(\phi\pi)^\ast{F^\vee}$,
Leibniz then implies $j_p^r{(\phi\pi)^\ast{F^\vee}}(p)=\bm{0}$ ``along'' the pullback line bundle $(\phi\pi)^\ast\mcal{O}_S(\deg{F^\vee})$.
But $\phi\pi\maps X'\to S$ is \'{e}tale at $p\in X'$,
so $j_{[\bm{c}]}^r{F^\vee}([\bm{c}])=\bm{0}$ ``along'' $\mcal{O}_S(\deg{F^\vee})$ itself,
over all points $[\bm{c}]\in(S\cap\PP{L^\perp})(\ol{\QQ})$.
Finally,
$S\cap\PP{L^\perp}$ is dense in $\PP{L^\perp}$, so the vanishing of the $r$th-order jet section $j^r F^\vee$ extends to all points $[\bm{c}]\in\PP{L^\perp}$, as desired.
\end{proof}

\begin{remark}
\label{RMK:friendly-etale-local-calculus}
In the friendly setting above,
our \'{e}tale morphisms (such as $\phi\pi\maps X'\to S$),
after base change to an algebraically closed field,
always induce isomorphisms on completed local rings.
So we could do calculus purely in terms of formal power series.
\end{remark}

\subsection{The diagonal case}
\label{SUBSEC:diagonal-dual-example}

Say $m$ is even and $F$ is diagonal,
and write $F = F_1x_1^3 + \dots + F_mx_m^3$.
Then we can explicitly verify all the theory above.
Here $[\grad{F}]\maps [\bm{x}]\mapsto [3F_1x_1^2,\dots,3F_mx_m^2]$
is multi-quadratic Galois of degree $2^{m-1}$.
We proceed by studying $\Upsilon$ and $F^\vee$ combinatorially.



For combinatorial purposes, let $[n]\defeq \set{1,2,\dots,n}$ for each integer $n\geq 1$.

\begin{definition}
\label{DEFN:diagonal-permissible-pairing}
Let $\mcal{J}=(\mcal{J}(k))_{k\in \mcal{K}}$ denote an \emph{ordered set partition} of $[m]$:
a list of pairwise disjoint nonempty sets $\mcal{J}(k)\belongs [m]$ covering $[m]$, indexed by a set $\mcal{K}\in \set{[1], [2], [3], \ldots}$.

\begin{enumerate}
    \item Call $\mcal{J}$, $\mcal{J}'$ \emph{equivalent} if they define the same unordered partition of $[m]$ (i.e.~if $\mcal{K}=\mcal{K}'$ and $\set{\mcal{J}(k): k\in \mcal{K}}=\set{\mcal{J}'(k): k\in \mcal{K}'}$).

    \item Call $\mcal{J}$ a \emph{pairing} if
    $\card{\mcal{J}(k)} = 2$ for all $k\in \mcal{K}$.

    \item Call $\mcal{J}$ \emph{permissible} if
    for all $k\in \mcal{K}$ and $i,j \in\mcal{J}(k)$, we have $F_j/F_i\in (\QQ^\times)^3$.
    For a permissible $\mcal{J}$,
    let $\mcal{R}_{\mcal{J}}
    \defeq \set{\bm{c}\in\ZZ^m
    : \textnormal{if $k\in \mcal{K}$ and $i,j \in\mcal{J}(k)$, then $c_i/F_i^{1/3}=c_j/F_j^{1/3}$}}$;
    and given $\bm{c}\in \mcal{R}_{\mcal{J}}$,
    define $c\maps \mcal{K}\to \RR$ so that
    if $k\in \mcal{K}$ and $i \in\mcal{J}(k)$, then $c_i/F_i^{1/3}=c(k)$.
\end{enumerate}
\end{definition}




We now recall the well-known construction of the $\frac{m}{2}$-dimensional vector spaces $L/\CC$ with $F\vert_L = 0$.
Each equivalence class of pairings $\mcal{J}$ yields $3^{m/2}$ distinct $L/\CC$,
obtained by setting $F_ix_i^3+F_jx_j^3=0$ for each part $\mcal{J}(k)=\set{i,j}$.
Over $\QQ$,
we must set $x_i+(F_j/F_i)^{1/3}x_j=0$---which is
valid when $F_i\equiv F_j\bmod{(\QQ^\times)^3}$.
The vectors spaces $L/\CC$ and $L/\QQ$ thus constructed are known to be the only possibilities \cite{wang2022thesis}*{Remark~6.3.8}.
Hence the following holds:

\begin{proposition}
\label{PROP:characterize-diagonal-L/Q}
There is a canonical bijection,
between $\Upsilon$ and the set of equivalence classes of permissible pairings $\mcal{J}$,
characterized by $L\cap\ZZ^m = \mcal{R}_{\mcal{J}}^\perp$ (an equality of sublattices of $\ZZ^m$).
\end{proposition}


Next, we turn to $F^\vee$.
For convenience,
fix square roots $F_i^{1/2}\in \ol{\QQ}^\times$.
For some constant $\beta_{F_1,\dots,F_m}\in \QQ^\times$,
the polynomial $F^\vee(\bm{c})$ factors in $\ol{\QQ}[\bm{c}^{1/2}]$ as follows:
\begin{equation}
\label{EQN:diagonal-discriminant-factorization}
F^\vee(\bm{c}) = \beta_{F_1,\dots,F_m} \cdot
\prod_{\bm{\eps}}(\eps_1F_1^{-1/2}c_1^{3/2}+\eps_2F_2^{-1/2}c_2^{3/2}+\dots+\eps_mF_m^{-1/2}c_m^{3/2}) \in \QQ[\bm{c}],
\end{equation}
with the product taken over $\bm{\eps} = (\eps_1,\dots,\eps_m)$ with $\eps_1=1$ and $\eps_2,\dots,\eps_m=\pm1$.
(This classical formula is a simple consequence of Definition~\ref{DEFN:define-Vdual-as-scheme-theoretic-image} and the Jacobian criterion.)



Let $\mscr{L}_{\bm{\eps}}\defeq \eps_1F_1^{-1/2}c_1^{3/2}+\dots+\eps_mF_m^{-1/2}c_m^{3/2} \in \ol{\QQ}[\bm{c}^{1/2}]$.
For each $\bm{c}\in \QQ^m\setminus \set{\bm{0}}$, fix square roots $c_i^{1/2}\in \ol{\QQ}$.
Using formal power series calculus over variables $c_i\neq 0$
(by Remark~\ref{RMK:friendly-etale-local-calculus}, adapted to $\Aff^1_{\ol{\QQ}}\to\Aff^1_{\ol{\QQ}},\;t\mapsto t^2$ away from the origin),
we will prove the following result,
which precisely characterizes the order of vanishing of $F^\vee$ at $\bm{c}$:

\begin{proposition}
\label{PROP:characterize-diagonal-Disc-vanishing-order}
Fix $r\geq 0$ and $\bm{c}\in \QQ^m\setminus \set{\bm{0}}$.
Then the following are equivalent:
(1) $j^r{F^\vee}$ vanishes at $\bm{c}$,
and (2) there exist at least $r+1$ distinct $\bm{\eps}$ with $\mscr{L}_{\bm{\eps}}=0$.
\end{proposition}

\begin{proof}
We first prove (1)$\Rightarrow$(2), by induction on $r\geq 0$.
The base case $r=0$ follows directly from \eqref{EQN:diagonal-discriminant-factorization}.
Now fix $r\geq 1$,
and assume the implication (1)$\Rightarrow$(2) holds for $r-1$.
Suppose (1) holds; we wish to prove that (2) holds.

By the inductive hypothesis, $\#\set{\bm{\eps}: \mscr{L}_{\bm{\eps}}=0}\geq r$.
To rule out the possibility that $\#\set{\bm{\eps}: \mscr{L}_{\bm{\eps}}=0} = r$,
we work with ``pure'' derivatives $\partial_{c_i}^{\leq r}$, for just a single index $i$ with $c_i\neq 0$.
For example,
if $c_1\neq 0$, and $\#\set{\bm{\eps}: \mscr{L}_{\bm{\eps}}=0} = r$,
then the product rule and \eqref{EQN:diagonal-discriminant-factorization} imply
\begin{equation*}
\partial_{c_1}^r{F^\vee}(\bm{c})
= \beta_{F_1,\dots,F_m} \cdot r! \cdot (\tfrac32 F_1^{-1/2}c_1^{1/2})^r
\prod_{\bm{\eps}:\, \mscr{L}_{\bm{\eps}}\neq 0} \mscr{L}_{\bm{\eps}}
\neq 0,
\end{equation*}
contradicting (1).
Thus $\#\set{\bm{\eps}: \mscr{L}_{\bm{\eps}}=0}\geq r+1$.
This completes the induction.

It remains to prove (2)$\Rightarrow$(1).
To avoid confusion, rename our given $\bm{c}$ to $\bm{\xi}$.
Say $\xi_i=0$ for $i\in I$, and $\xi_i\neq 0$ for $i\in [m]\setminus I$.
We must take extra care over $i\in I$.

Let $\mscr{T}_I$ be the set of triples $(\bm{a},\bm{b},E)$, with $(\bm{a}, \bm{b})\in \ZZ_{\geq 0}^I \times \ZZ^{[m]\setminus I}$ and $E\belongs \set{\bm{\eps}\in \set{\pm1}^m: \eps_1=1}$, such that for each $i\in I$, the set $E\cup (-E)$ is invariant under the flip $\eps_i\mapsto -\eps_i$.
Consider the following ``formal analytic functions'' (inspired by \eqref{EQN:diagonal-discriminant-factorization}), indexed by $(\bm{a},\bm{b},E)\in \mscr{T}_I$:
\begin{equation}
\label{EXPR:formal-analytic-functions}
\biggl(\,\prod_{i\in I} c_i^{a_i}\biggr)
\biggl(\,\prod_{i\notin I} c_i^{b_i/2}\biggr)
\biggl(\,\prod_{\bm{\eps}\in E} \mscr{L}_{\bm{\eps}}\biggr)
\in \ol{\QQ}[c_i]_{i\in I}[c_i^{1/2},c_i^{-1/2}]_{i\notin I}.
\end{equation}
The functions \eqref{EXPR:formal-analytic-functions} span a vector space over $\ol{\QQ}$ that contains $F^\vee$ and is closed under differentiation in $\bm{c}$.
In fact,
differentiating \eqref{EXPR:formal-analytic-functions} in $c_i$ leads to
terms with $a_i\mapsto a_i-1$ or $(a_i,\card{E})\mapsto (a_i+2,\card{E}-2)$ if $i\in I$,
and to terms with $b_i\mapsto b_i-2$ or $(b_i,\card{E})\mapsto (b_i+1,\card{E}-1)$ if $i\notin I$.
In each case,
applying $\partial_{c_i}$ decreases $\min_{\bm{a},\bm{b},E}(\abs{\bm{a}}+\card{E})$ by at most $1$.

Now suppose (2) holds,
and fix $\bm{r}\in \ZZ_{\geq 0}^m$ with $\abs{\bm{r}}\leq r$.
Then $\partial_{\bm{c}}^{\bm{r}}{F^\vee}$ is a $\ol{\QQ}$-linear combination of functions \eqref{EXPR:formal-analytic-functions} indexed by
triples $(\bm{a},\bm{b},E)\in \mscr{T}_I$ with $\abs{\bm{a}}+\card{E}\geq 2^{m-1}-r$ (and thus $\card{E}\geq 2^{m-1}-r$ or $\abs{\bm{a}}\geq 1$).
By (2), each such function must vanish at our original given point $\bm{c}=\bm{\xi}$.
Thus $\partial_{\bm{c}}^{\bm{r}}{F^\vee}(\bm{\xi}) = 0$.
So (1) holds.
\end{proof}

\begin{remark}
A short computation yields the equality
\begin{equation}
\label{EQN:branching-interpretation-of-number-of-vanishing-L_eps}
\#\set{\bm{\eps}: \mscr{L}_{\bm{\eps}}=0}
= \sum_{[\bm{x}]\in \gamma(\overline{\QQ})^{-1}([\bm{c}])} 2^{\#\set{i\in [m]:\, x_i=0}},
\end{equation}
where $\gamma(\overline{\QQ})^{-1}([\bm{c}])\defeq
\set{[\bm{x}]\in V(\ol{\QQ}): [\grad{F}(\bm{x})]=[\bm{c}]}
= \set{\textnormal{singular $\ol{\QQ}$-points of $V_{\bm{c}}$}}$.
(Here $x_i$ corresponds to $\eps_iF_i^{-1/2}c_i^{1/2}$,
with some ambiguity or ``multiplicity'' in $\eps_i$ when $x_i=0$.)
Using \eqref{EQN:branching-interpretation-of-number-of-vanishing-L_eps}, one can formulate Proposition~\ref{PROP:characterize-diagonal-Disc-vanishing-order} more geometrically,
without reference to $\bm{\eps}$'s.
Does this geometric formulation extend somehow to more general $F$?
\end{remark}


Finally, we analyze the interaction between $\Upsilon$, $F^\vee$.
Fix $L\in \Upsilon$.
By Proposition~\ref{PROP:characterize-diagonal-L/Q}, $L$ corresponds to some permissible pairing $\mcal{J}$.
Proposition~\ref{PROP:characterize-diagonal-Disc-vanishing-order} has the following corollary:

\begin{corollary}
\label{COR:baseline-diagonal-Disc-vanishing}
For $L$ as above,
we have $(j^{2^{m/2-1}-1}{F^\vee})\vert_{L^\perp}=\bm{0}$.
\end{corollary}

\begin{proof}
For each part $\mcal{J}(k)=\set{i,j}$,
there are exactly two choices of signs $(\eps_i,\eps_j)\in\set{\pm1}^2$---or
only one choice if $1\in\mcal{J}(k)$---such that
$\eps_iF_i^{-1/2}c_i^{3/2}+\eps_jF_j^{-1/2}c_j^{3/2}$ vanishes over \emph{all} $\bm{c}\in L^\perp\cap\ZZ^m = \mcal{R}_{\mcal{J}}$ lying in a given orthant of $\RR^m$.
So given $\bm{c}\in L^\perp\setminus\set{\bm{0}}$,
we can apply Proposition~\ref{PROP:characterize-diagonal-Disc-vanishing-order}(2)$\Rightarrow$(1) with $r\defeq 2^{m/2-1}-1$.
Since $L^\perp\setminus\set{\bm{0}}$ is dense in $L^\perp$, the vanishing of $j^r{F^\vee}$ then extends to all of $L^\perp$.
\end{proof}

Thus we have explicitly verified the conclusion of Proposition~\ref{PROP:dual-linear-subvariety}.
The next result shows that in fact,
$F^\vee$ generally does not vanish to higher order along $L^\perp$
(and furthermore, \ref{OBS:higher-order-diagonal-Disc-vanishing-critera}(1) gives us a simple description of when higher vanishing occurs).
Let $s\maps \QQ\to \QQ,\;q\mapsto q^2$.

\begin{observation}
\label{OBS:higher-order-diagonal-Disc-vanishing-critera}
Given $L$, $\mcal{J}$ as above,
fix $\bm{c}\in L^\perp\cap\ZZ^m = \mcal{R}_{\mcal{J}}$.
Define $c(k)$ as in Definition~\ref{DEFN:diagonal-permissible-pairing}.
For each $k\in \mcal{K}$,
fix a square root $c(k)^{3/2}\in \ol{\QQ}$ of $c(k)^3$.
Then the following hold:
\begin{enumerate}
    \item $j^{2^{m/2-1}}{F^\vee}(\bm{c})=\bm{0}$
    if and only if
    there exist $l\geq 1$ distinct indices $k_1,\dots,k_l\in \mcal{K}$ such that
    $c(k_1)^{3/2}\pm\cdots\pm c(k_l)^{3/2}=0$ holds for some choice of signs.
    
    \item If $j^{2^{m/2-1}}{F^\vee}(\bm{c})=\bm{0}$,
    then $c(k_1)^3c(k_2)^3\in s(\QQ)$ for some distinct $k_1,k_2\in \mcal{K}$.
\end{enumerate}
\end{observation}

\begin{proof}
(1):
If $\bm{c}=\bm{0}$, then $j^{2^{m/2-1}}{F^\vee}(\bm{c})=\bm{0}$ (since $\deg{F^\vee} = 3\cdot 2^{m-2}\geq 1+2^{m/2-1}$), and $c(k)^3 = 0$ for all $k\in \mcal{K}$.
If $\bm{c}\neq \bm{0}$,
apply Proposition~\ref{PROP:characterize-diagonal-Disc-vanishing-order}(1)$\Leftrightarrow$(2) with $r\defeq 2^{m/2-1}$
(and then ``simplify''  the condition \ref{PROP:characterize-diagonal-Disc-vanishing-order}(2) using the fact that $\mcal{J}$ is a pairing).

(2):
Suppose $j^{2^{m/2-1}}{F^\vee}(\bm{c})=\bm{0}$.
Apply (1); choose $k_1,\dots,k_l\in \mcal{K}$ with $l$ minimal.
We now do casework on $l$, making use of minimality if $l\geq 2$.
\begin{itemize}
    \item If $l=1$,
    then $c(k)^3=0$ for some $k\in \mcal{K}$.
    
    \item If $l=2$,
    then $c(k_1)^3=c(k_2)^3\in\QQ^\times$ for some distinct $k_1,k_2\in \mcal{K}$.
    
    \item If $l\geq 3$,
    then $c(k_t)^3\in\QQ^\times$ for all $t\in[l]$,
    and by multi-quadratic field theory in characteristic $0$,
    the square classes $c(k_1)^3,\dots,c(k_l)^3\bmod{(\QQ^\times)^2}$ must all coincide.
    (In fact,
    say we fix $i_t\in \mcal{J}(k_t)$ for $t\in [l]$.
    Then $c_{i_t}/F_{i_t}=x_{i_t}^2d$
    for some $d, x_{i_t}\in \QQ^\times$ such that $F_{i_1}x_{i_1}^3+\dots+F_{i_l}x_{i_l}^3=0$.
    So $c(k_t)^3 = F_{i_t}^2x_{i_t}^6d^3$.)
\end{itemize}
In each case,
$c(k_1)^3c(k_2)^3\in s(\QQ)$ holds for some distinct $k_1,k_2\in \mcal{K}$.
\end{proof}

\begin{remark}
\label{RMK:char-p-diagonal-Disc-analysis}
Though written in characteristic $0$,
the main results of \S\ref{SUBSEC:diagonal-dual-example}
carry over to arbitrary fields of characteristic $p\nmid (6^m)!F_1\cdots F_m$,
under \eqref{INEQ:diagonal-Fdual-primitive-scaling-convention}.
Such extensions of
Proposition~\ref{PROP:characterize-diagonal-Disc-vanishing-order},
\eqref{EQN:branching-interpretation-of-number-of-vanishing-L_eps},
Corollary~\ref{COR:baseline-diagonal-Disc-vanishing},
and Observation~\ref{OBS:higher-order-diagonal-Disc-vanishing-critera}(1)
to $\FF_p$
will prove useful in \S\ref{SEC:proving-bias-in-linear-S_c's}.
\end{remark}

\section{Poisson summation for the endgame}
\label{SEC:reverse-sums-by-duality}

Suppose $2\mid m$.
Fix $L\in \Upsilon$, and recall Definition~\ref{DEFN:primitive-sublattices}.
Then $F\vert_{\Lambda} = 0$,
and $\Lambda$, $\Lambda^\perp$ are primitive rank-$\frac{m}{2}$ sublattices of $\ZZ^m$.
Now choose bases $\bm{\Lambda}$, $\bm{\Lambda^\perp}$ of $\Lambda$, $\Lambda^\perp$, respectively.
Identify $\bm{\Lambda}$, $\bm{\Lambda^\perp}$ with
$m\times\frac{m}{2}$ and $\frac{m}{2}\times m$ integer matrices, respectively,
so that $\Lambda=\bm{\Lambda}\ZZ^{m/2}$ and $\Lambda^\perp=\ZZ^{m/2}\bm{\Lambda^\perp}$
(where we view $\Lambda$ as a ``column space'' and $\Lambda^\perp$ as a ``row space'').

We seek to prove
Lemma~\ref{LEM:reverse-integral-averaging} and Proposition~\ref{PROP:reverse-exponential-sum-averaging} below, about certain averages over a given shifted dilate of $\Lambda^\perp$
(for the ``endgame'' of \S\ref{SEC:full-proof-outline}).

For the rest of \S\ref{SEC:reverse-sums-by-duality},
let $\bm{x}$, $\bm{h}$, $\bm{x}'$ denote column vectors and $\bm{c}$, $\bm{b}$, $\bm{v}$ row vectors.
In particular,
the dot product $\bm{c}\cdot \bm{x}$ then coincides with standard matrix multiplication.

\subsection{Preliminaries}

We have $\Lambda = (\Lambda^\perp)^\perp$,
i.e.~$\Lambda = \set{\bm{x}\in\ZZ^m:\bm{\Lambda^\perp}\bm{x}=\bm{0}}$.
Let
\begin{equation}
\label{EQN:real-density-linear-spaces}
\sigma_{\infty,L^\perp,w}
\defeq\lim_{\eps\to 0}{(2\eps)^{-m/2} \int_{\bm{\Lambda^\perp}\tilde{\bm{x}}\in [-\eps,\eps]^{m/2}} d\tilde{\bm{x}}\,w(\tilde{\bm{x}})}.
\end{equation}
Then the formula \eqref{EQN:intro-Lambda-sum-re-interpretation-formula} holds,
by Poisson summation over $L\cap\ZZ^m=\Lambda$
(or, at least morally,
by the circle method applied to $\bm{\Lambda^\perp}\bm{x}=\bm{0}$).
In particular, $\sigma_{\infty,L^\perp,w}$ does not depend on the choice of $\bm{\Lambda^\perp}$.
For an alternative interpretation of $\sigma_{\infty,L^\perp,w}$, see the second part of Lemma~\ref{LEM:reverse-integral-averaging}.

For calculations to come,
it will help to extend $\bm{\Lambda^\perp}$ to a basis of $\ZZ^m$.

\begin{definition}
By primitivity of $\Lambda^\perp$,
choose $\Gamma$ (itself primitive) such that $\ZZ^m=\Lambda^\perp\oplus\Gamma$.
Then choose a $\frac{m}{2}\times m$ basis matrix $\mathbf{\Gamma}$ so that $\Gamma=\ZZ^{m/2}\bm{\Gamma}$.
\end{definition}

The rows of the $m\times m$ matrix $\inmat{\bm{\Lambda^\perp} \\ \bm{\Gamma}}$ form a basis of $\ZZ^m$.
Therefore,
\begin{equation}
\label{EQN:key-matrix-is-unimodular}
    \det\inmat{\bm{\Lambda^\perp} \\ \bm{\Gamma}} = \pm 1.
\end{equation}


Let $R$ denote a ring, e.g.~$\RR$ or $\ZZ/n\ZZ$.
Given a $\ZZ$-module $A$, let $A_R\defeq A\otimes R$.
To study Fourier transforms over $\bm{x}\in R^m$ at $\bm{c}\in \Lambda^\perp$,
it will help to rewrite $\bm{c}\cdot \bm{x}$.

\begin{definition}
\label{DEFN:c^star,psi,psi_R}
Using the definition of $\bm{\Lambda^\perp}$ as a basis matrix,
let $\psi\maps \Lambda^\perp\to \ZZ^{m/2}$ be the $\ZZ$-linear isomorphism $\bm{c}\mapsto \bm{c}^\star$ defined by the equation $\bm{c} = \bm{c}^\star\bm{\Lambda^\perp}$.
Let $\psi_R\defeq \psi\otimes R$ be the $R$-linear isomorphism $(\Lambda^\perp)_R\to R^{m/2}$ induced by $\psi$.
\end{definition}

For each $\bm{x}\in R^m$, let $\inmat{\bm{h} \\ \bm{x}'}\defeq \inmat{\bm{\Lambda^\perp} \\ \bm{\Gamma}}\bm{x}$.
Equivalently, let
\begin{equation}
\label{EQN:define-h,x'}
\bm{h}\defeq \bm{\Lambda^\perp}\bm{x}\in R^{m/2},
\qquad \bm{x}'\defeq \bm{\Gamma}\bm{x}\in R^{m/2}.
\end{equation}

\begin{example}
If $F=x_1^3+\dots+x_m^3$ and
$\Lambda^\perp=\mcal{R}_{\mcal{J}}$ (in the notation of Proposition~\ref{PROP:characterize-diagonal-L/Q}),
then we can choose $\bm{\Lambda^\perp}$ so that $h_k = \sum_{i\in \mcal{J}(k)} x_i$ for all $k\in \mcal{K}$.
\end{example}

We need two more algebraic propositions.
Over $R$, let $\bm{c}\cdot \bm{x}\maps (\Lambda^\perp)_R\times R^m\to R$ denote the ``obvious'' $R$-bilinear map induced by the usual dot product $\bm{c}\cdot \bm{x}\maps \Lambda^\perp\times \ZZ^m\to \ZZ$.

\begin{proposition}
\label{PROP:Lambda^perp-dot-product-coordinate-change}
If $\bm{c}\in (\Lambda^\perp)_R$ and $\bm{x}\in R^m$, then
$\bm{c}\cdot \bm{x} = \psi_R(\bm{c})\cdot \bm{h}$.
\end{proposition}

\begin{proof}
By $R$-linearity properties, reduce to the case $R=\ZZ$, which is trivial.
\end{proof}

Let $\Lambda\cdot R\belongs R^m$ be the $R$-module generated by $\Lambda$ (via the composite map $\Lambda\to \ZZ^m\to R^m$).

\begin{proposition}
\label{PROP:h=0-cuts-out-Lambda}
Let $\bm{x}\in R^m$.
Then $\bm{h}=\bm{0}$ if and only if $\bm{x}\in \Lambda\cdot R$.
\end{proposition}

\begin{proof}
The multiplication map $\bm{\Lambda^\perp}\maps \ZZ^m\to \ZZ^{m/2}$ (given by $\bm{x}\mapsto \bm{\Lambda^\perp}\bm{x}$) is surjective by \eqref{EQN:key-matrix-is-unimodular},
and thus defines an exact sequence $\Lambda\to\ZZ^m\to\ZZ^{m/2}\to0$.
The right exact functor $\otimes R$ therefore gives an exact sequence $\Lambda_R\to R^m\to R^{m/2}\to0$.
So
\begin{equation*}
    \ker(\bm{\Lambda^\perp}\maps R^m\to R^{m/2}) = \im(\Lambda_R\to R^m).
\end{equation*}
But $\im(\Lambda_R\to R^m) = \Lambda\cdot R$, since $\Lambda_R$
is $R$-linearly generated by $\Lambda$.
\end{proof}


\subsection{Averaging the oscillatory integrals}

Recall $I_{\bm{c}}(n)$ from \eqref{EQN:define-I_c(n)}.
Let $R\defeq \RR$ in \eqref{EQN:define-h,x'}.
Since $\bm{x}\mapsto (\bm{h},\bm{x}')$ is unimodular by \eqref{EQN:key-matrix-is-unimodular}, we have
\begin{equation}
\label{EQN:I_c(n)-unimodular-coordinate-change}
I_{\bm{c}}(n)
= \int_{(\bm{h},\bm{x}')\in \RR^m} d\bm{h}\,d\bm{x}'\,
w(\bm{x}/X)
h(n/Y, F(\bm{x})/Y^2)
e(-\bm{c}\cdot \bm{x}/n),
\end{equation}
an integral in which we view $\bm{x}$ as a function of $(\bm{h},\bm{x}')$.
For each $\bm{h}\in \ZZ^{m/2}$, let
\begin{equation*}
J(\bm{h};n)\defeq \int_{\bm{x}'\in \RR^{m/2}} d\bm{x}'\,
w(\bm{x}/X) h(n/Y, F(\bm{x})/Y^2).
\end{equation*}
(The integral $J(\bm{h};n)$ corresponds to the integral $J_q(\bm{j})$ on \cite{heath1998circle}*{p.~692}.
But in the present paper, we do not make any serious use of $J(\bm{h};n)$ for $\bm{h}\neq \bm{0}$.)

\begin{lemma}
\label{LEM:reverse-integral-averaging}
Suppose $n = n_0n_1$ (where $n_0,n_1\geq 1$ are integers) and $\bm{b}\in \Lambda^\perp$.
Then
\begin{equation}
\label{EQN:Poisson-sum-result-over-Lambda^perp}
n_1^{-m/2} \sum_{\bm{c}\in \bm{b}+n_0\Lambda^\perp} I_{\bm{c}}(n)
= \sum_{\bm{h}\in n_1\ZZ^{m/2}} e_n(-\bm{b}^\star\cdot \bm{h}) J(\bm{h};n).
\end{equation}
Here $J(\bm{0};n) = \sigma_{\infty,L^\perp,w}X^{m/2}\cdot h(n/Y,0)$.
Furthermore,
if $n_1 \geq M_1 X$ for a sufficiently large positive real $M_1 \ll_{\bm{\Lambda^\perp},w} 1$,
then $J(\bm{h};n)=0$ for all nonzero $\bm{h}\in n_1\ZZ^{m/2}$.
\end{lemma}

\begin{proof}
Write $\bm{c}=\bm{b}+n_0\bm{v}$ for $\bm{v}\in \Lambda^\perp$.
Define $\bm{c}^\star,\bm{b}^\star,\bm{v}^\star\in\ZZ^{m/2}$ using Definition~\ref{DEFN:c^star,psi,psi_R}.
Proposition~\ref{PROP:Lambda^perp-dot-product-coordinate-change} then delivers the equality
$\bm{c}\cdot\bm{x}/n
= \bm{c}^\star\cdot\bm{h}/n
= \bm{b}^\star\cdot\bm{h}/n
+ \bm{v}^\star\cdot\bm{h}/n_1$.
So by \eqref{EQN:I_c(n)-unimodular-coordinate-change}, the integral $I_{\bm{c}}(n)=I_{\bm{b}+n_0\bm{v}}(n)$ is the Fourier transform at $\bm{v}^\star/n_1\in\RR^{m/2}$ of
the function
\begin{equation*}
\RR^{m/2}\to \RR,
\quad
\bm{h}
\mapsto
e(-\bm{b}^\star\cdot\bm{h}/n)
J(\bm{h};n).
\end{equation*}
Poisson summation over $\bm{h}\in n_1\ZZ^{m/2}$ hence yields
\begin{equation*}
\sum_{\bm{h}\in n_1\ZZ^{m/2}} e(-\bm{b}^\star\cdot\bm{h}/n) J(\bm{h};n)
= n_1^{-m/2} \sum_{\bm{v}^\star\in \ZZ^{m/2}} I_{\bm{b}+n_0\bm{v}}(n),
\end{equation*}
which implies \eqref{EQN:Poisson-sum-result-over-Lambda^perp} (upon recalling the correspondence between $\bm{v}^\star\in \ZZ^{m/2}$ and $\bm{v}\in \Lambda^\perp$).

We now compute $J(\bm{0};n)$.
Let $\tilde{\bm{x}}\defeq \bm{x}/X$.
Since $F\vert_{\Lambda\cdot\RR}=0$,
we have
\begin{equation*}
J(\bm{0};n) = X^{m/2}\cdot h(n/Y,0)\cdot
\int_{(\tilde{\bm{h}},\tilde{\bm{x}}')\in\set{\bm{0}}\times\RR^{m/2}}d\tilde{\bm{x}}'\,w(\tilde{\bm{x}}).
\end{equation*}
But an $\eps$-thickening in $\tilde{\bm{h}}$,
followed by an application of \eqref{EQN:real-density-linear-spaces}, yields
\begin{equation*}
\int_{(\tilde{\bm{h}},\tilde{\bm{x}}')\in\set{\bm{0}}\times\RR^{m/2}}d\tilde{\bm{x}}'\,w(\tilde{\bm{x}})
= \lim_{\eps\to 0}{\int_{(\tilde{\bm{h}},\tilde{\bm{x}}')\in\RR^m} d\tilde{\bm{h}}\,d\tilde{\bm{x}}'\,w(\tilde{\bm{x}})
\cdot \frac{\bm{1}_{\tilde{\bm{h}}\in [-\eps,\eps]^{m/2}}}{(2\eps)^{m/2}}}
= \sigma_{\infty,L^\perp,w}.
\end{equation*}

Finally, suppose $\bm{h}\in n_1\ZZ^{m/2}$ and $J(\bm{h};n)\neq 0$.
Take $\bm{x}'\in \RR^{m/2}$ with $w(\bm{x}/X)\neq 0$.
If $n_1\geq M_1X$,
then $n_1\mid \bm{h}=\bm{\Lambda^\perp}\bm{x}\ll X$,
so $\bm{h}=\bm{0}$ if $M_1$ is sufficiently large.
\end{proof}

\subsection{Vertically averaging the exponential sums}

Recall $S_{\bm{c}}(n)$ from \eqref{EQN:define-S_c(n)}, for each tuple $\bm{c}\in \ZZ^m$ and integer $n\geq 1$.
Let $\phi(n)$ denote Euler's totient function.

\begin{proposition}
\label{PROP:reverse-exponential-sum-averaging}
The quantity $S_{\bm{c}}(n)$ is a function of $n$ and $\bm{c}\bmod{n}$.
If $\bm{j}\in \ZZ^{m/2}$,
then
\begin{equation}
\label{EQN:general-j-c^star-Fourier-average}
\EE_{\bm{c}\in\Lambda^\perp/n\Lambda^\perp}[S_{\bm{c}}(n)
e_n(-\bm{c}^\star\cdot\bm{j})]
= \sum_{1\leq a\leq n:\, \gcd(a,n)=1}\,
\sum_{1\leq \bm{x}\leq n}
e_n(aF(\bm{x}))\cdot \bm{1}_{n\mid\bm{h}-\bm{j}}
\end{equation}
(where the left-hand side is well-defined).
In particular,
$\EE_{\bm{c}\in\Lambda^\perp/n\Lambda^\perp}[S_{\bm{c}}(n)] = \phi(n)n^{m/2}$.
\end{proposition}

\begin{proof}
The first sentence is clear by definition of $S_{\bm{c}}(n)$.
(Here $\bm{c}\bmod{n}$ refers to $\bm{c}\bmod{n\ZZ^m}$.)

A priori, if $\bm{c}\in \Lambda^\perp$, then $\bm{c}\bmod{n\Lambda^\perp}$ determines $\bm{c}\bmod{n}$.
So \eqref{EQN:general-j-c^star-Fourier-average} makes sense.
Since $\Lambda^\perp$ is primitive,
one can say more about ``congruence modulo $n$'',
but we need not do so.

Now fix $\bm{j}\in \ZZ^{m/2}$.
Let $R\defeq \ZZ/n\ZZ$.
Via the map $\psi_R$ from Definition~\ref{DEFN:c^star,psi,psi_R},
elements $\bm{c}\in \Lambda^\perp/n\Lambda^\perp = (\Lambda^\perp)_R$
correspond isomorphically to $\bm{c}^\star\in R^{m/2}=(\ZZ/n\ZZ)^{m/2}$,
so that if $\bm{x}\in R^m$, then $\bm{c}\cdot\bm{x} = \bm{c}^\star\cdot\bm{h}$ by Proposition~\ref{PROP:Lambda^perp-dot-product-coordinate-change}.
By \eqref{EQN:define-S_c(n)}, it follows that
\begin{equation*}
\EE_{\bm{c}\in\Lambda^\perp/n\Lambda^\perp}[S_{\bm{c}}(n)
e_n(-\bm{c}^\star\cdot\bm{j})]
= \sum_{a\in R^\times}
\sum_{\bm{x}\in R^m}
e_n(aF(\bm{x}))\cdot \EE_{\bm{c}^\star\in R^{m/2}}[e_n(\bm{c}^\star\cdot\bm{h}-\bm{c}^\star\cdot\bm{j})].
\end{equation*}
Summing over $\bm{c}^\star$ gives \eqref{EQN:general-j-c^star-Fourier-average}.
Furthermore, if we take $\bm{j}=\bm{0}$, then by Proposition~\ref{PROP:h=0-cuts-out-Lambda} and the identity $F\vert_{\Lambda\cdot R} = 0$, the right-hand side of \eqref{EQN:general-j-c^star-Fourier-average} simplifies to $\phi(n)n^{m/2}$.
\end{proof}

Let $T(\bm{j};n)$ denote the right-hand side of \eqref{EQN:general-j-c^star-Fourier-average}.
Then $T(\bm{j};n)$ corresponds to the sum $T_q(\bm{j})$ on \cite{heath1998circle}*{p.~692}.
It would be interesting to better understand $T(\bm{j};n)$ in general.
The sums $T(\bm{j};n)$ are multiplicative in $n$, and related to the singular series of
certain $\frac{m}{2}$-variable affine quadrics varying with $\bm{j}$.
Below, however, we only make use of the $\bm{j}=\bm{0}$ case of \eqref{EQN:general-j-c^star-Fourier-average}.

\section{General upper bounds}
\label{SEC:worst-case-Disc-0-sparse-bound}

We first recall some background on \eqref{EQN:delta-method}.
Let $\norm{\bm{c}}\defeq \max(\abs{c_1},\dots,\abs{c_m})$.
Recall \eqref{DEFN:normalize-S_c(n)}.

\begin{proposition}
[See e.g.~\cite{wang2022thesis}*{\S3.1},
or \cite{heath1998circle}*{(2.3)--(3.7)}]
\label{PROP:standard-delta-method-n,c-cutoffs}
The following hold:
\begin{enumerate}
    \item $I_{\bm{c}}(n)$ vanishes over $n\geq M_2 Y$, for some positive $M_2\ll_{F,w} 1$ independent of $\bm{c}\in \ZZ^m$.
    
    \item For any fixed $\eps,A>0$, we have $I_{\bm{c}}(n)\ll_{F,w,\eps,A} X^{-A}$ over $\norm{\bm{c}}\geq X^{1/2+\eps}$.
    
    \item We have
    \begin{equation*}
        \sum_{n\geq 1} \sum_{\bm{c}\in \ZZ^m} n^{(1-m)/2}\abs{S^\natural_{\bm{c}}(n)I_{\bm{c}}(n)}
        = \sum_{n\geq 1} \sum_{\bm{c}\in \ZZ^m} n^{-m}\abs{S_{\bm{c}}(n)I_{\bm{c}}(n)}
        < \infty.
    \end{equation*}
\end{enumerate}
\end{proposition}

Now fix a set $\mcal{S}\belongs \set{\bm{c}\in \ZZ^m: F^\vee(\bm{c})=0}$ and a real $\delta\geq 0$.
Suppose $\mcal{S}\cap [-C,C]^m$ has size $O(C^{m/2-\delta})$ for all reals $C\geq 1$.
Let
\begin{equation}
\label{EQN:sparse-absolute-contribution-f(S)-to-delta-method}
f(\mcal{S})\defeq Y^{-2} \sum_{\bm{c}\in \mcal{S}\setminus \set{\bm{0}}}
\sum_{n\geq 1} \abs{I_{\bm{c}}(n)} \cdot n^{1-m/2}
\cdot \sum_{n_\star\mid n} n_\star^{-1/2}\abs{S^\natural_{\bm{c}}(n_\star)}.
\end{equation}
At several points
in \S\ref{SEC:full-proof-outline},
Lemma~\ref{LEM:sparse-bound} will let us cleanly discard $f(\mcal{S})$ for various choices of $\mcal{S}$.

\begin{lemma}
\label{LEM:sparse-bound}
Assume $F$ is diagonal.
Suppose $\set{q\bm{c}: (q, \bm{c})\in \QQ^\times \times \mcal{S}} \cap \ZZ^m = \mcal{S}$.
Assume $m\leq 6$ and $\delta\leq \frac12(m-2)$.
Then $f(\mcal{S})\ll_\eps X^{(m-\delta)/2+\eps}$.
\end{lemma}


The work \cite{heath1998circle}*{pp.~688--689} proves a simplified version of Lemma~\ref{LEM:sparse-bound}
with $(m,\delta)\in\set{(4,1),(6,0)}$,
and with $\bm{1}_{n_\star=n}$ instead of $\sum_{n_\star\mid n}$.
The method of \cite{heath1998circle} should directly extend to Lemma~\ref{LEM:sparse-bound}.
However, we transpose Heath-Brown's argument a bit in order to highlight an intermediate result that one may hope to generalize: Proposition~\ref{PROP:S_c(n)-average-absolute-bound-over-n-for-fixed-lacunary-c}.
Proposition~\ref{PROP:S_c(n)-average-absolute-bound-over-n-for-fixed-lacunary-c} offers a potential partial alternative to axiom~(3) in Remark~\ref{RMK:non-diagonal-F}.

Let $v_p(-)$ denote the usual $p$-adic valuation.
For integers $c\neq 0$,
let $\map{sq}(c)\defeq \prod_{p^2\mid c} p^{v_p(c)}$ and $\map{cub}(c)\defeq \prod_{p^3\mid c} p^{v_p(c)}$; and for convenience, let $\map{sq}(0)\defeq 0$ and $\map{cub}(0)\defeq 0$.
A positive integer $n$ is said to be \emph{square-full} if $n = \map{sq}(n)$,
and \emph{cube-full} if $n = \map{cub}(n)$.
In the absence of a deeper algebro-geometric understanding of $S_{\bm{c}}$,
one relies heavily on the following bound of Hooley and Heath-Brown, valid for diagonal $F$ (for all $n\geq 1$ and $\bm{c}\in \ZZ^m$):
\begin{equation}
\label{INEQ:crude-pointwise-bound-on-diagonal-S_c(n)}
n^{-1/2} S^\natural_{\bm{c}}(n)
\ll_{F,\eps} n^\eps
\prod_{1\leq j\leq m} \gcd\bigl(\map{cub}(n)^2, \gcd(\map{cub}(n),\map{sq}(c_j))^3\bigr)^{1/12}.
\end{equation}
(See \cite{wang2023_large_sieve_diagonal_cubic_forms}*{Proposition~4.9} or \cite{wang2022thesis}*{Proposition~3.3.3} for precise references.)

Suppose $m\leq 6$, and assume $F$ is diagonal.
Let $C, N\in \set{1,2,4,8,\ldots}$.
For each $i\in [m]$, the $c_i^{\deg F^\vee}$ coefficient of the homogeneous polynomial $F^\vee\in \ZZ[\bm{c}]$ is nonzero, by \eqref{EQN:diagonal-discriminant-factorization}.
Suppose $\bm{c}\in \mcal{S}$ with $C\leq \norm{\bm{c}}<2C$;
then $\abs{c_i}\gg C$ for at least two indices $i\in [m]$.
In particular,
\begin{equation}
\label{INEQ:H-B-lacunary-pointwise-integral-estimate}
    I_{\bm{c}}(n) \ll_\eps X^{m+\eps}(XC/N)^{(2-m)/4},
\end{equation}
by \cite{heath1998circle}*{p.~688, (7.3)}.
We also have the following result (for $\bm{c}$, $N$ as above):

\begin{proposition}
\label{PROP:S_c(n)-average-absolute-bound-over-n-for-fixed-lacunary-c}
Here $\sum_{N\leq n<2N} \sum_{n_\star\mid n} n_\star^{-1/2} \abs{S^\natural_{\bm{c}}(n_\star)}\ll_\eps (CN)^\eps \map{cub}(\gcd(\bm{c}))^{1/6} N$.
\end{proposition}

\begin{proof}
Let $M_F\defeq \prod_{1\leq i\leq m} F_i$.
Write $\bm{c} = g\bm{c}'$ with $g\defeq \gcd(\bm{c})>0$, so that $\bm{c}'$ is primitive.
For each prime $p$,
the equation $F^\vee(\bm{c}') = 0$ implies, via \eqref{EQN:diagonal-discriminant-factorization}, that $\#\set{i\in [m]: v_p(c'_i)\leq v_p(M_F)}\geq 2$.
The bound \eqref{INEQ:crude-pointwise-bound-on-diagonal-S_c(n)} for $n_\star\mid n$ now implies
\begin{equation*}
\sum_{n_\star\mid n} n_\star^{-1/2}\abs{S^\natural_{\bm{c}}(n_\star)}
\ll_\eps n^\eps
\map{cub}(n)^{(m-2)/6}
\gcd(\map{cub}(n), \map{sq}(M_F\cdot g))^{2/4}
\end{equation*}
(since $\map{cub}(n_\star)\mid \map{cub}(n)$,
and $v_p(\map{sq}(c_i)) \leq v_p(\map{sq}(M_F\cdot g))$ whenever $v_p(c'_i)\leq v_p(M_F)$).
Thus
\begin{equation}
\label{INEQ:key-lacunary-reduction}
\sum_{N\leq n<2N}\, \sum_{n_\star\mid n} n_\star^{-1/2} \abs{S^\natural_{\bm{c}}(n_\star)}
\ll_\eps \sum_{d\mid \map{sq}(M_F\cdot g)}\,
\sum_{\substack{q<2N:\, q=\map{cub}(q)}} (N/q)\cdot q^{(m-2)/6} d^{1/2} \bm{1}_{d\mid q}.
\end{equation}
However, for any integer $d\geq 1$ and real $\eps>0$, the sum $\sum_{\substack{q<2N:\, q=\map{cub}(q)}} (2N)^{-\eps} q^{-1/3} d^{1/2} \bm{1}_{d\mid q}$ is
\begin{equation*}
\le \sum_{q\geq 1:\, q=\map{cub}(q)} \frac{d^{1/2} \bm{1}_{d\mid q}}{q^{1/3+\eps}}
= \biggl(\,\prod_{p\mid d} \frac{O(p^{v_p(d)/2})}{p^{\max(1, v_p(d)/3)}}\biggr)
\prod_{p\nmid d} (1 + O(p^{-1-3\eps}))
\ll_\eps d^\eps \map{cub}(d)^{1/6}.
\end{equation*}
Since $(m-2)/6\leq 2/3$, the right-hand side of \eqref{INEQ:key-lacunary-reduction} is therefore $\ll_\eps (CN)^\eps \map{cub}(M_F\cdot g)^{1/6} N$.
But $\map{cub}(M_F\cdot g)\mid M_F^3 \map{cub}(g)$, so Proposition~\ref{PROP:S_c(n)-average-absolute-bound-over-n-for-fixed-lacunary-c} follows.
\end{proof}

\begin{proof}
[Proof of Lemma~\ref{LEM:sparse-bound}]
By Proposition~\ref{PROP:standard-delta-method-n,c-cutoffs}, \eqref{INEQ:H-B-lacunary-pointwise-integral-estimate}, Proposition~\ref{PROP:S_c(n)-average-absolute-bound-over-n-for-fixed-lacunary-c}, and dyadic decomposition of $\norm{\bm{c}}, n$, the quantity $f(\mcal{S})$ (see \eqref{EQN:sparse-absolute-contribution-f(S)-to-delta-method}) is
\begin{equation*}
\ll_{\eps,A} X^{-A}
+ Y^{-2} X^{m+\eps} \sum_{\substack{C,N\in \set{1,2,4,8,\ldots}: \\ 1\leq C\leq X^{1/2+\eps},\; 1\leq N\leq M_2 Y}}
(XC/N)^{(2-m)/4} N^{2-m/2} \sum_{\substack{\bm{c}\in \mcal{S}: \\ C\leq \norm{\bm{c}}<2C}} \map{cub}(\gcd(\bm{c}))^{1/6}.
\end{equation*}
Our hypotheses on $\mcal{S}$ imply
\begin{equation}
\label{INEQ:bound-cub(gcd(c))-on-average-over-sparse-set-S}
\sum_{\bm{c}\in \mcal{S}:\, C\leq \norm{\bm{c}}<2C} \map{cub}(\gcd(\bm{c}))^{1/6}
\leq \sum_{g<2C} \map{cub}(g)^{1/6}\cdot (2C/g)^{m/2-\delta}.
\end{equation}
But the Dirichlet series $\sum_{g\geq 1} \map{cub}(g)^{1/6} g^{-s}$ converges absolutely for $\Re(s) > 1$.
Since $m/2-\delta\geq 1$, the right-hand side of \eqref{INEQ:bound-cub(gcd(c))-on-average-over-sparse-set-S} is therefore $\ll_\eps C^{m/2-\delta+\eps}$.
So
\begin{equation}
\label{INEQ:final-dyadic-bound-for-f(S)}
f(\mcal{S}) \ll_{\eps,A} X^{-A}
+ Y^{-2} X^{m+2\eps} \sum_{\substack{C,N\in \set{1,2,4,8,\ldots}: \\ 1\leq C\leq X^{1/2+\eps},\; 1\leq N\leq M_2 Y}}
(XC/N)^{(2-m)/4} N^{2-m/2} C^{m/2-\delta}.
\end{equation}

The total exponent of $N$ in \eqref{INEQ:final-dyadic-bound-for-f(S)} is $\frac{m-2}{4} + 2-\frac{m}{2} = \frac{6-m}{4}\geq 0$,
and the total exponent of $C$ in \eqref{INEQ:final-dyadic-bound-for-f(S)} is $\frac{2-m}{4} + \frac{m}{2}-\delta = \frac{m+2}{4}-\delta\geq 0$ (since $\frac{m+2}{4}\geq \frac{m-2}{2}$ for $m\leq 6$).
It follows that
\begin{equation*}
f(\mcal{S})\ll_\eps Y^{-2} X^{m+3\eps} (X^{3/2+\eps}/Y)^{(2-m)/4} Y^{2-m/2} (X^{1/2+\eps})^{m/2-\delta}.
\end{equation*}
Since $Y=X^{3/2}$, we get $f(\mcal{S})\ll_\eps X^{m+O(m\eps)} X^{-3m/4} X^{m/4-\delta/2} = X^{m/2-\delta/2+O(m\eps)}$.
\end{proof}

Now drop the earlier assumptions ``$m\le 6$'' and ``$F$ is diagonal''.
The axioms in Remark~\ref{RMK:non-diagonal-F} would allow us to prove a version of Lemma~\ref{LEM:sparse-bound} without assuming $F$ is diagonal.

\begin{lemma}
\label{LEM:axiomatic-bound-on-sparse-contribution-f(S)}
Assume \ref{RMK:non-diagonal-F}(2)--(3).
Say $\delta\leq \min(1, \frac23(10-m))$.
Then $f(\mcal{S})\ll_\eps X^{m/2-\delta/4+\eps}$.
\end{lemma}
\noindent When $m\le 9$, this beats the square-root threshold $X^{m/2}$, corresponding to linear subspaces.
(Some degenerate ranges of $\bm{c}$, $n$
seem to prevent us from handling $m\geq 10$.)

\begin{proof}
[Proof sketch for Lemma~\ref{LEM:axiomatic-bound-on-sparse-contribution-f(S)}]
Under \ref{RMK:non-diagonal-F}(2),
it is known that
\begin{equation}
\label{INEQ:Hessian-free-integral-bound}
I_{\bm{c}}(n)
\ll_\eps X^{m+\eps} (1 + X\norm{\bm{c}}/n)^{1-m/2};
\end{equation}
see e.g.~\cite{hooley2014octonary}*{p.~252, (31)}.
It is also known that if $n$ is cube-free, then
\begin{equation}
\label{INEQ:cube-free-bound}
    n^{-1/2} S^\natural_{\bm{c}}(n)\ll_\eps n^\eps;
\end{equation}
see e.g.~\cite{hooley2014octonary}*{Lemmas~8--9}.
Fix $\eps>0$.
For integers $C,N\geq 1$, let
\begin{equation*}
B(C,N)\defeq (C^{m/2-\delta})^{1/2} (C^{m/2} + N^{m/6})^{1/2}
+ (C^{(m-1)/2} N^{1/6} + N^{m/6}).
\end{equation*}
By \eqref{COND:first-moment-axiom-over-W=0}, plus Cauchy on \eqref{COND:second-moment-axiom-over-nonzero-W} over $\mcal{S}$, it follows that if $X$ is sufficiently large, then
\begin{equation*}
f(\mcal{S})
\ll_\eps Y^{-2} X^{m+\eps} (1 + XC/N)^{1-m/2} N^{1-m/2} (N/N_\star)
\cdot N_\star^{1/3} B(C,N_\star)
\end{equation*}
for some choice of $C, N, N_\star\in \set{1,2,4,8,\ldots}$ with $C\leq X^{1/2+\eps}$ and $N_\star\leq N\leq Y$.
Optimizing $N_\star$ over $1\leq N_\star\leq N$ yields
$f(\mcal{S})
\ll_\eps X^{m-3+\eps} (N+XC)^{1-m/2} (N\cdot B(C,1) + N^{1/3}\cdot B(C,N))$.
But $B(C,1)\ll C^{m/2-\delta/2}$ (since $\delta\leq 1$) and $B(C,N)\ll (C^{1/2}+N^{1/6})^m$, so
\begin{equation*}
f(\mcal{S})
\ll_\eps X^{m-3+\eps} (N+XC)^{1-m/2} (N C^{m/2-\delta/2} + N^{1/3} C^{m/2} + N^{1/3+m/6}).
\end{equation*}
Let $M\defeq N+XC$, and plug in the bounds $N\leq M$ and $C\leq M/X$, to get
\begin{equation*}
f(\mcal{S})/X^{m-3+\eps}
\ll_\eps M^{2-\delta/2}/X^{m/2-\delta/2}
+ M^{4/3}/X^{m/2}
+ M^{4/3-m/3}.
\end{equation*}
But $X\leq M\leq 2X^{3/2+\eps}$,
so $f(\mcal{S})\ll_\eps X^{3\eps} (X^{m/2-\delta/4} + X^{m/2-1} + X^{2m/3-5/3})$.
Here $\frac{m}{2}-1\leq \frac{m}{2}-\frac{\delta}{4}$ and $\frac{2m}{3}-\frac{5}{3}\leq \frac{m}{2}-\frac{\delta}{4}$, since $\delta\leq \min(1, \frac23(10-m))$.
\end{proof}


\section{Estimates at the center}
\label{SEC:c=0-singular-series}

In this section, we collect some standard facts we need about the quantities $I_{\bm{0}}(n)$, $S_{\bm{0}}(n)$ defined in \eqref{EQN:define-I_c(n)} and \eqref{EQN:define-S_c(n)}.
By \cite{heath1996new}*{Lemma~16}, we have
\begin{equation}
\label{INEQ:I_0(n)-bounded}
I_{\bm{0}}(n)\ll X^m.
\end{equation}
By \eqref{INEQ:I_0(n)-bounded} and \cite{heath1996new}*{Lemma~13}, we have (uniformly over $X,n\ge 1$)
\begin{equation}
\label{EQ:I0-approximation}
X^{-m} I_{\bm{0}}(n) = \sigma_{\infty,F,w} + O_A((n/Y)^A)
\end{equation}
(for all $A>0$),
where
\begin{equation}
\label{EQ:define-real-density}
\sigma_{\infty,F,w}
\defeq \lim_{\eps\to 0}{(2\eps)^{-1} \int_{\abs{F(\bm{x})}\leq\eps} d\bm{x}\, w(\bm{x})}
\ll_{F,w} 1.
\end{equation}
Aside from the real density $\sigma_{\infty,F,w}$, we also have the singular series
\begin{equation}
\label{EQ:define-singular-series}
\mathfrak{S}_F
\defeq \sum_{n\geq1} n^{-m}S_{\bm{0}}(n)
= \sum_{n\geq1} n^{(1-m)/2}S^\natural_{\bm{0}}(n),
\end{equation}
which converges absolutely for $m\geq 5$ (as one can show using Lemma~\ref{LEM:sum-S_0(n)-trivially}, for instance).

\begin{lemma}
\label{LEM:sum-S_0(n)-trivially}
Assume $m\ge 4$.
If $N\in \set{1,2,4,8,\ldots}$,
then
\begin{equation*}
\sum_{N\le n<2N} n^{1-m/2}
\sum_{d\mid n} d^{-1/2}\abs{S^\natural_{\bm{0}}(d)}
\ll_\eps N^{(4-m)/3+\eps}.
\end{equation*}
\end{lemma}

\begin{proof}
We have $S^\natural_{\bm{0}}(d)
\ll_\eps d^{1/2+\eps}\map{cub}(d)^{m/6}$ by \eqref{INEQ:crude-pointwise-bound-on-diagonal-S_c(n)} for diagonal $F$, and by \cite{hooley1988nonary}*{p.~95, (170)} in general.
Thus $\sum_{d\mid n} d^{-1/2}\abs{S^\natural_{\bm{0}}(d)}\ll_\eps n^\eps \map{cub}(n)^{m/6}$.
But
\begin{equation*}
\sum_{N\le n<2N} \map{cub}(n)^{m/6}
\ll \sum_{n_3<2N:\, n_3=\map{cub}(n_3)} n_3^{m/6}\cdot (N/n_3)
\ll_\eps N \cdot N^{m/6-2/3+\eps} \sum_{n_3=\map{cub}(n_3)} n_3^{-1/3-\eps},
\end{equation*}
since $m/6-2/3\ge 0$.
Observe that $\sum_{n_3=\map{cub}(n_3)} n_3^{-1/3-\eps} \ll_\eps 1$, and multiply by $N^{1-m/2+\eps}$.
\end{proof}


Since $I_{\bm{0}}(n)=0$ for $n\ge M_2Y$ (by Proposition~\ref{PROP:standard-delta-method-n,c-cutoffs}),
a routine calculation\footnote{with the same numerics as the diagonal treatment in \cite{wang2023_large_sieve_diagonal_cubic_forms}*{\S6}} using \eqref{EQ:I0-approximation}, Lemma~\ref{LEM:sum-S_0(n)-trivially}, and \eqref{EQ:define-singular-series} gives, for $m\geq 5$, the equality
\begin{equation}
\label{INEQ:isolate-singular-series-from-center}
Y^{-2}\sum_{n\geq1}
n^{-m}S_{\bm{0}}(n)I_{\bm{0}}(n)
= \sigma_{\infty,F,w}\mathfrak{S}_F X^{m-3}
+ O_\eps(X^{m/2-1+\eps}).
\end{equation}

\begin{remark}
The presence of $\sum_{d\mid n}$ in Lemma~\ref{LEM:sum-S_0(n)-trivially} is not important for \eqref{INEQ:isolate-singular-series-from-center} (where $\bm{1}_{d=n}$ would suffice in place of $\sum_{d\mid n}$),
but rather for Lemma~\ref{LEM:handle-n_1>=Y/P-for-extra-c's}
(to help separate $\bm{c}=\bm{0}$ from $\bm{c}\neq \bm{0}$).
\end{remark}

\section{Establishing bias in exponential sums}
\label{SEC:proving-bias-in-linear-S_c's}

In this section, we will realize \eqref{EQN:illustrate-bias-philosophy} from \S\ref{SEC:motivating-background-and-intro-to-structured-part-detection}.
We will start with general theory, and gradually impose restrictions.
It will be convenient to work over the $p$-adic integers, $\ZZ_p$.

Let $p$ be a prime.
Certain hyperplane sections govern the behavior of $S_{\bm{c}}(p^l)$ for $l\ge 1$, as $\bm{c}$ ranges over $\ZZ^m$ or more generally, $\ZZ_p^m$.\footnote{The formula for $S_{\bm{c}}(p^l)$ in \eqref{EQN:define-S_c(n)} still makes sense for $\bm{c}\in \ZZ_p^m$, because $\bm{c}\bmod{p^l}\in \ZZ/p^l\ZZ$.}
Let $\mcal{V}$ and $\mcal{V}_{\bm{c}}$ denote the closed subschemes of $\PP^{m-1}_{\ZZ_p}$
defined by the equations $F(\bm{x}) = 0$ and $F(\bm{x}) = \bm{c}\cdot\bm{x} = 0$, respectively.
Throughout the present \S\ref{SEC:proving-bias-in-linear-S_c's} only,
let $V\defeq\mcal{V}_{\FF_p}$ and $V_{\bm{c}}\defeq(\mcal{V}_{\bm{c}})_{\FF_p}$,
so that $V$ and $V_{\bm{c}}$
live over $\FF_p$ (not $\QQ$).

We now recall some standard background recorded (with references) in \cite{wang2022thesis}*{\S3.2}.
Let $m_\ast\defeq m-3$.
In terms of the point counts $\card{V(\FF_p)}$ and $\card{V_{\bm{c}}(\FF_p)}$ over $\FF_p$, let
\begin{equation*}
E(p)\defeq \card{V(\FF_p)} - (p^{m-1}-1)/(p-1),
\quad
E_{\bm{c}}(p)\defeq \card{V_{\bm{c}}(\FF_p)} - (p^{m-2}-1)/(p-1).
\end{equation*}
Then let $E^\natural(p) \defeq p^{-(m_\ast+1)/2}E(p)$
and $E^\natural_{\bm{c}}(p) \defeq p^{-m_\ast/2}E_{\bm{c}}(p)$.
By Deligne's resolution of the Weil conjectures,
\begin{equation}
\label{INEQ:E(p)-Weil-conjectures-bound}
    E^\natural(p) \ll_F 1.
\end{equation}
Furthermore,
whenever $p\nmid\bm{c}$, we have $S_{\bm{c}}(p) = p^2 E_{\bm{c}}(p) - p E(p)$, or equivalently,
\begin{equation}
\label{EQN:rewrite-S_c(p)-via-E_c}
S^\natural_{\bm{c}}(p) = E^\natural_{\bm{c}}(p) - p^{-1/2}E^\natural(p).
\end{equation}

For prime powers $p^l$ with $l\geq 2$, a different flavor of geometry, based on Hensel lifting, comes into play; cf.~\cite{hooley1986HasseWeil}*{pp.~65--66}.
For the sake of other work (\cite{wang2023_HLH_vs_RMT}), we prove more than we presently need;
we encourage the reader to skip ahead to Corollary~\ref{COR:primitive-case-projective-balanced-rewrite-of-S_c} on a first reading.
We start by identifying a clean source of cancellation in \eqref{EQN:define-S_c(n)}.
For each integer $d\ge 0$, let
\begin{equation}
\label{EQN:define-script-S(c,p^d)-level-p^d-critical-locus}
\mscr{S}(\bm{c},p^d)\defeq \bigcup_{\lambda\in \ZZ_p^\times} \set{\bm{x}\in \ZZ_p^m:
p^d\mid \grad{F}(\bm{x})-\lambda\bm{c}}.
\end{equation}

\begin{lemma}
\label{LEM:prelim-reduction-to-halfway-gradient-proportionality}
Let $\bm{c},\bm{x}_0\in\ZZ_p^m$.
Let $l\geq 2$ and $d\in [0, (l+\bm{1}_{p\mid\bm{x}_0})/2]$.
Then $d\leq l-1$, and
\begin{equation}
\label{EQN:basic-cancellation-lemma-goal}
\sum_{1\le a\le p^l:\, p\nmid a}\,
\sum_{1\le \bm{x}\le p^l}
\bm{1}_{\bm{x}\equiv\bm{x}_0\bmod{p^{l-d}}}
\cdot \bm{1}_{\bm{x}\notin \mscr{S}(\bm{c},p^d)}
\cdot e_{p^l}(aF(\bm{x}) + \bm{c}\cdot\bm{x})
= 0.
\end{equation}
\end{lemma}

\begin{proof}
If $l=2$, then $d\leq 1\leq 2l/3$;
if $l\geq3$, then $d\leq (l+1)/2\leq 2l/3$.
Thus $d\le \floor{2l/3}\le l-1$.

On the left-hand side of \eqref{EQN:basic-cancellation-lemma-goal}, write $\bm{x}=\bm{x}_0+p^{l-d}\bm{r}$ (with $\bm{r}\in \ZZ_p^m$ running over a complete set of residues modulo $p^d$).
Trivially, $\bm{c}\cdot\bm{x}
\equiv \bm{c}\cdot\bm{x}_0
+ \bm{c}\cdot p^{l-d}\bm{r}\bmod{p^l}$.
Also, since $F$ is homogeneous of degree $3$, Taylor expansion (using $\min(\bm{1}_{p\mid\bm{x}_0} + 2(l-d), 3(l-d))\geq l$) gives
\begin{equation*}
F(\bm{x})
\equiv F(\bm{x}_0)
+ \grad{F}(\bm{x}_0)
\cdot p^{l-d}\bm{r}\bmod{p^l}.
\end{equation*}
Furthermore, $\grad{F}(\bm{x})\equiv\grad{F}(\bm{x}_0)\bmod{p^d}$
(since $\grad{F}$ is ``homogeneous of degree $2$'',
and $\min(\bm{1}_{p\mid\bm{x}_0}
+ (l-d),
2(l-d))\geq d$),
whence $\bm{1}_{\bm{x}\notin \mscr{S}(\bm{c},p^d)} = \bm{1}_{\bm{x}_0\notin \mscr{S}(\bm{c},p^d)}$ (by \eqref{EQN:define-script-S(c,p^d)-level-p^d-critical-locus}).

If $\bm{x}_0\in \mscr{S}(\bm{c},p^d)$, then the left-hand side of \eqref{EQN:basic-cancellation-lemma-goal} directly vanishes.
Now suppose $\bm{x}_0\notin \mscr{S}(\bm{c},p^d)$.
Then for each $a$ in \eqref{EQN:basic-cancellation-lemma-goal}, we have $p^d\nmid a\grad{F}(\bm{x}_0) + \bm{c}$, whence the sum of $e_{p^l}(aF(\bm{x}) + \bm{c}\cdot\bm{x})$
over $\bm{x}\equiv\bm{x}_0\bmod{p^{l-d}}$ (i.e.~over $\bm{r}\bmod{p^d}$) vanishes.
Summing over $a$ gives \eqref{EQN:basic-cancellation-lemma-goal}.
\end{proof}

We now show that to understand $S_{\bm{c}}(p^l)$ for $l\ge 2$, it suffices to understand
\begin{equation}
\label{EQN:define-primitive-x-sum-S'_c}
S'_{\bm{c}}(p^l)
\defeq \sum_{1\le a\le p^l:\, p\nmid a}\,
\sum_{1\le \bm{x}\le p^l:\, p\nmid \bm{x}}
e_{p^l}(aF(\bm{x}) + \bm{c}\cdot\bm{x}).
\end{equation}
Before proceeding, note that by \eqref{EQN:define-S_c(n)}, \eqref{EQN:define-primitive-x-sum-S'_c}, we have
\begin{equation}
\label{EQN:decompose-S_c-as-primitive-x-plus-non-primitive-x}
S_{\bm{c}}(p^l)-S'_{\bm{c}}(p^l)
= \sum_{1\le a\le p^l:\, p\nmid a}\,
\sum_{1\le \bm{x}\le p^l:\, p\mid \bm{x}}
e_{p^l}(aF(\bm{x}) + \bm{c}\cdot\bm{x}).
\end{equation}

\begin{lemma}
\label{LEM:reduction-to-primitive-x's}
Fix a tuple $\bm{c}\in\ZZ_p^m$ and an integer $l\geq 2$.
Then $S_{\bm{c}}(p^l)-S'_{\bm{c}}(p^l)$
equals
\begin{enumerate}
    \item $\bm{1}_{p\mid\bm{c}}
    \cdot \phi(p^2)p^m$ if $l=2$,
    and
    
    \item $\bm{1}_{p^2\mid\bm{c}}
    \cdot [\phi(p^l)/\phi(p^{l-3})]
    \cdot p^{2m}
    \cdot S_{\bm{c}/p^2}(p^{l-3})$ if $l\geq 3$.
\end{enumerate}
In particular,
if $p\nmid\bm{c}$,
then $S_{\bm{c}}(p^l)=S'_{\bm{c}}(p^l)$.
\end{lemma}

\begin{proof}
Let $d\defeq \min(2, \floor{(l+1)/2})$;
then Lemma~\ref{LEM:prelim-reduction-to-halfway-gradient-proportionality} applies whenever $p\mid \bm{x}_0$.
Summing \eqref{EQN:basic-cancellation-lemma-goal} over $\set{1\le \bm{x}_0\le p^{l-d}: p\mid \bm{x}_0}$, we get (by \eqref{EQN:decompose-S_c-as-primitive-x-plus-non-primitive-x})
\begin{equation*}
S_{\bm{c}}(p^l)-S'_{\bm{c}}(p^l)
= \sum_{1\le a\le p^l:\, p\nmid a}\,
\sum_{1\le \bm{x}\le p^l:\, p\mid \bm{x}}
\bm{1}_{\bm{x}\in \mscr{S}(\bm{c},p^d)}
\cdot e_{p^l}(aF(\bm{x}) + \bm{c}\cdot\bm{x}).
\end{equation*}
Here $p^2\mid\grad{F}(\bm{x})$,
and $d\leq 2$,
so $\bm{1}_{\bm{x}\in \mscr{S}(\bm{c},p^d)} = \bm{1}_{p^d\mid\bm{c}}$ (by \eqref{EQN:define-script-S(c,p^d)-level-p^d-critical-locus}).
So if $p^d\nmid \bm{c}$, then $S_{\bm{c}}(p^l)-S'_{\bm{c}}(p^l) = 0$, which suffices (since $d=1$ if $l=2$, and $d=2$ if $l\ge 3$).
Now suppose $p^d\mid \bm{c}$.

If $l=2$,
then $p^l\mid F(\bm{x}),\bm{c}\cdot\bm{x}$,
so $S_{\bm{c}}(p^l)-S'_{\bm{c}}(p^l)
= \phi(p^l) \sum_{1\le \bm{x}\le p^l:\, p\mid \bm{x}} 1
= \phi(p^2)p^m$.

Now suppose $l\ge 3$.
Then $d=2$, and $\bm{c}'\defeq\bm{c}/p^2\in\ZZ_p^m$.
Now write $\bm{x}=p\bm{x}'$ to get
\begin{equation*}
S_{\bm{c}}(p^l)-S'_{\bm{c}}(p^l)
= \sum_{1\le a\le p^l:\, p\nmid a}\,
\sum_{1\le \bm{x}'\le p^{l-1}}
e_{p^l}(ap^3F(\bm{x}') + p^3\bm{c}'\cdot\bm{x}').
\end{equation*}
Now let $l'\defeq l-3\ge 0$.
Then $e_{p^l}(-)$ is determined by $(a,\bm{x}')\bmod{p^{l'}}$.
Therefore
\begin{equation*}
S_{\bm{c}}(p^l)-S'_{\bm{c}}(p^l)
= [\phi(p^l)/\phi(p^{l'})]
\cdot p^{2m}
\cdot \sum_{1\le a\le p^{l'}:\, p\nmid a}\,
\sum_{1\le \bm{x}'\le p^{l'}}
e_{p^{l'}}(aF(\bm{x}') + \bm{c}'\cdot\bm{x}'),
\end{equation*}
which equals $[\phi(p^l)/\phi(p^{l'})]
\cdot p^{2m}\cdot S_{\bm{c}'}(p^{l'})$ by \eqref{EQN:define-S_c(n)}.
This completes the proof.
\end{proof}

The general study of $S'_{\bm{c}}(p^l)$ needs some setup.
For any vector $\bm{b}\in \ZZ_p^m$, let $v_p(\bm{b})\defeq v_p(\gcd(b_1,\dots,b_m))\in [0,\infty]$.
Given $\bm{c}\in \ZZ_p^m\setminus \set{\bm{0}}$ and integers $r,s,d\ge 0$, let
\begin{equation}
\begin{split}
\label{EQN:define-S_r,s,d,S^ast_r,s,d}
    \mcal{S}_{r,s}(\bm{c},d)
    &\defeq \set{\bm{x}\in \ZZ_p^m\cap \mscr{S}(\bm{c},p^d):
    p\nmid \bm{x},\; p^r\mid F(\bm{x}),\; p^{s+v_p(\bm{c})}\mid \bm{c}\cdot \bm{x}}, \\
    \mcal{S}^\ast_{r,s}(\bm{c},d)
    &\defeq \set{\bm{x}\in \mcal{S}_{r,s}(\bm{c},d): v_p(\grad{F}(\bm{x})) = v_p(\bm{c})}.
\end{split}
\end{equation}
Let $\mu_p$ denote the usual Haar measure on $\ZZ_p^m$, so that for all $l\ge \max(1,r,s,d)$,
we have
\begin{equation}
\label{EQN:trivial-counting-interpretation-of-S_r,s(c,d)-measure}
\mu_p(\mcal{S}_{r,s}(\bm{c},d)) = p^{-lm} \cdot
\card{\set{\bm{x}\in \mscr{S}(\bm{c},p^d):
1\le \bm{x}\le p^l,\; p\nmid \bm{x},\; p^r\mid F(\bm{x}),\; p^{s+v_p(\bm{c})}\mid \bm{c}\cdot \bm{x}}}.
\end{equation}

\begin{lemma}
Suppose $\bm{c}\in \ZZ_p^m\setminus \set{\bm{0}}$ and $g=v_p(\bm{c})$.
Let $u, v\ge 1$ and $d\ge g$ be integers.
Then
\begin{align}
\mu_p(\mcal{S}_{u,v}(\bm{c},d))
&= p^{u-v}\cdot \mu_p(\mcal{S}_{u,u}(\bm{c},d)),
\textnormal{ if }
v\ge d\textnormal{ and }d\le u\le g+v;
\label{EQN:equalize-if-d<=u<=g+v} \\
\mu_p(\mcal{S}^\ast_{u,v}(\bm{c},d))
&= p^{-1}\cdot \mu_p(\mcal{S}^\ast_{u-1,v}(\bm{c},d)),
\textnormal{ if }
u-1-g \ge \max(1+g, d, v).
\label{EQN:decrease-u-if-u>=1+g+v}
\end{align}
\end{lemma}

\begin{proof}
Let $r,s\ge 1$ be integers.
Using \eqref{EQN:define-S_r,s,d,S^ast_r,s,d} and the congruence
\begin{equation*}
    F(\bm{x}+p^{\max(d,r-g)}\bm{h})\equiv F(\bm{x}) + \grad{F}(\bm{x})\cdot p^{\max(d,r-g)}\bm{h} \bmod{p^r}
\end{equation*}
(valid since $2\max(d,r-g)\ge 2\max(g,r-g)\ge r$), we find that
\begin{enumerate}
    \item $\mcal{S}_{r,s}(\bm{c},d)$ is invariant under addition by any element of $p^{\max(d, r-g, s)}\ZZ_p^m$ (since $p^{g}\mid \grad{F}(\bm{x})$ for all $\bm{x}\in \mcal{S}_{r,s}(\bm{c},d)$, by \eqref{EQN:define-script-S(c,p^d)-level-p^d-critical-locus}), and therefore

    \item $\mcal{S}^\ast_{r,s}(\bm{c},d)$ is invariant under $p^{\max(1+g, d, r-g, s)}\ZZ_p^m$.
\end{enumerate}

\emph{Case~1: $v\ge d$ and $d\le u\le g+v$.}
Then $\max(d,u-g,v)=v$ and $\max(d,u-g,u)=u$.
\begin{itemize}
    \item If $u\le v$, then the inclusion $\mcal{S}_{u,v}(\bm{c},d)\to \mcal{S}_{u,u}(\bm{c},d)$ descends to a map
    $\mcal{S}_{u,v}(\bm{c},d)/p^{v}\ZZ_p^m\to \mcal{S}_{u,u}(\bm{c},d)/p^u\ZZ_p^m$;
    and this map has fibers of size $p^{(m-1)(v-u)}$, so \eqref{EQN:equalize-if-d<=u<=g+v} holds.
    
    \item If $u\ge v$, then the inclusion $\mcal{S}_{u,u}(\bm{c},d)\to \mcal{S}_{u,v}(\bm{c},d)$ descends to a map
    $\mcal{S}_{u,u}(\bm{c},d)/p^u\ZZ_p^m\to \mcal{S}_{u,v}(\bm{c},d)/p^{v}\ZZ_p^m$;
    and this map has fibers of size $p^{(m-1)(u-v)}$, so \eqref{EQN:equalize-if-d<=u<=g+v} holds.
\end{itemize}

\emph{Case~2: $u-1-g \ge \max(1+g, d, v)$.}
Then $u\ge 2$, and $\max(1+g,d,r-g,v) = r-g$ for all $r\ge u-1$.
So the inclusion $\mcal{S}^\ast_{u,v}(\bm{c},d)\to \mcal{S}^\ast_{u-1,v}(\bm{c},d)$ descends to a map
\begin{equation*}
    \mcal{S}^\ast_{u,v}(\bm{c},d)/p^{u-g}\ZZ_p^m\to \mcal{S}^\ast_{u-1,v}(\bm{c},d)/p^{u-g-1}\ZZ_p^m.
\end{equation*}
This map has fibers of size $p^{m-1}$, since for all $\bm{x}_0\in \mcal{S}^\ast_{u-1,v}(\bm{c},d)$ and $\bm{r}\in \ZZ_p^m$, we have $v_p(\grad{F}(\bm{x}_0)) = g$ and the ``lifting congruence''
\begin{equation*}
    F(\bm{x}_0+p^{u-g-1}\bm{r}) \equiv F(\bm{x}_0) + \grad{F}(\bm{x}_0)\cdot p^{u-g-1}\bm{r}\bmod{p^{u}}.
\end{equation*}
(This congruence holds because $2(u-g-1)\ge u$.)
Thus \eqref{EQN:decrease-u-if-u>=1+g+v} holds.
\end{proof}

The next result synthesizes a lot of old and new Hensel work.

\begin{proposition}
\label{PROP:general-affine-balanced-rewrite-of-S'_c}
Let $\bm{c}\in \ZZ_p^m\setminus \set{\bm{0}}$, and let $g = v_p(\bm{c})$.
Let $d\ge g$ if $V$ is smooth, and let $d\ge 1+g$ if $V$ is singular.
Let $l\ge \max(2+2g, 2d)$.
Then
\begin{equation}
\label{EQN:affine-balanced-rewrite-S'_c-goal}
p^{-lm}\phi(p^l)S'_{\bm{c}}(p^l)
= p^{2l+g} \mu_p(\mcal{S}_{l,l}(\bm{c},d))
- p^{2l-2+g} \mu_p(\mcal{S}_{l-1,l-1}(\bm{c},d)).
\end{equation}
\end{proposition}

\begin{proof}
Lemma~\ref{LEM:prelim-reduction-to-halfway-gradient-proportionality} applies, since $d\le l/2$.
Summing \eqref{EQN:basic-cancellation-lemma-goal} over $\set{1\le \bm{x}_0\le p^{l-d}: p\nmid \bm{x}_0}$ gives
\begin{equation*}
S'_{\bm{c}}(p^l)
= \sum_{1\le a\le p^l:\, p\nmid a}\,
\sum_{1\le \bm{x}\le p^l:\, p\nmid \bm{x}}
\bm{1}_{\bm{x}\in \mscr{S}(\bm{c},p^d)}
\cdot e_{p^l}(aF(\bm{x}) + \bm{c}\cdot\bm{x}).
\end{equation*}
Replacing $\bm{c}$ with $\lambda\bm{c}$ for $\lambda\in \ZZ_p^\times$, and summing over $\lambda$, we get
(via the scalar symmetries $S'_{\bm{c}}(p^l) = S'_{\lambda \bm{c}}(p^l)$ and $\mscr{S}(\bm{c},p^d) = \mscr{S}(\lambda\bm{c},p^d)$ that follow from \eqref{EQN:define-primitive-x-sum-S'_c} and \eqref{EQN:define-script-S(c,p^d)-level-p^d-critical-locus}, respectively)
\begin{equation}
\label{EQN:simplify-S'_c-by-averaging-over-1-dim-scalar-subspace}
    p^{-lm}\phi(p^l)S'_{\bm{c}}(p^l)
    = \sum_{1\le \lambda\le p^l:\, p\nmid \lambda} p^{-lm}S'_{\lambda \bm{c}}(p^l)
    = \sum_{l-1\le u,v\le l} (-p)^{u+v}\mu_p(\mcal{S}_{u,v-g}(\bm{c},d)),
\end{equation}
by a short calculation using $\sum_{1\le b\le p^l:\, p\nmid b} = \sum_{1\le b\le p^l} - \sum_{1\le b\le p^l:\, p\mid b}$ (for $b=a,\lambda$) and \eqref{EQN:trivial-counting-interpretation-of-S_r,s(c,d)-measure}.

Let $r,s\ge 1$.
Before proceeding, we prove (by casework) that
\begin{equation}
\label{EQN:S^ast-equals-S-under-mild-conditions}
\mcal{S}^\ast_{r,s}(\bm{c},d) = \mcal{S}_{r,s}(\bm{c},d).
\end{equation}

\emph{Case~1: $d=g\ge 1$.}
Then $V$ is smooth, so $\set{\bm{x}\in \ZZ_p^m: p\nmid \bm{x},\; p\mid F(\bm{x}),\; p\mid \grad{F}(\bm{x})} = \emptyset$.
Since $p\mid \bm{c}$, we conclude by \eqref{EQN:define-script-S(c,p^d)-level-p^d-critical-locus}, \eqref{EQN:define-S_r,s,d,S^ast_r,s,d} that $\mcal{S}_{1,1}(\bm{c},1) = \emptyset$.
So both sides of \eqref{EQN:S^ast-equals-S-under-mild-conditions} are empty.

\emph{Case~2: $d=g=0$.}
Then $V$ is smooth, and $p\nmid \bm{c}$.
So by \eqref{EQN:define-S_r,s,d,S^ast_r,s,d}, we have \eqref{EQN:S^ast-equals-S-under-mild-conditions}.

\emph{Case~3: $d\ge g+1$.}
Then $v_p(\grad{F}(\bm{x})) = g$ for all $\bm{x}\in \mscr{S}(\bm{c}, p^d)$.
So by \eqref{EQN:define-S_r,s,d,S^ast_r,s,d}, we have \eqref{EQN:S^ast-equals-S-under-mild-conditions}.

Having established \eqref{EQN:S^ast-equals-S-under-mild-conditions} in all cases, we now return to \eqref{EQN:simplify-S'_c-by-averaging-over-1-dim-scalar-subspace}.
Since $d\ge g$ and $l\ge \max(2+2g, 1+g+d, l)$, we may apply \eqref{EQN:decrease-u-if-u>=1+g+v} (with $u=l$ and $v=l-1-g$) and \eqref{EQN:S^ast-equals-S-under-mild-conditions} to get
\begin{equation*}
\mu_p(\mcal{S}_{l,l-1-g}(\bm{c},d)) = p^{-1}\cdot \mu_p(\mcal{S}_{l-1,l-1-g}(\bm{c},d)).
\end{equation*}
But since $d\ge g$ and $l-1-g\ge \max(1, d)$, we may use \eqref{EQN:equalize-if-d<=u<=g+v} to get
\begin{equation*}
\mu_p(\mcal{S}_{u,v}(\bm{c},d)) = p^{u-v}\cdot \mu_p(\mcal{S}_{u,u}(\bm{c},d))
\end{equation*}
for all $(u, v) \in \set{l-1,l} \times \set{l-1-g, l-g}$ such that $u\le g+v$.
Thus the right-hand side of \eqref{EQN:simplify-S'_c-by-averaging-over-1-dim-scalar-subspace} equals $p^{l+l+g}\cdot \mu_p(\mcal{S}_{l,l}(\bm{c},d)) - (p^{l+(l-1)-1+g} + p^{(l-1)+l+(g-1)} - p^{(l-1)+(l-1)+g}) \cdot \mu_p(\mcal{S}_{l-1,l-1}(\bm{c},d))$, which simplifies to the right-hand side of \eqref{EQN:affine-balanced-rewrite-S'_c-goal}.
So \eqref{EQN:simplify-S'_c-by-averaging-over-1-dim-scalar-subspace} implies \eqref{EQN:affine-balanced-rewrite-S'_c-goal}.
\end{proof}

\begin{corollary}
\label{COR:primitive-case-projective-balanced-rewrite-of-S_c}
Suppose $\bm{c}\in \ZZ_p^m\setminus \set{\bm{0}}$ is primitive and $V$ is smooth.
Let $l\ge 2$.
Then
\begin{equation*}
S_{\bm{c}}(p^l)
= p^{2l} \card{\mcal{V}_{\bm{c}}(\ZZ/p^l\ZZ)}
- p^{2l+m_\ast} \card{\mcal{V}_{\bm{c}}(\ZZ/p^{l-1}\ZZ)}.
\end{equation*}
\end{corollary}

\begin{proof}
Here $v_p(\bm{c})=0$, so \eqref{EQN:trivial-counting-interpretation-of-S_r,s(c,d)-measure} implies $\mu_p(\mcal{S}_{u,u}(\bm{c},0)) = p^{-um}\phi(p^u)\card{\mcal{V}_{\bm{c}}(\ZZ/p^u\ZZ)}$ for $u\ge 1$.
Plug this into Proposition~\ref{PROP:general-affine-balanced-rewrite-of-S'_c} (with $d=g=0$);
then note that $S_{\bm{c}}(p^l) = S'_{\bm{c}}(p^l)$ by Lemma~\ref{LEM:reduction-to-primitive-x's}.
\end{proof}

Now fix $L\in \Upsilon$.
We will build up to Lemma~\ref{LEM:stating-bias-for-generic-linear-c's} (realizing \eqref{EQN:illustrate-bias-philosophy} from \S\ref{SEC:motivating-background-and-intro-to-structured-part-detection}).
As Remark~\ref{RMK:non-diagonal-F} suggests,
Lemma~\ref{LEM:stating-bias-for-generic-linear-c's}
might extend to more general $F$.
But to maximize the accessibility of \S\ref{SEC:proving-bias-in-linear-S_c's},
we focus on the diagonal case.
We will use an ad hoc change of coordinates,
highlighting specific features (of diagonal forms)
that may be of independent interest.

So for the rest of \S\ref{SEC:proving-bias-in-linear-S_c's}, assume $F$ is diagonal with $m\in \set{4,6}$.
Let $\mcal{J}$ be a permissible pairing corresponding to $L$ in Proposition~\ref{PROP:characterize-diagonal-L/Q}.
Since $\mcal{J}$ is permissible, there exist unique cube-free integers $F_{(k)}$ (for $k\in \mcal{K}$) such that $F_i/F_{(k)}$ is an integer cube for all $k\in \mcal{K}$ and $i\in \mcal{J}(k)$.

Suppose $\bm{c}$ lies in $\Lambda^\perp$ or more generally, $\Lambda^\perp \otimes \ZZ_p$.
Assume
\begin{equation}
\label{COND:jet-is-non-degenerate-mod-p}
    p\nmid j^{2^{m/2-1}}{F^\vee}(\bm{c}).
\end{equation}
In particular,
by \eqref{INEQ:diagonal-Fdual-primitive-scaling-convention},
we have
\begin{equation}
\label{COND:p-not-divides-6F_1...F_m}
    p\nmid (6^m)!F_1\cdots F_m.
\end{equation}

\begin{proposition}
\label{PROP:nice-consequences-of-non-degenerate-diagonal-jet}
Under \eqref{COND:jet-is-non-degenerate-mod-p} and \eqref{COND:p-not-divides-6F_1...F_m}, the following hold:
\begin{enumerate}
    \item Each $c(k)^3$ lies in $\ZZ_p$.
    
    \item If $k_1<\dots<k_t$,
    then $p\nmid\prod(c(k_1)^{3/2}\pm\cdots\pm c(k_t)^{3/2})$.
    In particular,
    $p\nmid c(k)^3$, hence $c(k)^3\in\ZZ_p^\times$,
    for each $k$.
    Also,
    $p\nmid c(i)^3-c(j)^3$ when $i\neq j$.
    
    \item $V_{\bm{c}}$ has exactly $2^{m/2-1}$ singular $\ol{\FF}_p$-points.
\end{enumerate}
\end{proposition}

\begin{proof}
(1):
By Definition~\ref{DEFN:diagonal-permissible-pairing},
$c(k)^3 = c_i^3/F_i$ for all $i\in\mcal{J}(k)$.
Here $p\nmid F_i$ by \eqref{COND:p-not-divides-6F_1...F_m}.

Next,
we use some results of \S\ref{SUBSEC:diagonal-dual-example},
carried over from $\QQ$ to $\FF_p$ via Remark~\ref{RMK:char-p-diagonal-Disc-analysis}.

(2):
Use \eqref{COND:jet-is-non-degenerate-mod-p} and the $\FF_p$-analog of Observation~\ref{OBS:higher-order-diagonal-Disc-vanishing-critera}(1).

(3):
$p\mid j^{2^{m/2-1}-1}{F^\vee}(\bm{c})$ by Corollary~\ref{COR:baseline-diagonal-Disc-vanishing},
carried over to $\FF_p$.
But $p\nmid j^{2^{m/2-1}}{F^\vee}(\bm{c})$ by \eqref{COND:jet-is-non-degenerate-mod-p}.
So by the $\FF_p$-analogs of Proposition~\ref{PROP:characterize-diagonal-Disc-vanishing-order} and \eqref{EQN:branching-interpretation-of-number-of-vanishing-L_eps}, the scheme $V_{\bm{c}}$ has
at least, but also at most, $2^{m/2-1}$ singular $\ol{\FF}_p$-points.
(Note that each singular $\ol{\FF}_p$-point of $V_{\bm{c}}$ has all coordinates nonzero, by the Jacobian criterion, since $p\nmid c_1\cdots c_m$ by (2).)
\end{proof}


Assume, until further notice, that
\begin{equation}
\label{EQN:convenient-diagonal-change-of-variables-assumption}
    F_i = F_{(k)}
\end{equation}
for all $k\in \mcal{K}$ and $i\in\mcal{J}(k)$.
For convenience, assume $\mcal{J}(k) = \set{k, k+m/2}$ for each $k\in \mcal{K}=[m/2]$.
Then $\bm{c}\in \Lambda^\perp \otimes \ZZ_p$ implies $c_k=c_{k+m/2}$.
Let $c^\star_k\defeq c_k=c_{k+m/2}$, so that
\begin{equation}
\label{EQN:rewrite-c(k)}
    c(k)^3 = (c^\star_k)^3/F_{(k)}.
\end{equation}

Now consider the equations $F(\bm{x}) = 0$ and $\bm{c}\cdot\bm{x} = 0$ defining $\mcal{V}_{\bm{c}}$.
These equations become
\begin{equation}
\label{EQN:V_c-equations-in-h,y-coordinates}
\sum_{1\le k\le m/2} F_{(k)}\cdot h[k]y[k]^2
= -3\sum_{1\le k\le m/2} F_{(k)}\cdot h[k]^3
\qquad\textnormal{and}\qquad
\sum_{1\le k\le m/2} c^\star_k\cdot h[k] = 0
\end{equation}
after a linear change of variables over $\ZZ[1/6]$.
Explicitly,
if $\mcal{J}(k) = \set{i,j}$ with $i<j$,
then we take $h[k]\defeq x_i+x_j$ and $y[k]\defeq 3(x_i-x_j)$,
so that the equation $h[1]=\cdots=h[m/2]=0$ cuts out $\Lambda \otimes \ZZ_p$.
(We use the letter ``$h$'' in analogy with
van der Corput or Weyl differencing.
The definition of $\bm{h}\defeq(h[k])_{1\le k\le m/2}$ is compatible with that in \S\ref{SEC:reverse-sums-by-duality}.)

Geometrically, over $K\defeq\FF_p$,
the space $\set{[\bm{h}]\in\PP^{m/2-1}:\bm{c}^\star\cdot\bm{h}=0}\cong\PP^{m/2-2}$
parameterizes projective $\frac{m}{2}$-planes $\PP{H}\belongs\PP\bm{c}^\perp\cong\PP^{m-2}$ containing the fixed $(\frac{m}{2}-1)$-plane $\PP\Lambda_K$.
Over this $[\bm{h}]$-space, we have the ``quadratic fibration''
(related to the blow-up of $V_{\bm{c}}$ along $\PP\Lambda_K$)
\begin{equation*}
V_{\bm{c}}\setminus\PP\Lambda_K
\to
\PP[\bm{c}^\star]^\perp\cong\PP^{m/2-2},
\quad
[\bm{x}]\mapsto[\bm{h}].
\end{equation*}
Concretely,
each slice $V_{\bm{c}}\cap\PP{H}$
consists of $\PP\Lambda_K$ and
a (possibly singular) quadric hypersurface $Q_H\belongs\PP{H}$ of dimension $m/2-1$,
where $Q_H\setminus\PP\Lambda_K$ is the fiber of $V_{\bm{c}}\setminus\PP\Lambda_K$ over $\PP{H}$.

Below, let $(\frac{r}{p})$ denote the Legendre symbol (if $p$ is odd), and write $\chi(r)\defeq (\frac{r}{p})$.


\begin{lemma}
\label{LEM:point-count-bias-mod-p}
Under \eqref{COND:jet-is-non-degenerate-mod-p}, \eqref{COND:p-not-divides-6F_1...F_m}, and \eqref{EQN:convenient-diagonal-change-of-variables-assumption},
we have $E^\natural_{\bm{c}}(p) = p^{1/2} + O(1)$.
\end{lemma}

\begin{proof}
Let $C(V_{\bm{c}})$ denote the affine cone of $V_{\bm{c}}$:
the subscheme $F(\bm{x}) = \bm{c}\cdot \bm{x} = 0$ of $\Aff^m$.
Then $\card{C(V_{\bm{c}})(\FF_p)} = 1+(p-1)\card{V_{\bm{c}}(\FF_p)}$.
We must show that $\card{C(V_{\bm{c}})(\FF_p)} = p^{m-2} + p^{m/2} + O(p^{(m-1)/2})$.

We count solutions to $F(\bm{x})=\bm{c}\cdot\bm{x}=0$
using the $(\bm{h},\bm{y})$ coordinates in \eqref{EQN:V_c-equations-in-h,y-coordinates}.
The locus $\bm{h}=\bm{0}$ contributes $\card{\Lambda/p\Lambda} = p^{m/2}$ solutions to \eqref{EQN:V_c-equations-in-h,y-coordinates}.
Let
\begin{equation*}
{\textstyle
U = V_{\bm{c}}\cap \set{\bm{h}\neq \bm{0}},
\quad U' = V_{\bm{c}}\cap \set{\prod h[k]\neq 0},
\quad Z_k = U\cap \set{h[k]=0},
}
\end{equation*}
and by slight abuse of notation, define the corresponding cones with origins removed:
\begin{equation*}
{\textstyle
C(U) = C(V_{\bm{c}})\cap \set{\bm{h}\neq \bm{0}},
\quad C(U') = C(V_{\bm{c}})\cap \set{\prod h[k]\neq 0},
\quad C(Z_k) = C(U)\cap \set{h[k]=0}.
}
\end{equation*}
Recall \eqref{EQN:rewrite-c(k)}.
Since $p\nmid c^\star_k$ for all $k\in [m/2]$ (by Proposition~\ref{PROP:nice-consequences-of-non-degenerate-diagonal-jet}(2)),
the equation $\bm{c}^\star\cdot\bm{h}=0$ implies that $U=U'$ if $m=4$,
and that $Z_1$, $Z_2$, $Z_3$ are pairwise disjoint if $m=6$.

Suppose first that $m=4$.
Then $U=U'$ is covered by a single affine chart, say with $h[2] = c^\star_1$ and $h[1] = -c^\star_2$.
The (remaining) defining equation $F(\bm{x})=0$ becomes
\begin{equation*}
-F_{(1)}\cdot c^\star_2\cdot y[1]^2+F_{(2)}\cdot c^\star_1\cdot y[2]^2
= 3[F_{(1)}\cdot(c^\star_2)^3-F_{(2)}\cdot(c^\star_1)^3].
\end{equation*}
Since $p\nmid c(2)^3-c(1)^3$ (by Proposition~\ref{PROP:nice-consequences-of-non-degenerate-diagonal-jet}(2)),
we get $\card{U(\FF_p)} = p - \chi(c(1)^3c(2)^3)$ by comparing $U$ with $\PP^1$.
Thus $\card{C(U)(\FF_p)} = (p-1) \card{U(\FF_p)} = p^2+O(p)$.
So $\card{C(V_{\bm{c}})(\FF_p)}
= \card{\Lambda/p\Lambda} + \card{C(U)(\FF_p)}
= p^{m-2}+p^{m/2}+O(p)$, which is better than satisfactory.

Suppose next that $m=6$.
Plugging $h[k] = 0$ into \eqref{EQN:V_c-equations-in-h,y-coordinates} identifies $C(Z_k)$ as the product of $\Aff^1$ (with coordinate $y[k]$) with a cone of the shape ``$C(U)$ for $m=4$'' (which has $p^2+O(p)$ points by the previous paragraph).
Thus $\card{C(Z_k)(\FF_p)} = p^3+O(p^2)$ for each $k$,
and $\card{\bigcup_k C(Z_k)(\FF_p)} = \sum_k \card{C(Z_k)(\FF_p)} = 3p^3+O(p^2)$ in total over $1\le k\le m/2$.
    
For $U'(\FF_p)$, first consider an individual $\bm{h}\in \FF_p^{m/2}$ with $\prod h[k]\neq 0$.
The equation $F(\bm{x})=0$ (a ternary quadratic in $\bm{y}$) has
$\mscr{N}_{\bm{h}} = p^2 + p\cdot \chi{\left(F_{(1)}h[1]\cdots F_{(3)}h[3]\cdot 3\sum F_{(k)}h[k]^3\right)}$
solutions $\bm{y}\in \FF_p^{m/2}$.
Indeed, if $\sum F_{(k)}h[k]^3=0$, then $\set{\bm{y}\in \Aff^3: F(\bm{x})=0}$ is an affine cone over a smooth conic $Q\cong \PP^1$;
otherwise, $\set{\bm{y}\in \Aff^3: F(\bm{x})=0}$ is a \emph{non-degenerate} affine quadric (and is thus the complement of a smooth conic in a smooth quadric in $\PP^3$).

Now identify $U'$ with its affine chart $h[3]=1$.
The equation $\bm{c}^\star\cdot\bm{h}=0$ becomes $h[2] = -(c^\star_2)^{-1}(c^\star_1\cdot h[1]+c^\star_3)$.
Thus $\card{U'(\FF_p)}$ equals the sum of $\mscr{N}_{\bm{h}}$ over $t\defeq h[1]\in \FF_p^\times\setminus\set{-c^\star_3/c^\star_1}$ (where we restrict $t$ so that $h[1]\cdot h[2]\neq 0$):
\begin{equation}
\label{EQN:key-deg-5-hyperelliptic-affine-chart}
-\sum_{t\in\set{0,-c^\star_3/c^\star_1}} \mscr{N}_{\bm{h}}
+\sum_{t\in\FF_p} \mscr{N}_{\bm{h}}
= -2p^2
+ p^3 + p\cdot\left(\#\set{(z,t)\in \FF_p^2: z^2 = P_{\bm{c}}(t)} - p\right),
\end{equation}
where $P_{\bm{c}}(t)
\defeq -3F_{(1)}F_{(2)}F_{(3)}(c^\star_2)^{-1}
\cdot t
\cdot (c^\star_1\cdot t+c^\star_3)
\cdot [F_{(1)}t^3-F_{(2)}(c^\star_2)^{-3}(c^\star_1\cdot t+c^\star_3)^3+F_{(3)}]$.
The count $\card{C(U')(\FF_p)}$ is $p-1$ times the right-hand side of \eqref{EQN:key-deg-5-hyperelliptic-affine-chart}.

Here $\deg_t P_{\bm{c}} = 5$,
since $p\nmid c(1)^3-c(2)^3$ by Proposition~\ref{PROP:nice-consequences-of-non-degenerate-diagonal-jet}(2).
By a routine computer calculation,
the discriminant of the quintic polynomial $P_{\bm{c}}(t)$
simplifies---up to a harmless ``unit monomial'' in $3^{\pm1}$, $F_{(k)}^{\pm1}$, $(c^\star_k)^{\pm1}$---to
$[c(1)^3-c(3)^3]^2
\cdot [c(2)^3-c(3)^3]^2
\cdot \prod[c(1)^{3/2}\pm c(2)^{3/2}\pm c(3)^{3/2}]$,
which lies in $\ZZ_p^\times$ (by Proposition~\ref{PROP:nice-consequences-of-non-degenerate-diagonal-jet}(2)).
Thus $z^2 = P_{\bm{c}}(t)$ defines an affine hyperelliptic curve over $\FF_p$ of genus $2$.
All in all, we have
(by the Weil bound for $z^2 = P_{\bm{c}}(t)$)
\begin{equation*}
\card{C(V_{\bm{c}})(\FF_p)}
= \card{\Lambda/p\Lambda} + (3p^3+O(p^2))
- 2(p^3-p^2) + (p^4-p^3) + (p^2-p)\cdot O(p^{1/2}),
\end{equation*}
which simplifies to
$\card{\Lambda/p\Lambda}+p^4+O(p^{5/2}) = p^{m-2}+p^{m/2}+O(p^{(m-1)/2})$,
as desired.
\end{proof}

\begin{remark}
By Lang--Weil for \emph{curves},
we only need $z^2 = P_{\bm{c}}(t)$ to be
\emph{absolutely irreducible} over $\FF_p$---not necessarily smooth.
However,
$p\nmid c(k)^3$ and $p\nmid c(i)^3-c(j)^3$ remain essential
throughout the proof of Lemma~\ref{LEM:point-count-bias-mod-p};
without them,
the bias could increase.
\end{remark}


\begin{lemma}
\label{LEM:bias-mod-p^l}
Let $l\geq 2$ be an integer.
Under \eqref{COND:jet-is-non-degenerate-mod-p}, \eqref{COND:p-not-divides-6F_1...F_m}, and \eqref{EQN:convenient-diagonal-change-of-variables-assumption},
we have
\begin{equation}
\label{EQN:normalized-S_c-goal-formula}
S^\natural_{\bm{c}}(p^l)
= \bm{1}_{\chi(c(1)^3)=\cdots=\chi(c(m/2)^3)}
\cdot 2^{m/2-1}\phi(p^l)p^{-l/2}.
\end{equation}
\end{lemma}

\begin{proof}
By \eqref{DEFN:normalize-S_c(n)}, the desired formula \eqref{EQN:normalized-S_c-goal-formula} is equivalent to
\begin{equation}
\label{EQN:un-normalized-S_c-goal-formula}
p^{-2l}S_{\bm{c}}(p^l)
= \bm{1}_{\chi(c(1)^3)=\cdots=\chi(c(m/2)^3)}
\cdot 2^{m/2-1}(p-1)p^{l(m-2)/2-1}.
\end{equation}

By Proposition~\ref{PROP:nice-consequences-of-non-degenerate-diagonal-jet} and \eqref{EQN:rewrite-c(k)}, we have $p\nmid c(k)^3 = (c^\star_k)^3/F_{(k)}$
for all $k\in [m/2]$.
In particular, $p\nmid \bm{c}$, so Corollary~\ref{COR:primitive-case-projective-balanced-rewrite-of-S_c} applies, since $V$ is smooth by \eqref{COND:p-not-divides-6F_1...F_m}.
Consider the map
\begin{equation}
\label{MAP:reduction-map-for-V_c-points-mod-p^l}
    \mcal{V}_{\bm{c}}(\ZZ/p^l\ZZ)\to \mcal{V}_{\bm{c}}(\ZZ/p^{l-1}\ZZ).
\end{equation}
Whenever $[\bm{x}]\in \mcal{V}_{\bm{c}}(\ZZ/p^{l-1}\ZZ)$ lies over a smooth point of $V_{\bm{c}}$, the fiber of \eqref{MAP:reduction-map-for-V_c-points-mod-p^l} over $[\bm{x}]$ has size exactly $p^{m_\ast}$ (by Hensel's lemma).
Therefore, Corollary~\ref{COR:primitive-case-projective-balanced-rewrite-of-S_c} simplifies to
\begin{equation}
\label{EQN:S_c-via-B'}
    p^{-2l}S_{\bm{c}}(p^l) = \card{\mcal{B}'(l)} - p^{m_\ast}\card{\mcal{B}'(l-1)},
\end{equation}
where $\mcal{B}'(v)$ denotes the subset of $\mcal{V}_{\bm{c}}(\ZZ/p^v\ZZ)$ lying over the singular locus of $V_{\bm{c}}$.

If $\chi(c(i)^3)\neq\chi(c(j)^3)$,
i.e.~$\chi(c^\star_i/3F_{(i)})\ne \chi(c^\star_j/3F_{(j)})$,
for some $i,j\in \mcal{K}$,
then $\mcal{B}'(1)=\emptyset$:
in fact,
there are no $\bm{x}\in\FF_p^m\setminus\set{\bm{0}}$ with
$\grad{F}(\bm{x})$, $\bm{c}$ linearly dependent over $\FF_p$.
So unless
\begin{equation}
\label{EQN:quadratic-consistency}
\chi(c(1)^3)=\cdots=\chi(c(m/2)^3)
\end{equation}
holds, we have $S_{\bm{c}}(p^l) = 0$ by \eqref{EQN:S_c-via-B'}.
So, from now on, assume \eqref{EQN:quadratic-consistency} holds.

By the conditions \eqref{EQN:quadratic-consistency} and $p\neq 2$, there exists $\lambda\in\ZZ_p^\times$ such that
$\lambda\cdot c^\star_k/F_{(k)}\in (\ZZ_p^\times)^2$ for all $k\in [m/2]$.
Say $\lambda\cdot c^\star_k = F_{(k)}d(k)^2$ for some choices $d(k)\in\ZZ_p^\times$;
write $d_i = d(k)$ when $i\in \mcal{J}(k)$.
Then $\mcal{B}'(1)$ is the set of $\FF_p$-points $[\bm{x}] = [\pm \res{d_i}]_{i\in[m]}$ with $F(\bm{x})=0$.
But by Proposition~\ref{PROP:nice-consequences-of-non-degenerate-diagonal-jet}(3),
the scheme $V_{\bm{c}}$ has exactly $2^{m/2-1}$ singular points $[\bm{x}]\in V_{\bm{c}}(\ol{\FF}_p)$,
which---by the Jacobian criterion, and the fact that $\pm\res{d_i}\in\FF_p$ for all $i\in[m]$---must
all lie in $\mcal{B}'(1)$.
Explicitly,
these $2^{m/2-1}$ points $[\bm{x}]$ arise from
the sign choices for which $\bm{x}\in \Lambda\otimes \FF_p$.

We now seek to count $\mcal{B}'(l)$ for $l\geq 2$.
Recall \eqref{EQN:V_c-equations-in-h,y-coordinates}, expressing $\mcal{V}_{\bm{c}}$ in terms of the $(\bm{h},\bm{y})$ coordinates.
Fix a point $[\bm{x}]\in \mcal{B}'(1)$,
say given in $(\bm{h},\bm{y})$ coordinates by $\bm{h}\equiv\bm{0}\bmod{p}$ and $y[k]\equiv d(k)\bmod{p}$.
Upon writing $c^\star_k = F_{(k)}d(k)^2/\lambda$,
the system \eqref{EQN:V_c-equations-in-h,y-coordinates} modulo $p^l$ becomes
\begin{equation}
\label{EQN:special-system-mod-p^l}
\sum_{1\le k\le m/2}F_{(k)} h[k]y[k]^2
\equiv_{p^l} -3\sum_{1\le k\le m/2}F_{(k)} h[k]^3,
\qquad
\sum_{1\le k\le m/2}F_{(k)} d(k)^2h[k]
\equiv_{p^l} 0.
\end{equation}
Let $\mcal{B}'_{\bm{d}}(l)$ be the set of solutions $(\bm{h},\bm{y})$ to \eqref{EQN:special-system-mod-p^l} lying over our fixed $[\bm{x}]\in \mcal{B}'(1)$.

\emph{Fix an affine chart} (i.e.~representatives in $\mcal{B}'_{\bm{d}}(l)$)
by setting $y[m/2]=d(m/2)$ identically over $\ZZ_p$.
Write $h[k] = p^sh_s[k]$ and $y[k] = d(k)+p^sy_s[k]$
(with $s=1$ for now,
but all $s\geq 1$ to be relevant below),
so $y_s[m/2] = 0$.
Then \eqref{EQN:special-system-mod-p^l} becomes
\begin{equation*}
\sum_{k\neq m/2}F_{(k)}h_s[k](2d(k)y_s[k]
+ p^sy_s[k]^2)
\equiv_{p^{l-2s}} -3p^s\sum_{k}F_{(k)}h_s[k]^3,
\;
\sum_{k}F_{(k)}d(k)^2h_s[k] \equiv_{p^{l-s}} 0.
\end{equation*}
So $\card{\mcal{B}'_{\bm{d}}(l)} = p^{m-2}\card{\mcal{A}_1(l-2)}$,
where $\mcal{A}_s(l)$ is the (non-homogeneous, affine) system
\begin{equation*}
\sum_{k\neq m/2}F_{(k)}h_s[k](2d(k)y_s[k]
+ p^sy_s[k]^2)
\equiv_{p^l} -3p^s\sum_{k}F_{(k)}h_s[k]^3,
\quad
\sum_{k}F_{(k)}d(k)^2h_s[k] \equiv_{p^l} 0.
\end{equation*}
Fix $s\geq 1$.
Clearly $\card{\mcal{A}_s(0)} = 1$,
while $\mcal{A}_s(1)$ is isomorphic to a cone over a smooth\footnote{$h_s[m/2]$ is determined by the remaining $h_s[k]$,
and $\sum_{k\neq m/2}F_{(k)}d(k)\cdot h_s[k]y_s[k] = 0$ is smooth.}
quadric in $m-2$ variables (i.e.~in $\PP^{m_\ast}$, of even dimension $m_\ast-1$)
with discriminant in $(-1)^{m/2-1}(\FF_p^\times)^2$,
so $\card{\mcal{A}_s(1)} = p^{m_\ast} + (p-1)p^{(m_\ast-1)/2}$.
For $l\geq 2$,
the origin of the cone $\mcal{A}_s(1)$
contributes $p^{m-2}\card{\mcal{A}_{s+1}(l-2)}$ points to $\mcal{A}_s(l)$,
while points away from the origin (i.e.~smooth points!)
lift uniformly to a total of $(\card{\mcal{A}_s(1)}-1)\cdot p^{(l-1)m_\ast}$ points of $\mcal{A}_s(l)$.
Thus
\begin{equation}
\label{EQN:A_s(l)-system-recursion}
\card{\mcal{A}_s(l)}
= p^{m_\ast+1}\card{\mcal{A}_{s+1}(l-2)}
+ (\card{\mcal{A}_s(1)}-1)\cdot p^{(l-1)m_\ast}
\end{equation}
for $s\geq 1$ and $l\geq 2$.
(The same holds for $l=1$,
provided we interpret $\card{\mcal{A}_s(-1)}\defeq p^{-(m_\ast+1)}$.)
By induction on $l\geq 0$
(with base cases $l=0,1$),
we immediately find that $\card{\mcal{A}_s(l)}$ is independent of
the choice of $\lambda$ and the $p$-adic square roots $d(k)$;
furthermore,
$\card{\mcal{A}_s(l)} = \card{\mcal{A}_1(l)}$ for all $s\geq 1$,
i.e.~there is no $\bm{d}$-dependence or $s$-dependence!

Finally,
by symmetry,
$\card{\mcal{B}'(l)} = 2^{m/2-1}\card{\mcal{B}'_{\bm{d}}(l)} = 2^{m/2-1}p^{m-2}\card{\mcal{A}_1(l-2)}$
for all $l\geq 1$.
(For $l=1$,
recall $\card{\mcal{A}_1(-1)}\defeq p^{-(m_\ast+1)} = p^{-(m-2)}$.)
Thus \eqref{EQN:S_c-via-B'} gives
\begin{equation*}
p^{-2l}S_{\bm{c}}(p^l)
= \card{\mcal{B}'(l)}
- p^{m_\ast}\card{\mcal{B}'(l-1)}
= 2^{m/2-1}p^{m-2}(\card{\mcal{A}_1(l-2)} - p^{m_\ast}\card{\mcal{A}_1(l-3)})
\end{equation*}
for $l\geq 2$.
To prove \eqref{EQN:un-normalized-S_c-goal-formula},
it remains to show that
\begin{equation*}
\card{\mcal{A}_1(l)}
- p^{m_\ast}\card{\mcal{A}_1(l-1)}
= (p-1)p^{l(m-2)/2-1}
= (p-1)p^{l(m_\ast+1)/2-1}
\end{equation*}
for $l\geq 0$.
To this end,
we compute $\card{\mcal{A}_1(l)} - p^{m_\ast}\card{\mcal{A}_1(l-1)}$ (using \eqref{EQN:A_s(l)-system-recursion} if $l\geq2$) to get
\begin{enumerate}
    \item $1 - p^{m_\ast}\cdot p^{-(m_\ast+1)} = 1-p^{-1}$,
    i.e.~$(p-1)p^{-1}$,
    for $l=0$;
    
    \item $[p^{m_\ast} + (p-1)p^{(m_\ast-1)/2}] - p^{m_\ast}\cdot 1 = (p-1)p^{(m_\ast-1)/2}$,
    i.e.~$(p-1)p^{(m_\ast+1)/2-1}$,
    for $l=1$;
    
    
    \item $p^{m_\ast+1}(\card{\mcal{A}_1(l-2)} - p^{m_\ast}\card{\mcal{A}_1(l-3)})$ for $l\geq 2$,
    since $p^{(l-1)m_\ast} = p^{m_\ast}\cdot p^{(l-2)m_\ast}$.
    
\end{enumerate}
By induction on $l\geq 0$,
we are done,
since $p^{m_\ast+1}\cdot p^{(l-2)(m_\ast+1)/2-1} = p^{l(m_\ast+1)/2-1}$.
\end{proof}

\begin{remark}
By induction on $l\geq 0$
(with base cases $l=0,1$),
one could prove the explicit formula
$\card{\mcal{A}_s(l)}
= p^{lm_\ast}
+ (p-1)p^{l(m_\ast+1)/2-1} (p^{l(m_\ast-1)/2}-1)/(p^{(m_\ast-1)/2}-1)$
(also valid for $l=-1$).
One could then explicitly compute $\card{\mcal{B}'_{\bm{d}}(l)}
= p^{m_\ast+1}\card{\mcal{A}_1(l-2)}$ for $l\ge 1$.
\end{remark}

For the rest of the paper, drop the assumption \eqref{EQN:convenient-diagonal-change-of-variables-assumption}.
We can finally state and prove the main result of \S\ref{SEC:proving-bias-in-linear-S_c's}.
Let $\ZZ_{(p)}^\times\defeq \set{q\in \QQ^\times: v_p(q)=0} = \ZZ_p^\times \cap \QQ$.


\begin{lemma}
\label{LEM:stating-bias-for-generic-linear-c's}
Assume $F$ is diagonal, with $m\in \set{4,6}$.
Let $\bm{c}\in \Lambda^\perp$.
Assume \eqref{COND:jet-is-non-degenerate-mod-p}.
Then
\begin{equation}
\label{EQN:key-bias-mod-p}
S^\natural_{\bm{c}}(p)
= \phi(p)p^{-1/2}
+ O(1).
\end{equation}
Also, $c(k)^3\in \ZZ_{(p)}^\times$ for all $k\in \mcal{K}$.
Finally,
\begin{equation}
\label{EQN:key-bias-mod-p^l}
S^\natural_{\bm{c}}(p^l)
= \phi(p^l)p^{-l/2}
\cdot \prod_{1\leq k\leq m/2-1} \left(1 + \chi{\left(c(k)^3 c(k+1)^3\right)}\right)
\ll \phi(p^l)p^{-l/2}
\end{equation}
for all integers $l\geq 2$.
The implied constants in \eqref{EQN:key-bias-mod-p} and \eqref{EQN:key-bias-mod-p^l} depend only on $m$.
\end{lemma}

\begin{proof}
As we noted earlier, \eqref{COND:jet-is-non-degenerate-mod-p} implies \eqref{COND:p-not-divides-6F_1...F_m}.
Now consider the unique invertible $\ZZ[1/F_1\cdots F_m]$-linear map $\bm{x}\mapsto \bm{x}'$ such that $F_ix_i^3 = F_{(k)}(x'_i)^3$ for all $k\in \mcal{K}$ and $i\in \mcal{J}(k)$.
This map transforms $F(\bm{x})$ into $F'(\bm{x}') = F'_1(x'_1)^3+\dots+F'_m(x'_m)^3$, where $F'_i = F_{(k)}$ for all $k\in \mcal{K}$ and $i\in \mcal{J}(k)$.
If we let $\bm{c}\mapsto \bm{c}'$ denote the dual linear map, then the following hold:
\begin{itemize}
    \item if $p\nmid F_1\cdots F_m$, then $S_{\bm{c}}(p^l)$, defined using $F$ as in \eqref{EQN:define-S_c(n)}, equals $S_{\bm{c}'}(p^l)$, defined using $F'$ in place of $F$;
    
    \item one can define the polynomial $(F')^\vee(\bm{c}')$ to be $F^\vee(\bm{c})$
    times a power of $F_1\cdots F_m$;
    
    \item the vector space $L'$ corresponding to $L$ is still associated to $\mcal{J}$;
    and
    
    \item we have $c_i^3/F_i = (c'_i)^3/F'_i$ for all $i\in [m]$.
\end{itemize}
By \eqref{COND:p-not-divides-6F_1...F_m}, we may thus assume \eqref{EQN:convenient-diagonal-change-of-variables-assumption}  (when proving Lemma~\ref{LEM:stating-bias-for-generic-linear-c's}).

The claim \eqref{EQN:key-bias-mod-p} now follows upon plugging Lemma~\ref{LEM:point-count-bias-mod-p} and \eqref{INEQ:E(p)-Weil-conjectures-bound} into \eqref{EQN:rewrite-S_c(p)-via-E_c}.
For the claim $c(k)^3\in\ZZ_{(p)}^\times$,
see Proposition~\ref{PROP:nice-consequences-of-non-degenerate-diagonal-jet}(1)--(2).
Finally, Lemma~\ref{LEM:bias-mod-p^l} implies \eqref{EQN:key-bias-mod-p^l} for $l\geq 2$.
\end{proof}

\begin{remark}
Since $\EE_{\bm{c}\in \Lambda^\perp/n\Lambda^\perp}[S^\natural_{\bm{c}}(n)] = \phi(n)n^{-1/2}$ (for all $n\geq 1$) by Proposition~\ref{PROP:reverse-exponential-sum-averaging},
we have chosen to formulate Lemma~\ref{LEM:stating-bias-for-generic-linear-c's}
using $\phi(n)n^{-1/2}$, not $n^{1/2}$.
\end{remark}

\section{Main delta-method analysis}
\label{SEC:full-proof-outline}

Fix $F$ in Theorem~\ref{THM:contribution-from-generically-singular-c's}.
For each $L\in \Upsilon$, recall $\Lambda$, $\Lambda^\perp$ from Definition~\ref{DEFN:primitive-sublattices}.
By Proposition~\ref{PROP:dual-linear-subvariety}, $F^\vee(\bm{c})=0$ for all $\bm{c}\in \bigcup_{L\in \Upsilon} \Lambda^\perp$.
(Recall, from Definition~\ref{DEFN:orthogonally-linear-c's}, that we call such $\bm{c}$'s \emph{linear}.)

Since $F$ is diagonal,
Proposition~\ref{PROP:characterize-diagonal-L/Q}
characterizes the linear $\bm{c}$'s via certain pairings introduced in Definition~\ref{DEFN:diagonal-permissible-pairing}.
More precisely,
the identification $\Lambda^\perp=\mcal{R}_{\mcal{J}}$ defines a bijection between
$\Upsilon$ and the set of equivalence classes of \emph{permissible pairings} $\mcal{J}$ of $[m]$.

Consider the left-hand side of \eqref{EQN:main-theorem-singular-c's}.
Recall the sets $\mcal{E}_1$, $\mcal{E}_2$ from \eqref{EQN:define-exceptional-sets-E_1,E_2}.

\begin{proposition}
\label{PROP:almost-all-solutions-are-linear}
For reals $C\geq 1$, we have $\card{\mcal{E}_1\cap [-C,C]^m}\ll_\eps C^{m/2-1+\eps}$.
\end{proposition}

\begin{proof}
This follows from the combinatorics of \cite{heath1998circle}*{p.~687}.
(In Heath-Brown's notation, any exponent $\sum_k \frac{e_k}{2}$ coming from $\mcal{E}_1$ must lie in $\set{\frac22}$ if $m=4$, and in $\set{\frac42, \frac32, \frac{2+2}{2}, \frac22}$ if $m=6$.)
\end{proof}

\begin{proposition}
\label{PROP:almost-all-linear-solutions-are-baseline}
For reals $C\geq 1$, we have $\card{\mcal{E}_2\cap [-C,C]^m}\ll_\eps C^{m/2-1+\eps}$.
\end{proposition}

\begin{proof}
Let $s\maps \QQ\to \QQ,\;q\mapsto q^2$.
Let $M\defeq \max(\abs{F_1},\dots,\abs{F_m})$.
Suppose $L\in \Upsilon$ and $\bm{c}\in \mcal{E}_2\cap \Lambda^\perp$.
Then by Observation~\ref{OBS:higher-order-diagonal-Disc-vanishing-critera}(2),
there exist distinct $k_1,k_2\in \mcal{K}$ with $c(k_1)^3 c(k_2)^3\in s(\QQ)$.
Fix $i\in \mcal{J}(k_1)$ and $j\in \mcal{J}(k_2)$;
then $(c_i^3/F_i)(c_j^3/F_j)\in s(\QQ)$.
So $(F_ic_i)(F_jc_j)\in s(\ZZ)$,
whence there exists $a\in \ZZ$ with $0<\abs{a}\leq MC$ such that $F_ic_i/a \in s(\ZZ)$ and $F_jc_j/a \in s(\ZZ)$.
Since $\card{\mcal{K}} = m/2$, it follows (upon summing over all possibilities for $\mcal{J}$, $k_1$, $k_2$, $i$, $j$) that $\mcal{E}_2\cap [-C,C]^m$ has size
$\ll_m C^{m/2-2}
\sum_{0<\abs{a}\leq MC} (MC/\abs{a})^{1/2}\cdot (MC/\abs{a})^{1/2}
\ll_M C^{m/2-1} \log(2+C)$.
\end{proof}

Propositions~\ref{PROP:almost-all-solutions-are-linear} and~\ref{PROP:almost-all-linear-solutions-are-baseline} imply that $F^\vee$ is
unsurprising (in the sense of Definition~\ref{DEFN:unsurprising-Disc-zero-locus}).
Lemma~\ref{LEM:sparse-bound} then gives
the useful bound $f(\mcal{E}_1\cup \mcal{E}_2)\ll_\eps X^{(m-1)/2+\eps}$.
(Recall $f(\mcal{S})$ from \eqref{EQN:sparse-absolute-contribution-f(S)-to-delta-method}.)

\begin{corollary}
The equality \eqref{EQN:main-theorem-singular-c's} holds,
provided that for each $L\in \Upsilon$, we have
\begin{equation}
\label{EQN-GOAL:isolate-maximal-linear-subvariety}
Y^{-2} \sum_{\bm{c}\in \Lambda^\perp\setminus \mcal{E}_2}
\sum_{n\geq 1} n^{(1-m)/2}S^\natural_{\bm{c}}(n)I_{\bm{c}}(n)
= O_\eps(X^{m/2-1/4+\eps})
+ \sigma_{\infty,L^\perp,w} X^{m/2}.
\end{equation}
\end{corollary}

\begin{proof}
Assume \eqref{EQN-GOAL:isolate-maximal-linear-subvariety} for $L\in \Upsilon$.
Since $f(\mcal{E}_2)\ll_\eps X^{(m-1)/2+\eps}$,
the relation \eqref{EQN-GOAL:isolate-maximal-linear-subvariety} implies
\begin{equation*}
Y^{-2} \sum_{\bm{c}\in \Lambda^\perp\setminus \set{\bm{0}}}
\sum_{n\geq 1} n^{(1-m)/2}S^\natural_{\bm{c}}(n)I_{\bm{c}}(n)
= O_\eps(X^{m/2-1/4+\eps})
+ \sigma_{\infty,L^\perp,w} X^{m/2}.
\end{equation*}
Upon summing over the finite set $\Upsilon$
(handling intersections using Lemma~\ref{LEM:sparse-bound}),
we obtain
\begin{equation*}
Y^{-2} \sum_{\bm{c}\in \bigcup_{L\in \Upsilon} \Lambda^\perp\setminus \set{\bm{0}}}
\sum_{n\geq 1} n^{(1-m)/2}S^\natural_{\bm{c}}(n)I_{\bm{c}}(n)
= O_\eps(X^{m/2-1/4+\eps})
+ \sum_{L\in \Upsilon} \sigma_{\infty,L^\perp,w} X^{m/2}.
\end{equation*}
By \eqref{EQN:intro-Lambda-sum-re-interpretation-formula}, \eqref{DEFN:normalize-S_c(n)}, and the bound $f(\mcal{E}_1)\ll_\eps X^{(m-1)/2+\eps}$,
the desired \eqref{EQN:main-theorem-singular-c's} follows.
\end{proof}

So \eqref{EQN-GOAL:isolate-maximal-linear-subvariety} would imply Theorem~\ref{THM:contribution-from-generically-singular-c's}.
The rest of \S\ref{SEC:full-proof-outline} is devoted to the proof of \eqref{EQN-GOAL:isolate-maximal-linear-subvariety}.
Fix $L\in \Upsilon$, and recall Proposition~\ref{PROP:characterize-diagonal-L/Q}.
We first explain why Heath-Brown's approach for $m=4$ in \cite{heath1998circle} does not seem to directly extend to $m=6$;
we then describe our approach.

Using Lemma~\ref{LEM:reverse-integral-averaging} (with $n_0=n$ and $n_1=1$) and Proposition~\ref{PROP:reverse-exponential-sum-averaging}, one can show that (in terms of certain quantities $T(\bm{j};n), J(\bm{j};n)$ we briefly discussed in \S\ref{SEC:reverse-sums-by-duality})
\begin{equation}
\label{EQN:Heath-Brown-Poisson-sum-approach}
X^{-3} \sum_{\bm{c}\in \mcal{R}_{\mcal{J}}}
\sum_{n\geq 1} n^{-m}S_{\bm{c}}(n)I_{\bm{c}}(n)
= X^{-3} \sum_{n\geq 1} n^{-m/2} \sum_{\bm{j}\in \ZZ^{m/2}} T(\bm{j};n)J(\bm{j};n);
\end{equation}
cf.~\cite{heath1998circle}*{p.~692, Poisson summation underlying Lemma~8.2}.
When $m=4$,
Heath-Brown proves that $\bm{j}=\bm{0}$ in \eqref{EQN:Heath-Brown-Poisson-sum-approach} captures the ``$\mcal{J}$-diagonal'' contribution to \eqref{EQN:define-N(X)},
and that the locus $\bm{j}\neq \bm{0}$ in \eqref{EQN:Heath-Brown-Poisson-sum-approach} forms an ``error term'' of $\ll_\eps X^{3/2+\eps}$.

When $m=4$ and $\sigma_{\infty,L^\perp,w}\neq 0$,
the $\mcal{J}$-diagonal in \eqref{EQN:define-N(X)} strictly dominates the $\bm{c}=\bm{0}$ contribution to \eqref{EQN:Heath-Brown-Poisson-sum-approach}.
When $m=6$ and $\sigma_{\infty,F,w}\cdot \sigma_{\infty,L^\perp,w}\neq 0$, however,
$\bm{c}=\bm{0}$ in \eqref{EQN:Heath-Brown-Poisson-sum-approach} is comparable in size to the $\mcal{J}$-diagonal in \eqref{EQN:define-N(X)},
so that $\bm{j}\neq \bm{0}$ in \eqref{EQN:Heath-Brown-Poisson-sum-approach} is likely \emph{no longer an error term}.
Perhaps for typical $\bm{j}\neq \bm{0}$,
the sums $T(\bm{j};n)$ can be analyzed in terms of $L$-functions,
but it is not clear where the $\bm{c}=\bm{0}$ contribution to \eqref{EQN:Heath-Brown-Poisson-sum-approach} would arise for $m=6$.
To push \eqref{EQN:Heath-Brown-Poisson-sum-approach} further---perhaps by considering small and large $n$ separately---thus seems technical and possibly delicate,
though it could be enlightening.\footnote{Heath-Brown's argument for $m=4$ is already challenging,
and the geometry involved might become even more complicated as $m$ grows (though the parity of $m/2$ may also play some role).}

Our approach to \eqref{EQN-GOAL:isolate-maximal-linear-subvariety} delays Poisson summation to the ``endgame'', thus sidestepping
\eqref{EQN:Heath-Brown-Poisson-sum-approach}.
Consider the left-hand side of \eqref{EQN-GOAL:isolate-maximal-linear-subvariety}.
We open not with Poisson summation over $\bm{c}\in \Lambda^\perp$, but with local geometry (Lemma~\ref{LEM:stating-bias-for-generic-linear-c's}).
Lemma~\ref{LEM:stating-bias-for-generic-linear-c's} exposes a uniform bias in $S_{\bm{c}}(p^l)$ over $\bm{c}\in \Lambda^\perp\setminus \mcal{E}_2$,
allowing us to decompose $S_{\bm{c}}(n)$ into simpler pieces (see \eqref{EQN:decompose-S_c-as-convolution}).
Even then, tricky issues remain (especially regarding the excised loci $\set{\bm{0}}$ and $\mcal{E}_2\setminus \set{\bm{0}}$), but Lemma~\ref{LEM:stating-bias-for-generic-linear-c's} is undoubtedly the driving force in our argument.
(However, if we could compute $I_{\bm{c}}(n)$ to greater precision when $\bm{c}\in \Lambda^\perp\setminus \mcal{E}_2$,
that might reduce our reliance on Lemma~\ref{LEM:stating-bias-for-generic-linear-c's}.)


For $\bm{c}\in \Lambda^\perp$, consider the Dirichlet series $\Phi(\bm{c},s)\defeq \sum_{n\geq 1} n^{-s}S^\natural_{\bm{c}}(n)$.
Let
\begin{equation*}
\Psi(s) \defeq
\sum_{n\geq 1} n^{-s}\phi(n)n^{-1/2}
= \zeta(s+1/2)^{-1} \cdot \zeta(s-1/2)
\end{equation*}
be the Dirichlet series for $\phi(n)n^{-1/2}$.
By Lemma~\ref{LEM:stating-bias-for-generic-linear-c's},
$\Phi(\bm{c},s)$ should typically resemble $\Psi(s)$, to leading order.
So divide $\Phi(\bm{c},s)$ by $\Psi(s)$ to define the ``error series''
\begin{equation}
\label{EQN:define-multiplicative-error-Dirichlet-series}
\sum_{n\geq 1} n^{-s}S^\natural_{\bm{c},0}(n)
\defeq \Phi(\bm{c},s)/\Psi(s)
= \zeta(s-1/2)^{-1} \cdot \zeta(s+1/2) \cdot \Phi(\bm{c},s).
\end{equation}
Since $\zeta(s)^{-1} = \sum_{n\geq 1} n^{-s}\mu(n)$ and $\zeta(s) = \sum_{n\geq 1} n^{-s}$, it follows from \eqref{EQN:define-multiplicative-error-Dirichlet-series} that
\begin{equation}
\label{EQN:expand-multiplicative-error-as-convolution}
S^\natural_{\bm{c},0}(n)
= \sum_{d_0d_1d_2=n} \mu(d_0)d_0^{1/2}\cdot d_1^{-1/2}\cdot S^\natural_{\bm{c}}(d_2).
\end{equation}

For convenience,
let $S_{\bm{c},0}(n)\defeq n^{(1+m)/2}S^\natural_{\bm{c},0}(n)$.
The multiplicativity of $S_{\bm{c}}(n)$, $\phi(n)$ in $n$ leads to Euler products for $\Phi$, $\Psi$,
and then to multiplicativity of $S_{\bm{c},0}(n)$.

We will need some basic properties of $S_{\bm{c},0}(n)$ as a function of $\bm{c}\in \Lambda^\perp$ and $n\geq 1$.


\begin{proposition}
\label{PROP:S_c,0(n)-nice-algebraic-properties}
The quantity $S_{\bm{c},0}(n)$ is a function of $n$ and $\bm{c}\bmod{n}$.
Also,
\begin{equation}
\label{EQN:S_c,0(n)-error-average}
\EE_{\bm{c}\in \Lambda^\perp/n\Lambda^\perp}[S^\natural_{\bm{c},0}(n)] = \bm{1}_{n=1}.
\end{equation}
\end{proposition}

\begin{proof}
By the first sentence of Proposition~\ref{PROP:reverse-exponential-sum-averaging},
$S_{\bm{c}}(n)$ depends at most on $n$ and $\bm{c}\bmod{n}$.
But $S_{\bm{c},0}(n)$ depends at most on the list of values $(S_{\bm{c}}(d_2))_{d_2\mid n}$,
hence at most on $n$ and $\bm{c}\bmod{n}$.
Yet by \eqref{EQN:define-multiplicative-error-Dirichlet-series}, the Dirichlet series identity $\sum_{n\geq 1} n^{-s}S^\natural_{\bm{c}}(n)
= \Psi(s) \sum_{n\geq 1} n^{-s}S^\natural_{\bm{c},0}(n)$ holds for any $\bm{c}\in \Lambda^\perp$.
Averaging formally (coefficient-wise) over $\bm{c}\in \Lambda^\perp$, we get
\begin{equation}
\label{EQN:averaged-Dirichlet-series-factorization}
\sum_{n\geq 1} n^{-s}\EE_{\bm{c}\in\Lambda^\perp/n\Lambda^\perp}[S^\natural_{\bm{c}}(n)]
= \Psi(s)
\sum_{n\geq 1} n^{-s}\EE_{\bm{c}\in\Lambda^\perp/n\Lambda^\perp}[S^\natural_{\bm{c},0}(n)].
\end{equation}
But by the final sentence of Proposition~\ref{PROP:reverse-exponential-sum-averaging},
the left-hand side of \eqref{EQN:averaged-Dirichlet-series-factorization} equals $\Psi(s)$.
So \eqref{EQN:S_c,0(n)-error-average} follows
formally by division.
\end{proof}


We now provide some bounds on $S_{\bm{c},0}(n)$ for $\bm{c}\in \Lambda^\perp$.
By \eqref{EQN:expand-multiplicative-error-as-convolution}, we have
\begin{equation}
\label{EQN:divide-S_c,0-by-sqrt-n}
n^{-1/2} S^\natural_{\bm{c},0}(n)
= \sum_{d_0d_1d_2=n} \mu(d_0)\cdot d_1^{-1}\cdot d_2^{-1/2}S^\natural_{\bm{c}}(d_2).
\end{equation}
Let $\tau_3(n)\defeq \sum_{d_0d_1d_2=n} 1$.
For any $n\geq 1$,
the triangle inequality on \eqref{EQN:divide-S_c,0-by-sqrt-n}, followed by an application of \eqref{INEQ:cube-free-bound} to cube-free divisors of $n$, yields
\begin{equation}
\label{INEQ:(E1)-general}
\abs{S^\natural_{\bm{c},0}(n)}
\leq \tau_3(n) n^{1/2}
\cdot \max_{d\mid n} d^{-1/2}\abs{S^\natural_{\bm{c}}(d)}
\ll_\eps \tau_3(n) n^{1/2+\eps}
\cdot \max_{d\in \mcal{N}_{\ge}(3):\, d\mid n} d^{-1/2}\abs{S^\natural_{\bm{c}}(d)}.
\end{equation}
If $p\nmid j^{2^{m/2-1}}{F^\vee}(\bm{c})$ and $l\geq 2$, then \eqref{INEQ:(E1)-general} and \eqref{EQN:key-bias-mod-p^l} imply
\begin{equation}
\label{INEQ:(E2)-good-prime-power}
    \abs{S^\natural_{\bm{c},0}(p^l)}\ll \tau_3(p^l) p^{l/2}.
\end{equation}
If $p\nmid j^{2^{m/2-1}}{F^\vee}(\bm{c})$, then \eqref{EQN:expand-multiplicative-error-as-convolution} and \eqref{EQN:key-bias-mod-p} imply
\begin{equation}
\label{INEQ:(E3)-good-prime}
    S^\natural_{\bm{c},0}(p) = -p^{1/2} + p^{-1/2} + S^\natural_{\bm{c}}(p) = O(1).
\end{equation}


The bounds \eqref{INEQ:(E1)-general}, \eqref{INEQ:(E2)-good-prime-power}, \eqref{INEQ:(E3)-good-prime} are most useful in conjunction with multiplicativity.
Let
\begin{equation*}
\mcal{N}^{\bm{c}} \defeq \set{n\geq 1: p\mid n\Rightarrow p\nmid j^{2^{m/2-1}}{F^\vee}(\bm{c})},
\quad
\mcal{N}_{\bm{c}} \defeq \set{n\geq 1: p\mid n\Rightarrow p\mid j^{2^{m/2-1}}{F^\vee}(\bm{c})}
\end{equation*}
for each $\bm{c}\in \Lambda^\perp$.
For each integer $t\geq 1$, let
\begin{equation*}
\mcal{N}_{\le}(t) \defeq \set{n\geq 1: p\mid n\Rightarrow v_p(n)\le t},
\quad
\mcal{N}_{\ge}(t) \defeq \set{n\geq 1: p\mid n\Rightarrow v_p(n)\ge t}.
\end{equation*}

\begin{lemma}
[Cf.~\cite{wang2023_large_sieve_diagonal_cubic_forms}*{Lemma~3.4}]
\label{LEM:count-square-full-and-cube-full}
Let $N, t\geq 1$ be integers.
Then $\card{\set{N\le n<2N: n\in \mcal{N}_{\ge}(t)}}\ll_t N^{1/t}$.
Also, if $\bm{c}\in \Lambda^\perp\setminus \mcal{E}_2$, then
$\card{\set{N\le n<2N: n\in \mcal{N}_{\bm{c}}}}\ll_\eps N^\eps \norm{\bm{c}}^\eps$.
\end{lemma}

\begin{proof}
The first bound is familiar.
It remains to prove the second.
Say $\bm{c}\in \Lambda^\perp\setminus \mcal{E}_2$.
Then the vector $j^{2^{m/2-1}}{F^\vee}(\bm{c})$ is nonzero, and thus has a nonzero coordinate $R$.
But $\#\set{N\le n<2N: n\mid R^\infty}\ll_\eps (RN)^\eps$ (by Rankin's trick), and $R\ll \norm{\bm{c}}^{\deg F^\vee}$.
\end{proof}

\begin{lemma}
\label{LEM:bound-S_c,0-average-in-terms-of-S}
Let $\bm{c}\in \Lambda^\perp\setminus \mcal{E}_2$ and $N\in \set{1,2,4,8,\ldots}$.
Suppose $\norm{\bm{c}}\leq X^{10}$.
Then
\begin{equation*}
\sum_{N\leq n<2N} \abs{S^\natural_{\bm{c},0}(n)}
\ll_\eps (XN)^\eps \cdot N \sum_{d<2N:\, d\in \mcal{N}_{\ge}(3)\cap \mcal{N}_{\bm{c}}} d^{-1} \abs{S^\natural_{\bm{c}}(d)}.
\end{equation*}
\end{lemma}

\begin{proof}
Any integer $n\ge 1$ can be written (uniquely) as $n_1n_2n_3$, where $n_1$, $n_2$, $n_3$ are pairwise coprime integers satisfying
$n_1\in \mcal{N}^{\bm{c}}\cap \mcal{N}_{\le}(1)$,
$n_2\in \mcal{N}^{\bm{c}}\cap \mcal{N}_{\ge}(2)$,
$n_3\in \mcal{N}_{\bm{c}}$.
Upon writing $S_{\bm{c},0}(n) = \prod_{1\le i\le 3} S_{\bm{c},0}(n_i)$, and applying \eqref{INEQ:(E3)-good-prime} to primes $p\mid n_1$, \eqref{INEQ:(E2)-good-prime-power} to primes $p\mid n_2$, and \eqref{INEQ:(E1)-general} to $n_3$, we get (by dyadic summation over $n_1$, $n_2$, $n_3$)
\begin{equation*}
\sum_{N\leq n<2N} \abs{S^\natural_{\bm{c},0}(n)}
\ll_\eps \sum_{\substack{N_1,N_2,N_3\mid 2N: \\ N/4<N_1N_2N_3<2N}} N^\eps
\sum_{\substack{n_1,n_2\ge 1: \\ N_1\le n_1<2N_1, \\
N_2\le n_2<2N_2,\; n_2\in \mcal{N}_{\ge}(2)}} n_2^{1/2}
\sum_{\substack{n_3,d\ge 1: \\ N_3\le n_3<2N_3,\; n_3\in \mcal{N}_{\bm{c}}, \\
d\in \mcal{N}_{\ge}(3),\; d\mid n_3}} (n_3/d)^{1/2} \abs{S^\natural_{\bm{c}}(d)}.
\end{equation*}
Upon summing over $n_1$, $n_2$, $n_3$ (for each fixed $d$), we get, by Lemma~\ref{LEM:count-square-full-and-cube-full},
\begin{equation*}
\sum_{N\leq n<2N} \abs{S^\natural_{\bm{c},0}(n)}
\ll_\eps \sum_{\substack{N_1,N_2,N_3\mid 2N: \\ N_1N_2N_3<2N}} N^\eps N_1N_2
\sum_{\substack{d<2N_3: \\ d\in \mcal{N}_{\bm{c}}\cap \mcal{N}_{\ge}(3)}} (N_3\norm{\bm{c}})^\eps (N_3/d)^{1/2} \abs{S^\natural_{\bm{c}}(d)}.
\end{equation*}
But $(N_3/d)^{1/2} \ll N_3/d$ for $d<2N_3$.
And $\norm{\bm{c}}\leq X^{10}$.
So Lemma~\ref{LEM:bound-S_c,0-average-in-terms-of-S} follows.
\end{proof}

Given $n$ in \eqref{EQN-GOAL:isolate-maximal-linear-subvariety},
we may use \eqref{EQN:define-multiplicative-error-Dirichlet-series} to decompose $S^\natural_{\bm{c}}(n)$ as a Dirichlet convolution:
\begin{equation}
\label{EQN:decompose-S_c-as-convolution}
S^\natural_{\bm{c}}(n) = \sum_{n_0n_1 = n} S^\natural_{\bm{c},0}(n_0) \cdot \phi(n_1)n_1^{-1/2}.
\end{equation}
We will study the ranges $n_1\geq Y/P$ and $n_1<Y/P$ separately, for a parameter $P$ to be chosen in \eqref{EQN:choice-of-P=P_X}.
We first handle the range $n_1<Y/P$, using Lemma~\ref{LEM:bound-S_c,0-average-in-terms-of-S}.
It turns out we will not need the full strength of Lemma~\ref{LEM:bound-S_c,0-average-in-terms-of-S}
(which might however still be useful in the future).

\begin{definition}
\label{DEFN:Sigma_K,X,T}
For a real number $K>0$ and a set $\mcal{T}\belongs \ZZ^m$, let
\begin{equation*}
\Sigma_{<K}(X, \mcal{T})
\defeq Y^{-2} \sum_{\bm{c}\in \mcal{T}}
\sum_{n\geq 1}
n^{(1-m)/2}
I_{\bm{c}}(n)
\sum_{\substack{n_0n_1 = n: \\ n_1<K}} S^\natural_{\bm{c},0}(n_0)
\cdot \phi(n_1)n_1^{-1/2}.
\end{equation*}
Similarly define $\Sigma_{\ge K}(X, \mcal{T})$ (by replacing $n_1<K$ with $n_1\ge K$).
\end{definition}

Using \eqref{EQN:decompose-S_c-as-convolution}, one may rewrite the left-hand side of \eqref{EQN-GOAL:isolate-maximal-linear-subvariety} as
\begin{equation}
    \label{EXPR:decompose-LHS-of-goal-into-3-pieces}
    \Sigma_{<Y/P}(X, \Lambda^\perp\setminus \mcal{E}_2)
    + \Sigma_{\ge Y/P}(X, \Lambda^\perp)
    - \Sigma_{\ge Y/P}(X, \mcal{E}_2).
\end{equation}

\begin{lemma}
\label{LEM:handle-n_1<Y/P}
Uniformly over reals $X, P\geq 1$, we have
\begin{equation}
\label{INEQ:result-of-absolute-bound-for-n_1<Y/P}
\Sigma_{<Y/P}(X, \Lambda^\perp\setminus \mcal{E}_2)
\ll_\eps X^{m/2+\eps} P^{-1/2}.
\end{equation}
\end{lemma}

\begin{proof}
We proceed somewhat crudely.
By Proposition~\ref{PROP:standard-delta-method-n,c-cutoffs},
the absolute value of the left-hand side of \eqref{INEQ:result-of-absolute-bound-for-n_1<Y/P}
is at most $O_{\eps,A}(X^{-A})$ plus the quantity
\begin{equation}
\label{EXPR:crude-absolute-bound-quantity}
Y^{-2} \sum_{\substack{\bm{c}\in \Lambda^\perp\setminus \mcal{E}_2: \\ \norm{\bm{c}}\le X^{1/2+\eps}}}\,
\sum_{\substack{n_0n_1<M_2Y: \\ n_1<Y/P}}
(n_0n_1)^{(1-m)/2}\abs{I_{\bm{c}}(n_0n_1)}
\cdot n_1^{1/2}\abs{S^\natural_{\bm{c},0}(n_0)}.
\end{equation}

We now examine an individual $\bm{c}\in \Lambda^\perp\setminus \mcal{E}_2$.
By Observation~\ref{OBS:higher-order-diagonal-Disc-vanishing-critera}(1), $c_1\cdots c_m\neq 0$.
So
\begin{equation}
\label{INEQ:Heath-Brown-integral-bound}
I_{\bm{c}}(n)
\ll_\eps X^{m+\eps} (X\norm{\bm{c}}/n) \prod_{1\le i\le m} (X\abs{c_i}/n)^{-1/2}
\end{equation}
by \cite{heath1998circle}*{Lemma~3.2}.
Upon inserting \eqref{INEQ:Heath-Brown-integral-bound} into \eqref{EXPR:crude-absolute-bound-quantity}, dyadically decomposing $n_0$, and then applying Lemma~\ref{LEM:bound-S_c,0-average-in-terms-of-S}, we find that the quantity \eqref{EXPR:crude-absolute-bound-quantity} is
\begin{equation*}
\ll_\eps \sum_{\substack{\bm{c}\in \Lambda^\perp\setminus \mcal{E}_2: \\ \norm{\bm{c}}\le X^{1/2+\eps}}}\,
\sum_{\substack{N_0n_1<M_2Y: \\ N_0\in \set{1,2,4,8,\ldots},\; n_1<Y/P}}
\frac{Y^{-2} X^{m+\eps} X^{1-m/2}\norm{\bm{c}}}{(N_0n_1)^{1/2} \abs{c_1\cdots c_m}^{1/2}}
\cdot n_1^{1/2}\cdot N_0 \sum_{\substack{d<2N_0: \\ d\in \mcal{N}_{\ge}(3)}} d^{-1} \abs{S^\natural_{\bm{c}}(d)}.
\end{equation*}
But \eqref{INEQ:crude-pointwise-bound-on-diagonal-S_c(n)} yields
$S^\natural_{\bm{c}}(d) \ll_\eps d^{1/2+\eps} \prod_{1\le i\le m} \map{sq}(c_i)^{1/4}$
(when $c_1\cdots c_m\neq 0$).
Plugging this in, and noting that $\sum_{d\in \mcal{N}_{\ge}(3)} d^{-1/2} \ll 1$ by Lemma~\ref{LEM:count-square-full-and-cube-full}, we find that \eqref{EXPR:crude-absolute-bound-quantity} is
\begin{equation*}
\ll_\eps \sum_{\bm{c}\in \Lambda^\perp\setminus \mcal{E}_2:\, \norm{\bm{c}}\le X^{1/2+\eps}}\,
\sum_{n_1<Y/P}
Y^{-2} X^{1+m/2+\eps} \norm{\bm{c}}
\cdot (M_2Y/n_1)^{1/2} \prod_{k\in \mcal{K}} \frac{\map{sq}(c_{\min\mcal{J}(k)})^{1/2}}{\abs{c_{\min\mcal{J}(k)}}},
\end{equation*}
in the notation of Definition~\ref{DEFN:diagonal-permissible-pairing}.
But $\sum_{0<\abs{c}\le C} \map{sq}(c)^{1/2}/\abs{c}\ll_\eps C^\eps$ by Lemma~\ref{LEM:count-square-full-and-cube-full} (or a calculation with Euler products).
So \eqref{EXPR:crude-absolute-bound-quantity} is
\begin{equation*}
\ll_\eps \sum_{n_1<Y/P} Y^{-2} X^{1+m/2+\eps} X^{1/2+\eps} \cdot (Y/n_1)^{1/2}
\ll Y^{-2} X^{3/2+m/2+2\eps} Y^{1/2} (Y/P)^{1/2}.
\end{equation*}
Plugging in $Y=X^{3/2}$ leads to \eqref{INEQ:result-of-absolute-bound-for-n_1<Y/P}.
\end{proof}

\begin{remark}
We have used diagonality of $F$.
In general (leaving details to the interested reader),
one can prove \eqref{INEQ:result-of-absolute-bound-for-n_1<Y/P} for $m\le 11$ under the axioms \ref{RMK:non-diagonal-F}(2)--(3),
by replacing \eqref{INEQ:Heath-Brown-integral-bound} with \eqref{INEQ:Hessian-free-integral-bound}, then using \eqref{COND:linear-locus-degeneracy-contains-W=0}, \eqref{COND:second-moment-axiom-over-nonzero-W} instead of \eqref{INEQ:crude-pointwise-bound-on-diagonal-S_c(n)}, and then verifying the inequality
\begin{equation*}
\frac{Y^{-2} X^{m+\eps} N_0^{1/2} N_1 P^{1/2}}{(XC+N_0N_1)^{m/2-1}}
\sum_{d<2N_0:\, d\in \mcal{N}_{\ge}(3)} \frac{C^{m/4} (C^{m/4}+d^{m/12})}{d^{1/2}}
\ll_\eps X^{m/2+O(\eps)}
\end{equation*}
over reals $C, P, N_0, N_1\ge 1$ with $C\le X^{1/2+\eps}$ and $N_0, P\le 100(1+M_2)Y/N_1$.
\end{remark}

In terms of $M_1$ from Lemma~\ref{LEM:reverse-integral-averaging}, let
\begin{equation}
\label{EQN:choice-of-P=P_X}
P=P_X\defeq M_1^{-1} X^{1/2}.
\end{equation}
We now turn to the range $n_1\ge Y/P$ in \eqref{EXPR:decompose-LHS-of-goal-into-3-pieces}.
Recall $\Sigma_{\ge K}(X, \mcal{T})$ from Definition~\ref{DEFN:Sigma_K,X,T}.

\begin{lemma}
\label{LEM:handle-n_1>=Y/P-after-over-extending-c}
Uniformly over reals $X\geq 1$, we have
\begin{equation}
\label{EQN:over-extended-asymptotic-evaluation-for-n_1>=Y/P}
\Sigma_{\ge Y/P}(X, \Lambda^\perp)
= \sigma_{\infty,L^\perp,w} X^{m/2} (1 + O(P^{-1})).
\end{equation}
\end{lemma}

\begin{proof}
Suppose $n,n_0,n_1\ge 1$ are integers with $n_1\ge Y/P$.
By the first part of Proposition~\ref{PROP:S_c,0(n)-nice-algebraic-properties},
$S_{\bm{c},0}(n_0)$ only depends on $\bm{c}\bmod{n_0}$.
Because $n_1\ge Y/P = M_1X$, Lemma~\ref{LEM:reverse-integral-averaging} therefore implies
\begin{equation*}
\sum_{\bm{c}\in\Lambda^\perp}
S^\natural_{\bm{c},0}(n_0)
\cdot I_{\bm{c}}(n)
= \sum_{\bm{b}\in\Lambda^\perp/n_0\Lambda^\perp}
S^\natural_{\bm{b},0}(n_0)
\cdot n_1^{m/2}
\cdot \sigma_{\infty,L^\perp,w}X^{m/2}h(n/Y,0).
\end{equation*}
By \eqref{EQN:S_c,0(n)-error-average} and the equality $\card{\Lambda^\perp/n_0\Lambda^\perp} = n_0^{m/2}$,
we conclude that
\begin{equation*}
\sum_{\bm{c}\in\Lambda^\perp}
S^\natural_{\bm{c},0}(n_0)
\cdot I_{\bm{c}}(n)
= n^{m/2}
\cdot\sigma_{\infty,L^\perp,w}X^{m/2}h(n/Y,0)
\cdot\bm{1}_{n_0=1}.
\end{equation*}
The left-hand side of \eqref{EQN:over-extended-asymptotic-evaluation-for-n_1>=Y/P} (see Definition~\ref{DEFN:Sigma_K,X,T}) thus simplifies to
\begin{equation*}
Y^{-2} \sum_{n\ge Y/P} \sigma_{\infty,L^\perp,w}X^{m/2}h(n/Y,0)\cdot \phi(n),
\end{equation*}
which equals $\sigma_{\infty,L^\perp,w} X^{m/2} (1 + O(P^{-1}))$ by Proposition~\ref{PROP:phi(n)h(n/Y,0)-sum} (below).
\end{proof}

\begin{proposition}
\label{PROP:phi(n)h(n/Y,0)-sum}
$\sum_{n\geq Y/P} \phi(n) h(n/Y,0)
= Y^2(1 + O(P^{-1}))$.
\end{proposition}

\begin{proof}
By \cite{heath1998circle}*{final paragraph on p.~692, and second paragraph on p.~676}, we have (in terms of the function $\omega$ defined on \cite{heath1998circle}*{p.~676})
\begin{equation*}
\sum_{n\geq1}\phi(n)h(n/Y,0)
= Y\sum_{n\geq 1} \omega(n/Y)
= Y^2(1 + O_A(Y^{-A})).
\end{equation*}
Furthermore,
$h(x,0)\ll x^{-1}$ by \cite{heath1996new}*{Lemma~4},
so
\begin{equation*}
\sum_{n<Y/P}\phi(n)h(n/Y,0)
\ll \sum_{n<Y/P} \phi(n)\cdot Y/n
\leq \sum_{n<Y/P} Y
\ll Y^2/P.
\end{equation*}
Proposition~\ref{PROP:phi(n)h(n/Y,0)-sum} follows upon writing $\sum_{n\ge Y/P} = \sum_{n\ge 1} - \sum_{n<Y/P}$.
\end{proof}

\begin{lemma}
\label{LEM:handle-n_1>=Y/P-for-extra-c's}
Uniformly over reals $X\geq 1$, we have
\begin{equation}
\label{EQN:loose-end-for-n_1>=Y/P}
\Sigma_{\ge Y/P}(X, \mcal{E}_2)
\ll_\eps X^{(m-1)/2+\eps}.
\end{equation}
\end{lemma}

\begin{proof}
The main subtlety here is that we must treat $\bm{c}\neq \bm{0}$ and $\bm{c}=\bm{0}$ separately.

First, given $n$,
the bound \eqref{INEQ:(E1)-general}, applied directly to $S_{\bm{c},0}(n_0)$ for each $n_0\mid n$,
implies
\begin{equation}
\label{INEQ:loose-end-reduce-S_c,0-to-S}
\sum_{\substack{n_0n_1 = n: \\ n_1\ge Y/P}}\,
\abs{S^\natural_{\bm{c},0}(n_0)}\cdot \phi(n_1)n_1^{-1/2}
\ll_\eps \bm{1}_{n\ge Y/P} \cdot n^{1/2+\eps} \cdot \max_{d\mid n} d^{-1/2}\abs{S^\natural_{\bm{c}}(d)}.
\end{equation}
Inserting \eqref{INEQ:loose-end-reduce-S_c,0-to-S} into $\Sigma_{\ge Y/P}(X, \mcal{E}_2\setminus \set{\bm{0}})$ (see Definition~\ref{DEFN:Sigma_K,X,T}), and recalling \eqref{EQN:sparse-absolute-contribution-f(S)-to-delta-method}, we get
\begin{equation*}
\Sigma_{\ge Y/P}(X, \mcal{E}_2\setminus \set{\bm{0}})
\ll_\eps Y^{-2}\sum_{\bm{c}\in \mcal{E}_2\setminus \set{\bm{0}}}
\sum_{n\geq 1} n^{1-m/2+\eps} \abs{I_{\bm{c}}(n)}\cdot \max_{d\mid n} d^{-1/2}\abs{S^\natural_{\bm{c}}(d)}
\ll_\eps Y^\eps f(\mcal{E}_2),
\end{equation*}
where $f(\mcal{E}_2)\ll_\eps X^{(m-1)/2+\eps}$ by Proposition~\ref{PROP:almost-all-linear-solutions-are-baseline} and Lemma~\ref{LEM:sparse-bound}.
Similarly, \eqref{INEQ:loose-end-reduce-S_c,0-to-S} gives
\begin{equation*}
\Sigma_{\ge Y/P}(X, \set{\bm{0}})
\ll_\eps Y^{-2} \sum_{n\geq Y/P}
n^{1-m/2+\eps} \abs{I_{\bm{0}}(n)}\cdot \max_{d\mid n} d^{-1/2}\abs{S^\natural_{\bm{0}}(d)},
\end{equation*}
which is $\ll_\eps Y^{-2} X^{m+\eps} (Y/P)^{(4-m)/3} \asymp X^{m-3+(4-m)/3+\eps}$ by \eqref{INEQ:I_0(n)-bounded} and Lemma~\ref{LEM:sum-S_0(n)-trivially} (summed over $N\in \set{1,2,4,8,\ldots}$).
But $m-3+(4-m)/3 \le (m-1)/2$, since $4\le m\le 7$.
\end{proof}

\begin{remark}
We have used diagonality of $F$.
In general, under \ref{RMK:non-diagonal-F}(1)--(3), one can prove $\Sigma_{\ge Y/P}(X, \mcal{E}_2) \ll_\eps X^{m/2-1/4+\eps}$ (by using Lemma~\ref{LEM:axiomatic-bound-on-sparse-contribution-f(S)} in place of Lemma~\ref{LEM:sparse-bound}).
\end{remark}


By \eqref{INEQ:result-of-absolute-bound-for-n_1<Y/P}, \eqref{EQN:choice-of-P=P_X}, \eqref{EQN:over-extended-asymptotic-evaluation-for-n_1>=Y/P}, and \eqref{EQN:loose-end-for-n_1>=Y/P},
the quantity \eqref{EXPR:decompose-LHS-of-goal-into-3-pieces} simplifies to
\begin{equation*}
    O_\eps(X^{m/2-1/4+\eps}) + \sigma_{\infty,L^\perp,w} X^{m/2},
\end{equation*}
matching
the right-hand side of \eqref{EQN-GOAL:isolate-maximal-linear-subvariety}.
So \eqref{EQN-GOAL:isolate-maximal-linear-subvariety} holds, thus concluding \S\ref{SEC:full-proof-outline}.

\section*{Acknowledgements}

This paper is an important component of the thesis work described in \cite{wang2022thesis};
many of my acknowledgements there
apply here as well.
I also thank my advisor, Peter Sarnak, for many helpful suggestions and questions on
the exposition, references, assumptions, and scope
of (various drafts of) the present work.{\let\thefootnote\relax\footnote{This work was partially supported by NSF grant DMS-1802211, and the European Union's Horizon~2020 research and innovation program under the Marie Sk\l{}odowska-Curie Grant Agreement No.~101034413.}}
I am also grateful to Trevor Wooley for providing some helpful general comments on special subvarieties and the reference \cite{vaughan1995certain}.
I thank Tim Browning for inspiring part of the current title of the paper.
Finally, thanks are due to the referee for providing many helpful suggestions.





\bibliographystyle{amsxport}
\bibliography{master.bib}
\end{document}